\documentclass[10pt,twocolumn,twoside]{IEEEtran}

\newtheorem{defin}{Definition}

\newtheorem{prop}{Proposition}
\newtheorem{lemma}{Lemma}

\newtheorem{example}{Example}

\newtheorem{assumption}{Assumption}

\newcommand{\be}{\begin{equation}}
\newcommand{\ee}{\end{equation}}
\newcommand{\ba}{\begin{array}}
\newcommand{\ea}{\end{array}}
\newcommand{\bea}{\begin{eqnarray}}
\newcommand{\eea}{\end{eqnarray}}

\newcommand{\diag}{{\mbox{diag}}}

\newcommand{\tran}{^{\mbox{\scriptsize T}}}  

\newcommand{\asign}{{\mbox{$\colon\hspace{-2mm}=\hspace{1mm}$}}}
\newcommand{\ssum}[1]{\mathop{ \textstyle{\sum}}_{#1}}

\newcommand{\vbar}{\raisebox{.17ex}{\rule{.04em}{1.35ex}}}
\newcommand{\vbarind}{\raisebox{.01ex}{\rule{.04em}{1.1ex}}}
\newcommand{\D}{\ifmmode {\rm I}\hspace{-.2em}{\rm D} \else ${\rm I}\hspace{-.2em}{\rm D}$ \fi}
\newcommand{\T}{\ifmmode {\rm I}\hspace{-.2em}{\rm T} \else ${\rm I}\hspace{-.2em}{\rm T}$ \fi}
\newcommand{\B}{\ifmmode {\rm I}\hspace{-.2em}{\rm B} \else \mbox{${\rm I}\hspace{-.2em}{\rm B}$} \fi}
\newcommand{\Hil}{\ifmmode {\rm I}\hspace{-.2em}{\rm H} \else \mbox{${\rm I}\hspace{-.2em}{\rm H}$} \fi}
\newcommand{\C}{\ifmmode \hspace{.2em}\vbar\hspace{-.31em}{\rm C} \else \mbox{$\hspace{.2em}\vbar\hspace{-.31em}{\rm C}$} \fi}
\newcommand{\Cind}{\ifmmode \hspace{.2em}\vbarind\hspace{-.25em}{\rm C} \else \mbox{$\hspace{.2em}\vbarind\hspace{-.25em}{\rm C}$} \fi}
\newcommand{\Q}{\ifmmode \hspace{.2em}\vbar\hspace{-.31em}{\rm Q} \else \mbox{$\hspace{.2em}\vbar\hspace{-.31em}{\rm Q}$} \fi}
\newcommand{\Z}{\ifmmode {\rm Z}\hspace{-.28em}{\rm Z} \else ${\rm Z}\hspace{-.38em}{\rm Z}$ \fi}

\renewcommand{\Re}{\mbox {Re}}
\renewcommand{\Im}{\mbox {Im}}

\renewcommand{\vec}[1]{{\bf{#1}}}     

\usepackage[dvips]{graphicx}
\usepackage{amsmath}
\usepackage{caption}
\usepackage{subcaption}
\usepackage{cite}
\usepackage{url}
\usepackage{rotating}
\usepackage{color}
\usepackage{xcolor}
\usepackage{changepage}
\usepackage{amssymb}
\usepackage{psfrag}
\usepackage{caption}
\usepackage{booktabs}
\usepackage{array}
\usepackage{fixltx2e}
\usepackage{epstopdf}
\usepackage{epsfig}
\usepackage{tikz}
\usetikzlibrary{arrows,positioning}
\usepackage{cases}
\usepackage{multirow}

\newcommand{\BlueText}[1]{#1 }

\newcommand{\sr}[1]{^{\mbox{\tiny {#1}}}}  
\newcommand{\R}{\ifmmode {\rm I}\hspace{-.2em}{\rm R} \else ${\rm I}\hspace{-.2em}{\rm R}$ \fi} 
\newcommand{\bec}[1]{\bar{\vec{#1}}}  

\newcommand{\captionFontSize}{\footnotesize}
\newcommand{\subCaptionFontSize}{\footnotesize}

\AtBeginDocument{\setlength\abovedisplayskip{0.8mm}}
\AtBeginDocument{\setlength\belowdisplayskip{0.8mm}}
\newcommand\subparagraph{%
  \@startsection{subparagraph}{5}
  {\parindent}
  {3.25ex \@plus 1ex \@minus .2ex}
  {-1em}
  {\normalfont\normalsize\bfseries}}
\usepackage[compact]{titlesec}
\titlespacing*{\section}{0pt}{5pt}{1pt}

\begin{document} 
\title{A Distributed Approach for the Optimal Power Flow Problem Based on ADMM and Sequential Convex Approximations }
\author{S. Magn\'{u}sson,
        P. C. Weeraddana,~\IEEEmembership{Member,~IEEE}, and
        C.~Fischione,~\IEEEmembership{Member,~IEEE}
\thanks{The authors are with Electrical Engineering School, Access Linnaeus Center, KTH Royal Institute of Technology, Stockholm, Sweden. E-mail:\textit{ \{sindrim, chatw, carlofi\}@kth.se}.

The authors would like to thank Anders Forsgren for valuable comments.}}


\maketitle

\begin{abstract}
The optimal power flow (OPF) problem, which plays a central role in operating electrical networks is considered. The problem is nonconvex and is in fact NP hard. Therefore, designing efficient algorithms of practical relevance is crucial, though their global optimality is not guaranteed. \BlueText{Existing semi-definite programming relaxation based approaches are restricted to OPF problems where zero duality holds. In this paper, an efficient novel method to address the \emph{general nonconvex} OPF problem is investigated.}
The proposed method is based on alternating direction method of multipliers combined with sequential convex approximations. The global OPF problem is decomposed into smaller problems associated to each bus of the network, the solutions of which are coordinated via a light communication protocol. Therefore, the proposed method is highly scalable. The convergence properties of the proposed algorithm are mathematically substantiated. Finally, the proposed algorithm is evaluated on a number of test examples, where the convergence properties of the proposed algorithm are numerically substantiated and the performance is compared with a global optimal method.
\end{abstract}

\vspace{-0mm}
\begin{IEEEkeywords}\vspace{-0mm}
Optimal power flow, distributed optimization, smart grid.
\end{IEEEkeywords}


\section{Introduction} \label{sec:intro}

 The optimal power flow (OPF) problem in electrical networks determines optimally the amount of power to be generated at each generator. 
 Moreover, it decides how to dispatch the power such that a global network-wide objective criterion is optimized,
 while ensuring that the power demand of each consumer is met and that the related laws of physics are held.
 Traditionally, the OPF problem has only been solved in transmission networks.
However, the extensive information gathering of individual power consumption in the smart grid has made the problem relevant, not only in transmission networks, but also in distribution networks which deliver electricity to end users.

\subsection{Previous work}
\BlueText{The problem was originally presented by Carpentier \cite{opf_original_carpentier}, and has been extensively studied since then and become of great  importance in efficient operation of power systems \cite{Frank-Steponavice-Rebennack-OPF-bibligraphic-survey-I}.
 The OPF problem is nonconvex due to quadratic relationship between the powers and the voltages and because of a lower bound on the voltage magnitudes.
 In fact, the problem is NP-hard, see~\cite{Lavaei-2012-OPF-TPS}.
 Therefore, practical and general purpose algorithms must rely on some approximations or heuristics.
  We refer the reader to~\cite{Frank-Steponavice-Rebennack-OPF-bibligraphic-survey-I,Frank-Steponavice-Rebennack-OPF-bibligraphic-survey-II} for a contemporary survey of~OPF.}


It is well known that the OPF problem is equivalently reformulated as a rank constrained problem~\cite{Bai_first_SDPrelaxation_of_the_the_OPF}. As a result, classic convex approximation techniques are applied to handle nonconvexities of the rank constraint, which usually results in a semidefinite program~\cite{convex_boyd} (SDP). SDP relaxations to OPF have gained a lot of attention recently, see~\cite{Lavaei-2012-OPF-TPS,Bose-2011-OPF-Allerton,2011arXiv1107.1467Z,2012arXiv1204.4419L} and references therein. Authors in \cite{Lavaei-2012-OPF-TPS} show that SDP relaxation is equivalent to the dual problem of the original OPF. Moreover, sufficient conditions for zero duality and mechanisms to recover the primal solution by the dual problem are given. Thus, \cite{Lavaei-2012-OPF-TPS} identifies a class of OPF problems, where the global optimum is obtained efficiently by using convex tools.
 Some other classes of OPF problems, for which zero duality holds, are investigated in \cite{Bose-2011-OPF-Allerton,2011arXiv1107.1467Z,2012arXiv1204.4419L}.
 In particular, \cite{Bose-2011-OPF-Allerton} derives zero duality results for networks that have tree topology where over satisfaction of loads is allowed.
 On the other hand, \cite{2011arXiv1107.1467Z,2012arXiv1204.4419L} provide a graphically intuitive conditions for the zero duality gap for 2-bus networks, which are later generalized  to tree topologies.
%

\BlueText{The references above suggest that the applicability of SDP approaches are limited to \emph{special network classes}. Of course, SDP relaxations can always be used to compute a lower bound on the optimal value of the primal problem. However, in practice, what is crucial is a network operating point. In general, SDP relaxations fail to provide a network operating point (i.e., a  feasible point) due to nonzero duality gap~\cite{Lesieutre-Exa-the-Lim-of-App-of-Sem-Pro-to-Pow-Flo-Pro-2011}.    
Another drawback of SDP based methods is that, even when zero duality holds, if the objective functions are non-quadratic, the dual machinery employed in constructing primal feasible solutions is not applied.  Authors in~\cite{Lesieutre-Exa-the-Lim-of-App-of-Sem-Pro-to-Pow-Flo-Pro-2011,Gopalakrishnan-BandB-OPF-2012} explore limitations of SDP approaches and give practical examples where the sufficient conditions for zero duality does not hold. }

%

Centralized methods for OPF problem, of course, exhibit a poor scalability.  On the contrary, distributed and scalable OPF solution methods are less investigated, though they are highly desirable in the context of rapidly growing real world electrical networks. Unlike centralized methods, distributed OPF solution methods are also appealing in the context of privacy and security, because they do not entail collecting possibly sensitive problem data at a central node.
%
%
%
%
In other words, when solving in a centralized manner the OPF problem in the smart grid, the power companies must rely on private information, such as the load profile of their costumers~\cite{Quinn-Pri-and-the-new-Ene-Inf-2009,Sma-Pri-for-the-Sma-Gri-Cavoukian-2010}, which might be of interest to a third party. For example, government agencies might inquire the information to profile criminal activity and insurance companies might be interested in buying the information to determine if an individual is viable for an insurance~\cite{Salinas-Pri-Pre-Ene-The-Det-in-Sma-Gri-2013}. Therefore, gathering of private information at a centralized node has raised serious concerns about personal privacy, which in turn discourages the use of centralized approaches.
\BlueText{ Interestingly, the sparsity of most electrical networks brings out an appealing decomposition structure, and therefore it is worth investigating distributed methods for the OPF problem. }

Distributed methods for the OPF problem are first studied in~\cite{Kim-Coa-Gra-Dis-OPF-1997,Kim-A-Fas-Dis-Imp-of-OPF-1999,Kim-A-Com-of-Dis-OPF-Alg-2000}, where the transmission network is divided into regions and different decomposition methods, including auxiliary problem principle, predictor-corrector proximal multiplier method, and alternating direction method are explored to solve the problem distributively among these regions. The formulation is restricted to \emph{2-region} network decompositions and border variables cannot be shared among more than 2 regions.
%
 Another approach to decentralize the problem into regions is presented in \cite{Conejo-A-Dec-Pro-Bas-on-App-New-Dir,Nogales-A-Dec-Met-App-to-the-Mul-Are-OPF-Pro,Hug-Glanzmann-Dec-OPF-Con-for-Ove-Are-in-Pow-Sys} . The method is based on solving the Karush-Kuhn-Tucker (KKT) optimality conditions where a Newton procedure is adapted. The authors provide a sufficient condition for convergent which can be interpreted as a measurement of coupling between regions.
 However, when the condition is not satisfied they rely on the generalized minimal residual method to find the Newton direction, which involves a lot of communications between entities.
Methods presented in~\cite{Bakirtzis-A-Dec-Sol-to-the-DC-OPF-of-Int-Pow-Sys-2003} are limited to DC OPF. 

\BlueText{More recent distributed algorithms are found in~\cite{2011arXiv1109.5229L,2012arXiv1211.5856D,2012arXiv1204.5226Z,Kraning_2014,Bolognani_2014,Sulc_2014,Sun_2013}.
The decentralized methods in~\cite{2011arXiv1109.5229L,2012arXiv1211.5856D,2012arXiv1204.5226Z} capitalize on the SDP relaxation, which still has the drawbacks of being specific to special classes of networks and lack of flexibility with general objective functions. Another relaxation method is presented in \cite{Kraning_2014}, where instead of the original nonconvex constraint sets, the convex hull of those are used. However, the method can result in an infeasible point to the original unrelaxed problem, entailing local methods to help construct good local points. Other recent works consider distributed methods for optimal reactive power flow in distribution networks~\cite{Bolognani_2014,Sulc_2014}.  
 Both papers first make approximations that yield a convex OPF problem and then distribute the computation by using dual decomposition~\cite{Bolognani_2014} and ADMM~\cite{Sulc_2014}. The recent work in~\cite{Sun_2013} employees ADMM to the general nonconvex OPF problem to devise a scalable algorithm. A major drawback of the method in~\cite{Sun_2013} is that its convergence is very sensitive to the initialization. In fact, the authors of~\cite{Sun_2013}  always initialize their algorithm with a point which is close to the optimal solution. However, the optimal solution is not known a priori, limiting scope of the method.}


\BlueText{Almost all the methods above can be classified as those which are general yet not scalable and those which are scalable yet not general. However, methods, which are \emph{simultaneously} general and scalable are of crucial importance in theory, as well as in practice, and therefore deserve investigations.}

\subsection{Our Contributions}
\BlueText{
The main contributions of this paper are as follows:
\begin{enumerate}
\item We develop a distributed algorithm for the \emph{general nonconvex} OPF problem. Our approach is not restricted to any special classes of networks, where zero duality holds. It also handles non-quadratic convex objective functions, unlike the SDP based distributed algorithms.


\item  We capitalize on alternating direction method of multipliers (ADMM)~\cite{Boyd:2011:DOS:2185815.2185816} to accomplish the distributed implementation (among electrical network buses) of the proposed algorithm with a \emph{little} coordination of the {neighboring} entities. Thus, the proposed algorithm is highly \emph{scalable}, which is favorable in practice.

\item In the case of subproblems, we capitalize on sequential approximations, in order to gracefully manipulate the nonconvexity issues. The approach is adopted from an existing algorithm originally proposed in~\cite{opf_form_IV} in the context of centralized OPF problem.

\item The convergent properties of the proposed algorithm are mathematically and numerically substantiated. %

\item  A number of numerical examples are provided to evaluate the performance of the proposed algorithm. 
\end{enumerate}
}

%
%
%
%
%
%
%


\subsection{Organization and Notations}

The paper is organized as follows.
 Section~\ref{sec:model} describes the system model and problem formulation.
 The solution method is presented in Section~\ref{sec:ald_devlop}. 
 In Section \ref{sec:algorithm_properties}, we discuss some fundamental properties of the algorithm.
 Numerical results are provided in Section~\ref{sec:results}.
 Finally, Section~\ref{sec:conclusions} concludes the paper.

The imaginary unit is denote by $j$, i.e., $j=\sqrt{-1}$.
 Boldface lower case and upper case letters represent vectors and matrices, respectively, and \BlueText{calligraphic} letters represent sets.
 The cardinality of $\mathcal{A}$ is denoted by $|\mathcal{A}|$. 
 \mbox{len}(\vec{x}) denotes the length of $\vec{x}$.
 The set of real and complex $n$-vectors are denoted by $\R^n$ and $\C^n$, respectively, and the set of real and complex $m\times n$ matrices are denoted by $\R^{m\times n}$ and $\C^{m\times n}$.
 We denote the real and imaginary parts of the complex number $z\in\C$ by $\Re(z)$ and $\Im(z)$, respectively.
 The set of nonnegative integers is denoted by $\mathbb{N}$, i.e., $\mathbb{N}=\{0,1,\ldots\}$.
 The superscript $(\cdot)^{\tran}$ stands for transpose.
 We use {parentheses} to construct column vectors from comma separated lists, e.g., $({\vec a}, {\vec b}, {\vec c})=[{\vec a}\tran \ {\vec b}\tran \ {\vec c}\tran]\tran$.
 We denote the diagonal block matrix with $\vec{A}_1,\cdots,\vec{A}_N$ on the diagonal by $\mbox{diag}(\vec{A}_1,\cdots,\vec{A}_N)$.
 The Hadamard product of the matrices $\vec{A}$ and $\vec{B}$ is denoted by $\vec{A} \circ \vec{B}$.
 We denote by $||{\vec x}||_2$ the $\ell_2$-norm of the vector~${\vec x}$.
 We denote the gradient of the function $f$ in the point $\vec{x}$ by $\nabla_{\vec{x}}f$.

\section{System model and Problem Formulation}  \label{sec:model}
  Consider an electrical network with $N$ buses with $\mathcal{N}= \{1,2,\ldots, N\}$ denoting the set of buses and $\mathcal{L}\subseteq \mathcal{N}\times \mathcal{N}$ the set of flow lines. We let $i_k = i_k\sr{Re}+j i_k\sr{Im}$ be the current injection and $v_k=v_k\sr{Re}+jv_k\sr{Im}$ be the voltage at bus $k\in \mathcal{N}$. Let $p_k\sr{D}+j q_k\sr{D} \in \C$ and $p_k\sr{G}+j q_k\sr{G} \in \C$ denote the complex power demand and the complex power generated by bus $k\in \mathcal{N}$, respectively.  Thus, the complex power $p_k+j q_k \in \C$ injected to bus $k$ is given by $ p_k+j q_k=(p_k\sr{G}+j q_k\sr{G})-(p_k\sr{D}+jq_k\sr{D})$.

For notational compactness, we let $\vec{p}\sr{G}$, $\vec{q}\sr{G}$,
       $\vec{p}\sr{D}$, $\vec{q}\sr{D}$,
       $\vec{p}$, $\vec{q}$,
       $\vec{i}$, $\vec{i}\sr{Re}$, $\vec{i}\sr{Im}$,
      $\vec{v}$, $\vec{v}\sr{Re}$, and $\vec{v}\sr{Im}$
     denote the vectors
       $(p_k\sr{G})_{k\in\mathcal{N}}$, $(q_k\sr{G})_{k\in \mathcal{N}}$,
       $(p_k\sr{D})_{k\in\mathcal{N}}$, $(q_k\sr{D})_{k\in \mathcal{N}}$,
       $(p_k)_{k\in \mathcal{N}}$, $(q_k)_{k\in \mathcal{N}}$,
       $(i_k)_{k\in \mathcal{N}}$, $(i_k\sr{Re})_{k\in \mathcal{N}}$, $(i_k\sr{Im})_{k\in \mathcal{N}}$,
      $(v_k)_{k\in \mathcal{N}}$, $(v_k\sr{Re})_{k\in \mathcal{N}}$, and $(v_k\sr{Im})_{k\in \mathcal{N}}$, respectively.
We denote by $i_{ls}\sr{Re} + j i_{ls}\sr{Im}\in \C$ the complex current and by $p_{ls}+jq_{ls}\in \C$ the complex power transferred from bus $l$ to the rest of the network through the flow line $(l,s)\in \mathcal{L}$.
The admittance matrix $\vec{Y}\in \C^{N \times N}$ of the network is given by
 \begin{equation}    \label{eq:admittance_matrix_definition}
  \vec{Y} =  \begin{cases}
                         y_{ll}+\sum_{(l,t)\in \mathcal{L} } y_{lt},     &    \mbox{if $l=s$}, \\
                        -y_{ls},                                                                &    \mbox{if $(l,s) \in \mathcal{L}$ },\\
                         0,                                                                          &    \mbox{otherwise},
                 \end{cases}
 \end{equation}
where $y_{ls}=g_{ls}+j b_{ls} \in \C$ is the admittance in the flow line $(l,s)\in \mathcal{L}$, and $y_{ll}=g_{ll}+j b_{ll}\in \C$ is the admittance to ground at bus $l$.
 We let $\vec{G} \in \R^{N\times N}$ and $\vec{B}\in \R^{N\times N}$ denote the real and imaginary parts of $\vec{Y}$, respectively. In particular, $[\vec{G}]_{ls}=g_{ls}$ and $[\vec{B}]_{ls}=b_{ls}$ yielding $\vec{Y}=\vec{G}+j\vec{B}$.

\subsection{Centralized formulation}

\BlueText{For fixed power demands, $\vec{p}\sr{D}$ and $\vec{q}\sr{D}$, the goal of the OPF problem is to find the optimal way to tune the variables $\vec{p}\sr{G}$, $\vec{q}\sr{G}$,  $\vec{p}$, $\vec{q}$, $\vec{i}\sr{Re}$, $\vec{i}\sr{Im}$, $\vec{v}\sr{Re}$,  $\vec{v}\sr{Im}$, ensuring that the relationships among the variables are held and system limitations are respected. 
The objective function differs between applications. In this paper we consider the minimization of  a convex cost function of real power generation. We denote by $f_k\sr{G}$ the cost of generating power at bus $k\in \mathcal{G}$, where $\mathcal{G}\subseteq \mathcal{N}$ denotes the set of generator buses. 
 The OPF problem  can now be expressed as}\footnote{\BlueText{Formulation~\eqref{eq:main_formulation} is equivalent to the OPF formulation in~\cite{Lavaei-2012-OPF-TPS}, and one can easily switch between the two formulation by using simple transformations. We use formulation~\eqref{eq:main_formulation}, because it is convenient, in terms of notations, when describing the contents in subsequent sections.}}
\begin{subequations}  \label{eq:main_formulation}
      \begin{align}
        & {\text{min}}
          & &  \sum_{k\in \mathcal{G} } f_k\sr{G} (p_k\sr{G})\\
            & \text{s. t.}
            & & \vec{i}\sr{Re}{+}j \vec{i}\sr{Im} {=}  \vec{G} \vec{v}\sr{Re} {-} \vec{B} \vec{v}\sr{Im} {+} j\left(  \vec{B} \vec{v}\sr{Re} {+} \vec{G} \vec{v}\sr{Im} \right),  \label{eq:main_formulation_current_as_a_function_of_voltages} \displaybreak[1]  \\
            & & &   p_k    {+}j q_k =  p_k\sr{G}{-}p_k\sr{D}   {+}j \left(   q_k\sr{G}{-}q_k\sr{D} \right),  k\in \mathcal{N}, \label{eq:conservation-of-power-flow1} \displaybreak[1]  \\
       & & &      i_{ls}\sr{Re}{+}j i_{ls}\sr{Im} {=} \left( \vec{c}_{ls}\tran{+} j \vec{d}_{ls}\tran\right) (v_l\sr{Re}{,}v_s\sr{Re}{,}v_l\sr{Im}{,}v_s\sr{Im}),~ (l{,}s){\in}\mathcal{L}, \label{eq:main_formulation_current_as_a_function_of_voltage_flow_line}    \displaybreak[1] \\
       & & & \vec{p}{+}j\vec{q} {=}     \vec{v}\sr{Re} {\circ} \vec{i}\sr{Re}{+} \vec{v}\sr{Im} {\circ} \vec{i}\sr{Im} {+} j\left( \vec{v}\sr{Im} {\circ} \vec{i}\sr{Re} {-} \vec{v}\sr{Re} {\circ} \vec{i}\sr{Im}   \right), \label{eq:main_formulation_power_as_a_function_of_voltage_and_current} \displaybreak[1] \\
       & & &      p_{ls} {+}jq_{ls}  {=}   v_{l}\sr{Re} i_{ls}\sr{Re}   {+}  v_l\sr{Im}  i_{ls}\sr{Im}   {+} j{(}       v_l\sr{Im} i_{ls}\sr{Re}   {-}  v_l\sr{Re}  i_{ls}\sr{Im} {)}, \notag \\
       &&& \hspace{4cm}(l,s)\in\mathcal{L}, \label{eq:main_formulation_power_as_a_function_of_voltage_and_current_flow_line}   \displaybreak[1] \\
        & & &  p_k\sr{G,min} {\leq} p_k \sr{G} {\leq}    p_k\sr{G,max},   k\in \mathcal{N}, \label{eq:main_formulation_real_power_generation_limits}    \\
       & & &   q_k\sr{G,min}                  {\leq} q_k \sr{G} {\leq}    q_k\sr{G,max},                                                                                      k\in \mathcal{N}, \displaybreak[1]\\
      & & &      (i_{ls}\sr{Re})^2{+}(i_{ls}\sr{Im})^2  {\leq}  (i_{ls}\sr{max})^2,                                                                  (l,s)\in\mathcal{L}, \label{eq:main_formulation_current_limit_flow_line} \displaybreak[1] \\
       & & &      p_{ls}^2{+} q_{ls}^2  {\leq}  (s_{ls}\sr{max})^2,                                                                                           (l,s)\in\mathcal{L},  \label{eq:main_formulation_power_limit_flow_line}  \displaybreak[1]\\
       & & &     |p_{ls}| {\leq} p_{ls}\sr{max},                                                                                  (l,s)\in\mathcal{L},  \label{eq:main_formulation_real_power_limit_flow_line} \displaybreak[1]\\
        & & &   (v_k\sr{min})^2   {\leq} (v_k\sr{Re})^2 {+} (v_k\sr{Im})^2   {\leq}      (v_k\sr{max})^2,                    k\in \mathcal{N},  \label{eq:main_formulation_voltage_limits} \displaybreak[1]
      \end{align}
\end{subequations}
%
%
where the variables are $\vec{p}\sr{G}$, $\vec{q}\sr{G}$,  $\vec{p}$, $\vec{q}$, $\vec{i}\sr{Re}$, $\vec{i}\sr{Im}$, $\vec{v}\sr{Re}$, $\vec{v}\sr{Im}$, and $i_{ls}\sr{Re}$, $i_{ls}\sr{Im}$, $p_{ls}$, $q_{ls}$ for $(l,s)\in \mathcal{L}$. Here constraint~\eqref{eq:main_formulation_current_as_a_function_of_voltages} is from that $\vec{i} = \vec{Y} \vec{v}$, \eqref{eq:conservation-of-power-flow1} is from that the conservation of power flow holds, \eqref{eq:main_formulation_current_as_a_function_of_voltage_flow_line} is from that $i_{ls}\sr{Re}  = \Re(y_{ls}(v_l-v_s) )$, \eqref{eq:main_formulation_power_as_a_function_of_voltage_and_current} is from that complex power is $(\vec{v} \circ \vec{i}^{*})$, and $i_{ls}\sr{Im}  = \Im(y_{ls}(v_l-v_s) )$ with $\vec{c}_{ls} =  (g_{ls},-g_{ls},b_{ls},-b_{ls})$ and $ \vec{d}_{ls} =  (b_{ls},-b_{ls},-g_{ls},g_{ls})$, and \eqref{eq:main_formulation_power_as_a_function_of_voltage_and_current_flow_line} is from that $p_{ls} = \Re(v_l i_{ls}^*)$ and $q_{ls} = \Im(v_l i_{ls}^*)$. Note that~\eqref{eq:main_formulation_current_as_a_function_of_voltages}-\eqref{eq:main_formulation_power_as_a_function_of_voltage_and_current_flow_line} correspond to the constraints imposed by the laws of physics associated with the electrical network. In addition, \eqref{eq:main_formulation_real_power_generation_limits}-\eqref{eq:main_formulation_voltage_limits} correspond to the constraints imposed by operational limitations, where the lower bound problem data $(\cdot)\sr{min}$ and the upper bound problem data $(\cdot)\sr{max}$ determine the boundaries of the feasible regions of power, current, as well as voltages in the network. Note that if a bus~$k$ is not a generator bus, then there is no power generation at that bus, and thus $p_k\sr{G}+j q_k\sr{G}=0$. Such situations can be easily modeled by letting
\be
p_k\sr{G,min}= p_k\sr{G,max}= q_k\sr{G,min}= q_k\sr{G,max}=0, \qquad k\in\mathcal{N} \setminus \mathcal{G} \ .
\ee

The constraints \eqref{eq:main_formulation_power_as_a_function_of_voltage_and_current}, \eqref{eq:main_formulation_power_as_a_function_of_voltage_and_current_flow_line}, and  \eqref{eq:main_formulation_voltage_limits} are nonconvex, which in turn make problem~\eqref{eq:main_formulation} nonconvex. In fact, the problem is NP-hard~\cite{Lavaei-2012-OPF-TPS}. Thus, it hinders efficient algorithms for achieving optimality. However, in the sequel, we design an efficient algorithm to address problem~\eqref{eq:main_formulation} in a decentralized manner. 

\subsection{Distributed formulation}

In this section, we derive an equivalent formulation of problem~\eqref{eq:main_formulation}, where all the constraints except for a \emph{single consistency constraint} are decoupled among the buses. In particular, the resulting formulation is in the form of general consensus problem~\cite[\S~7.2]{Boyd:2011:DOS:2185815.2185816}, where fully decentralized implementation can be realized, without any coordination of a central authority. More generally, the proposed formulation can be easily adapted to accomplish decoupling among subsets of buses, each of which corresponds to buses located in a given area, e.g., multi-area~OPF \cite{Kim-Coa-Gra-Dis-OPF-1997}.

We start by identifying the coupling constraints of problem~\eqref{eq:main_formulation}.  From constraint~\eqref{eq:main_formulation_current_as_a_function_of_voltages}, note that the current injection of each bus is affected by the voltages of its neighbors and by its own voltage. Therefore, constraint~\eqref{eq:main_formulation_current_as_a_function_of_voltages} introduces coupling between neighbors.
To decouple constraint~\eqref{eq:main_formulation_current_as_a_function_of_voltages}, we let each bus maintain local copies of the neighbors' voltages and then enforce them to agree by introducing consistency constraints.

To formally express the idea above, we first denote by $\mathcal{N}_k$ the set of bus $k$ itself and its neighboring buses, i.e., $\mathcal{N}_k= \{ k \} \cup \{ n | (k,n) \in  \mathcal{L}\}$. 
 Copies of real and imaginary parts of the voltages corresponding to buses in $\mathcal{N}_k$ is denoted by $\BlueText{ \vec{v}_k\sr{Re} } \in \R^{|\mathcal{N}_k|}$ and $\BlueText{\vec{v}_k\sr{Im}} \in \R^{|\mathcal{N}_k|}$, respectively. 
  For notational convenience, we let $\BlueText{(\vec{v}_k\sr{Re})_1} = v_k\sr{Re}$ and $\BlueText{(\vec{v}_k\sr{Im})_1} = v_k\sr{Im}$. We refer to $\vec{v}\sr{Re}$ and $\vec{v}\sr{Im}$ as \emph{real and imaginary net variables}, respectively. 
  Note that the copies of either the net variable $v_k\sr{Re}$ or $v_k\sr{Im}$ are shared among $|\mathcal{N}_k|$ entities, which we call the degree of net variable $v_k\sr{Re}$ or $v_k\sr{Im}$. The consistency  constraints are given by
 \begin{equation} \label{eq:consistency_constraints}
  \BlueText{ \vec{v}_k\sr{Re} }= \vec{E}_k \vec{v}\sr{Re}, \qquad  \BlueText{ \vec{v}_k\sr{Im} }=\vec{E}_k\vec{v}\sr{Im},
 \end{equation}
 where $\vec{E}_k \in \R^{|\mathcal{N}_k| \times N}$ is given by
 \begin{equation}
   (\vec{E}_k)_{ls} {=}
       \begin{cases}
                1       &  \text{if \BlueText{$(\vec{v}_k\sr{Re})_l$} is a local copy of $v_s\sr{Re}$ }\\
                0       & \text{otherwise} .
      \end{cases}
 \end{equation}
Note that~\eqref{eq:consistency_constraints} ensures the agreement of the copies of the net variables and that, for any bus~$k$, either \BlueText{$\vec{v}_k\sr{Re}$} or \BlueText{$\vec{v}_k\sr{Im}$} is local in the sense that they depend only on neighbors.

 The constraints \eqref{eq:main_formulation_current_as_a_function_of_voltages}-\eqref{eq:main_formulation_voltage_limits} of problem~\eqref{eq:main_formulation} can be written by using local variables \BlueText{$\vec{v}_k\sr{Re}$ and $\vec{v}_k\sr{Im}$}. In particular, we can equivalently list them as follows:
\begin{subequations}  \label{eq:dist_formulation_list}
      \begin{align}
            &  i_k\sr{Re}{+}j i_k\sr{Im} =     \vec{g}_k\tran \BlueText{\vec{v}_k\sr{Re}} {-} \vec{b}_k\tran \BlueText{\vec{v}_k\sr{Im}} {+} j\left(  \vec{b}_k\tran \BlueText{\vec{v}_k\sr{Re}} + \vec{g}_k\tran \BlueText{\vec{v}_k\sr{Im}} \right),  \label{eq:dist_formulation_current_as_a_function_of_voltages} \\
          & p_k    +j q_k =  p_k\sr{G}-p_k\sr{D}   +j \left(   q_k\sr{G}-q_k\sr{D} \right),      &  \label{eq:dist_formulation_power_injections_generator_bus}\\
        &     \bec{i}_k\sr{Re}+j  \bec{i}_k\sr{Im} = \vec{C}_k \BlueText{\vec{v}_k\sr{Re}} + \vec{D}_k \BlueText{\vec{v}_k\sr{Im}} + j\left(  \vec{D}_k \BlueText{\vec{v}_k\sr{Re}} - \vec{C}_k \BlueText{\vec{v}_k\sr{Im} }\right),    & \label{eq:dist_formulation_psysical_5} \\
       &   p_k{+}j q_k {=}    (\BlueText{\vec{v}_k\sr{Re}})_1 i_k\sr{Re} {+} (\BlueText{\vec{v}_k\sr{Im}})_1 i_k\sr{Im} {+} j\left( \BlueText{(\vec{v}_k\sr{Im})_1} i_k\sr{Re} {-} \BlueText{ (\vec{v}_k\sr{Re})_1 } i_k\sr{Im}   \right),  \label{eq:dist_formulation_power_as_a_function_of_voltage_and_current}\\
       &   \bec{p}_k {+}j\bec{q}_k  {=}  \BlueText{ (\vec{v}_k\sr{Re})_1 } \bec{i}_k\sr{Re}  {+} \BlueText{ (\vec{v}_k\sr{Im})_1 }  \bec{i}_k\sr{Re}   {+} j\left(    \BlueText{   (\vec{v}_k\sr{Im})_1}  \bec{i}_k\sr{Re}  {-} \BlueText{  (\vec{v}_k\sr{Re})_1 }    \bec{i}_k\sr{Im} \right),                                       \label{eq:dist_formulation_power_as_a_function_of_voltage_and_current_flow_line} \\
       &   p_k\sr{G,min}                  \leq p_k \sr{G} \leq    p_k\sr{G,max},         \label{eq:dist_formulation_real_power_generation_limits} \\
      &   q_k\sr{G,min}                  \leq q_k \sr{G} \leq    q_k\sr{G,max},                                                                                  &  \label{eq:dist_formulation_reac_power_generation_limits}\\
     &   (\bec{i}_{k}\sr{Re})^2_r+(\bec{i}_{k}\sr{Im})^2_r \leq (\bec{i}_{k}\sr{max})^2_r,       \quad  r{=}1,\ldots,|\mathcal{N}_k|{-}1, \label{eq:dist_formulation_current_circle} \\
     &   (\bec{p}_{k})^2_r+(\bec{q}_{k})^2_r \leq (\bec{s}_{k}\sr{max})^2_r,    \quad  r{=}1,\ldots,|\mathcal{N}_k|{-}1 ,  \label{eq:dist_formulation_power_circle}\\
       \displaybreak[1]
       &   |(\bec{p}_{k})_r| \leq (\bec{p}_{k}\sr{max})_r,     \quad  r{=}1,\ldots,|\mathcal{N}_k|{-}1,  \label{eq:dist_formulation_powerr_limits} \\
       &  \BlueText{(\vec{v}_k\sr{min})^2_r {\leq} (\vec{v}_k\sr{Re})^2_r {+} (\vec{v}_k\sr{Im})^2_r {\leq} (\vec{v}_k\sr{max})^2_r,}  \quad r{=}1,\ldots,|\mathcal{N}_k|,  \label{eq:dist_formulation_voltage_limits}
      \end{align}
\end{subequations}
where $k\in\mathcal{N}$, $\bec{i}_k\sr{Re}=(i_{kl}\sr{Re})_{l\in\mathcal{N}_k \setminus \{k\}}$, $\bec{i}_k\sr{Im}=(i_{kl}\sr{Im})_{l\in\mathcal{N}_k \setminus \{k\}}$, $\bec{p}_k=(p_{kl})_{l\in\mathcal{N}_k \setminus \{k\}}$, and  $\bec{q}_k=(q_{kl})_{l\in\mathcal{N}_k \setminus \{k\}}$, with the order kept preserved as in \BlueText{$(\vec{v}_k\sr{Re})_{1:|\mathcal{N}_k|}$} and \BlueText{$(\vec{v}_k\sr{Im})_{1:|\mathcal{N}_k|}$}. In addition, $\vec{g}_k$ (or $\vec{b}_k$) in constraint~\eqref{eq:dist_formulation_current_as_a_function_of_voltages} are obtained by first extracting the $k$-th column of $\vec{G}$ (respectively, $\vec{B}$) and then extracting the rows corresponding to the buses in $\mathcal{N}_k$, where the order of the components are preserved as in \BlueText{$\vec{v}_k\sr{Re}$} and \BlueText{$\vec{v}_k\sr{Im}$}. In addition, $\vec{C}_k\in \R^{(|\mathcal{N}_k|-1)\times |\mathcal{N}_k|}$  and $\vec{D}_k\in \R^{(|\mathcal{N}_k|-1)\times |\mathcal{N}_k|}$ in constraint~\eqref{eq:dist_formulation_psysical_5} are given by
\begin{align}
  \vec{C}_k &{=} {\left( \begin{array}{cccc}
                                                       (\vec{g}_k)_2                  & -(\vec{g}_k)_2  & \cdots & 0 \\
                                                       (\vec{g}_k)_3                  & 0                        & \cdots & 0 \\
                                                      \vdots                               & \vdots               & \ddots & \vdots \\
                                                       (\vec{g}_k)_{|\mathcal{N}_k|-1} & 0         & \cdots & -(\vec{g}_k)_{|\mathcal{N}_k|}
                                        \end{array} \right)},\\
  \vec{D}_k &{=} {\left( \begin{array}{cccc}
                                                       (\vec{b}_k)_2                  & -(\vec{b}_k)_2                    & \cdots & 0 \\
                                                       (\vec{b}_k)_3                  & 0                        & \cdots & 0 \\
                                                      \vdots                                 & \vdots                               & \ddots & \vdots \\
                                                       (\vec{b}_k)_{|\mathcal{N}_k|-1} & 0                                        & \cdots & -(\vec{b}_k)_{|\mathcal{N}_k|}
                                        \end{array} \right)}.
\end{align}
Moreover, $\BlueText{\vec{v}_k\sr{min}}, \BlueText{\vec{v}_k\sr{max}}, \bec{i}_{k}\sr{max}, \bec{s}_{k}\sr{max}$, and $\bec{p}_{k}\sr{max}$ of constraints~\eqref{eq:dist_formulation_real_power_generation_limits}-\eqref{eq:dist_formulation_voltage_limits} are chosen in a straightforward manner~[\emph{cf} \eqref{eq:main_formulation_real_power_generation_limits}-\eqref{eq:main_formulation_voltage_limits}].

Finally, for notational convenience, associated with each bus~$k$, we denote by
\be \label{eq:local_vec}
\vec{z}_k=({p}_k\sr{G}, {q}_k\sr{G}, {p}_k, {q}_k, {i}_k\sr{Re}, {i}_k\sr{Im},\BlueText{ \vec{v}_k\sr{Re}}, \BlueText{\vec{v}_k\sr{Im}}, \bec{i}_{k}\sr{Re}, \bec{i}_{k}\sr{Im}, \bec{p}_{k}, \bec{q}_{k})
\ee
the local variables of bus~$k$, by $\boldsymbol{\alpha}_k(\vec{z}_k)=0$ the affine constraints~\eqref{eq:dist_formulation_current_as_a_function_of_voltages}-\eqref{eq:dist_formulation_psysical_5}, by $\boldsymbol{\lambda}_k(\vec{z}_k)=0$ the nonlinear equality constraint~\eqref{eq:dist_formulation_power_as_a_function_of_voltage_and_current}, by $\boldsymbol{\mu}_k(\vec{z}_k)=0$ the nonlinear equality constraint~\eqref{eq:dist_formulation_power_as_a_function_of_voltage_and_current_flow_line}, by $\boldsymbol{\beta}_k(\vec{z}_k)\leq 0$ the linear convex inequality constraints~\eqref{eq:dist_formulation_real_power_generation_limits}, \eqref{eq:dist_formulation_reac_power_generation_limits}, \eqref{eq:dist_formulation_powerr_limits}, and by $\boldsymbol{\gamma}_k(\vec{z}_k)\leq 0$ the nonlinear convex inequality constraints~\eqref{eq:dist_formulation_current_circle} and \eqref{eq:dist_formulation_power_circle} as we will see next.\footnote{\BlueText{The functions $\beta_k$ and $\gamma_k$ depend on $ p_k\sr{G,min}$, $ p_k\sr{G,max}$, $q_k\sr{G,min}$, $q_k\sr{G,max}$, $\bec{i}_{k}\sr{max}$, $\bec{s}_{k}\sr{max}$, and $\bec{p}_{k}\sr{max}$ which are intentionally omitted for clarity and space reasons.  }   }

 Now we can express the distributed formulation of problem~\eqref{eq:main_formulation} as
\begin{subequations}  \label{eq:dist_formulation}
      \begin{align}
        & {\text{min}}
          & &  \sum_{k\in \mathcal{G} } f_k\sr{G} (p_k\sr{G})  \label{eq:dist_formulation_obj}\\
            & \text{s. t.}
          & &  \vec{z}_k{=}{(}{p}_k\sr{G}{,} {q}_k\sr{G}{,} {p}_k{,} {q}_k{,} {i}_k\sr{Re}{,} {i}_k\sr{Im}{,} \vec{x}_k\sr{Re}{,} \vec{x}_k\sr{Im}{,} \notag \\
          & & &   \hspace{2.4cm} \bec{i}_{k}\sr{Re}, \bec{i}_{k}\sr{Im}, \bec{p}_{k}, \bec{q}_{k}), \quad k\in \mathcal{N},  \label{eq:local_variables_set} \displaybreak[2] \\
       & & &    \big( \boldsymbol{\alpha}_k(\vec{z}_k), \boldsymbol{\lambda}_k(\vec{z}_k), \boldsymbol{\mu}_k(\vec{z}_k)\big)=\vec{0}, \quad k\in \mathcal{N}, \displaybreak[1] \\
       & & &  \big(\boldsymbol{\beta}_k(\vec{z}_k), \boldsymbol{\gamma}_k(\vec{z}_k) \big)\leq \vec{0},  \hspace{1.6cm} k\in \mathcal{N}, \displaybreak[1] \\
       & & & \BlueText{ (\vec{v}_k\sr{min})^2_r \leq (\vec{v}_k\sr{Re})^2_r + (\vec{v}_k\sr{Im})^2_r \leq (\vec{v}_k\sr{max})^2_r}, \notag  \\
       & & &    \hspace{2cm} r=1,\ldots,|\mathcal{N}_k|, \hspace{0.6cm}   k\in \mathcal{N}, \label{eq:nonconvex_inequality}  \displaybreak[1] \\
       & & &  \BlueText{ \vec{v}_k\sr{Re} } {+} j \BlueText{\vec{v}_k\sr{Im}} {=} \vec{E}_k \vec{v}\sr{Re} {+} j \vec{E}_k \vec{v}\sr{Im}, \hspace{1cm}   k \in \mathcal{N}, \label{eq:dist_formulation_consistency_constraints}
      \end{align}
\end{subequations}
where the variables are ${p}_k\sr{G}$, ${q}_k\sr{G}$, ${p}_k$, ${q}_k$, ${i}_k\sr{Re}$, ${i}_k\sr{Im}$, \BlueText{$\vec{v}_k\sr{Re}$}, \BlueText{$\vec{v}_k\sr{Im}$}, $\bec{i}_{k}\sr{Re}$, $\bec{i}_{k}\sr{Im}$, $\bec{p}_{k}$, $\bec{q}_{k}$, $\vec{z}_k$ for $k\in \mathcal{N}$ and $\vec{v}\sr{Re}$ and $\vec{v}\sr{Im}$. Note that \eqref{eq:dist_formulation_consistency_constraints} establishes the consistency constraints~[\emph{cf} \eqref{eq:consistency_constraints}], which affirms the consistency among neighbor voltages. The coupling in the original centralized formulation~\eqref{eq:main_formulation} has been subsumed in the consistency constraint~\eqref{eq:dist_formulation_consistency_constraints}, which results in the form of general consensus problem~\cite[\S~7.2]{Boyd:2011:DOS:2185815.2185816}, where decomposition methods can gracefully be~applied.

\section{Distributed  solution method}\label{sec:ald_devlop}
In this section, we present our distributed algorithm to the OPF problem~\eqref{eq:dist_formulation}. In particular, we use the ADMM method as basis for our algorithm development, where we have fast convergence properties, compared to the dual decomposition \cite{Boyd:2011:DOS:2185815.2185816}.  The use of ADMM method is promising in the sense that it works on many nonconvex problems as a good heuristic \cite[\S~9]{Boyd:2011:DOS:2185815.2185816}.
 Once the solution method is established, we investigate the properties in Section \ref{sec:algorithm_properties}.

\subsection{Outline of the algorithm}

For notational simplicity, we let  \BlueText{$\vec{v}_k$} and $\bec{E}_k$ denote $(\BlueText{\vec{v}_k\sr{Re}},\BlueText{\vec{v}_k\sr{Im}})$ and $\mbox{diag}(\vec{E}_k, \vec{E}_k)$, respectively, for each $k\in \mathcal{N}$. Moreover, we let $\vec{v}$ denote $(\vec{v}\sr{Re}, \vec{v}\sr{Im})$.
The ADMM essentially minimizes the augmented Lagrangian associated with the problem in an iterative manner. Particularized to our problem~\eqref{eq:dist_formulation}, The partial augmented Lagrangian with respect to the consistency constraints~\eqref{eq:dist_formulation_consistency_constraints} (i.e., $\BlueText{\vec{v}_k} = \bec{E}_k \vec{v}$) is given by
 \begin{multline}\label{eq:partial_augmented_Lagrangian}
    L_{\rho}(\vec{p}\sr{G},(\BlueText{\vec{v}_k})_{k\in\mathcal{N}},\vec{v},(\vec{y}_k)_{k\in \mathcal{N}}) = \sum_{k\in\mathcal{G}} f_k\sr{G} (p_k\sr{G}) \\
     + \sum_{ k \in \mathcal{N}} \left( \vec{y}_k\tran(\BlueText{\vec{v}_k}-\bec{E}_k\vec{v})  {+} \frac{\rho}{2}||\BlueText{\vec{v}_k}-\bec{E}_k \vec{v} {||}_2^2 \right),
 \end{multline}
where $\vec{y}_k$ is dual variable associated with~\eqref{eq:dist_formulation_consistency_constraints} and $\rho$ is called the \emph{penalty parameter}. Together with the separability of~\eqref{eq:partial_augmented_Lagrangian} among $k\in\mathcal{N}$, steps of ADMM is formally expressed below.


\noindent\rule{\linewidth}{0.3mm}
\\
\emph{Algorithm 1}: \ \textsc{\small{ADMM} for distributed OPF (ADMM-DOPF)}
\begin{enumerate}
\item Initialization:  Set $n=0$ and initialize $\vec{y}_k^{(n)}$ and $\vec{v}^{(n)}$.
\item \label{al:ADMM_private_variable_update}
          Private variable update:  Set $\vec{y}_k=\vec{y}_k^{(n)}$ and $\vec{v}=\vec{v}^{(n)}$. Each bus $k\in \mathcal{N}$ updates $\vec{x}_k$ locally, where we let $(\vec{z}^{(n+1)}_k,\vec{u}^{(n+1)}_k)$ be the primal and dual (possibly) optimal variables achieved for the following problem:
					 \begin{subequations}  \label{eq:ADMM_private_variable_update_generator_bus}
					      \begin{align}
					        & {\text{min}}
					          & &   {f_k\sr{G} (p_k\sr{G}){+}\vec{y}_k\tran(\BlueText{\vec{v}_k}{-}\bec{E}_k\vec{v}){+} \frac{\rho}{2}||\BlueText{\vec{v}_k}{-}\bec{E}_k \vec{v} ||_2^2}, \label{eq:ADMM_step2_obj}\\
					            & \text{s. t.}
          & & {\vec{z}}_k{=}({p}_k\sr{G}{,} {q}_k\sr{G}{,} {p}_k{,} {q}_k{,} {i}_k\sr{Re}{,} {i}_k\sr{Im}{,} \BlueText{\vec{v}_k\sr{Re}}{,} \BlueText{\vec{v}_k\sr{Im}}{,} \notag \\ 
          & & &   \hspace{2.6cm}   \bec{i}_{k}\sr{Re}, \bec{i}_{k}\sr{Im}, \bec{p}_{k}, \bec{q}_{k}) && \label{eq:local_variables_set_sub} \\
          & & &  \boldsymbol{\alpha}_k(\vec{z}_k)=\vec{0}     &&    \label{eq:affine_sub}\\
        & & &       \boldsymbol{\lambda}_k(\vec{z}_k)=\vec{0}  && \label{eq:nonlinear_equality1_sub} \\
       & & &\boldsymbol{\mu}_k(\vec{z}_k)=\vec{0} &&  \label{eq:nonlinear_equality2_sub} \displaybreak[1]\\
       & & &    \boldsymbol{\beta}_k(\vec{z}_k)\leq \vec{0}                                      &&  \label{eq:convex_inequality_sub} \\
       & & &    \boldsymbol{\gamma}_k(\vec{z}_k)\leq \vec{0}                                      &&  \label{eq:convex_nonlinear_inequality_sub} \displaybreak[1]\\
       & & &  \BlueText{ (\vec{v}_k\sr{min})^2_r {\leq} (\vec{v}_k\sr{Re})^2_r {+} (\vec{v}_k\sr{Im})^2_r {\leq} (\vec{v}_k\sr{max})^2_r,}  \notag \\
       & & & \hspace{2.4cm}               r=1,\ldots,|\mathcal{N}_k|,&  \label{eq:nonconvex_inequality_sub}
					      \end{align}
					\end{subequations}
%
where the variables are ${p}_k\sr{G}$, ${q}_k\sr{G}$, ${p}_k$, ${q}_k$, ${i}_k\sr{Re}$, ${i}_k\sr{Im}$, \BlueText{$\vec{v}_k\sr{Re}$}, \BlueText{$\vec{v}_k\sr{Im}$}, $\bec{i}_{k}\sr{Re}$, $\bec{i}_{k}\sr{Im}$, $\bec{p}_{k}$, $\bec{q}_{k}$, and $\vec{z}_k$. We denote by \BlueText{$\vec{v}^{(n+1)}_k$} the part of $\vec{z}^{(n+1)}_k$ corresponding to $(\vec{v}_k\sr{Re},\vec{v}_k\sr{Im})$.

\item Net variable update: We let $\vec{v}^{(n+1)}$ be the solution to the problem
                                                       \begin{multline} \label{eq:ADMM_net_variable_update}                                         
					         {\text{min}}
					           {\sum_{k\in \mathcal{N}}}  \vec{y}_k\tran( \BlueText{\vec{v}_k^{(n{+}1)} }{-}\bec{E}_k\vec{v})  {+}
                                                                                                              \frac{\rho}{2}|| \BlueText{ \vec{v}_k^{(n+1)} } {-}\bec{E}_k \vec{v} ||_2^2 ,
                                          \end{multline}
					where the variable is $\vec{v}$.
\item Dual variable update: Each bus $k\in \mathcal{N}$ updates its dual variable $\vec{y}_k$ as
                                                    \begin{equation} \label{eq:ADMM_dual_variable_update}
                                                      \vec{y}_k^{(n+1)} = \vec{y}_k^{(n)} + \rho ( \BlueText{ \vec{v}_k^{(n+1)} } - \bec{E}_k\vec{v}^{(n+1)} ) .
                                                    \end{equation}
\item Stopping criterion: Set $n\asign n+1$. If stopping criterion is not met go to step \ref{al:ADMM_private_variable_update}, otherwise \textsc{stop} and return $\big(\vec{z}^{(n)}, \vec{v}^{(n)}, \vec{u}^{(n)}, \vec{y}^{(n)}\big)=\big((\vec{z}^{(n)}_k)_{k\in\mathcal{N}}, \vec{v}^{(n)}, (\vec{u}^{(n)}_k)_{k\in\mathcal{N}}, (\vec{y}^{(n)}_k)_{k\in\mathcal{N}}\big)$.
\end{enumerate}
\vspace{-3mm}
\rule{\linewidth}{0.3mm}

\noindent
 The first step initializes the net and dual variables.
 In the second step each bus solves a nonconvex optimization problem in order to update its private variable (see Section~\ref{sec:ald_devlop_private_variable_update}).
 In the third step the net variable is updated by solving the unconstrained quadratic optimization problem~\eqref{eq:ADMM_net_variable_update}, which has a close form solution.
 The net variable update can be done in a distributed fashion with a light communication protocol (see Section~\ref{sec:ald_devlop_net_and_dual_variable_update}).
 The \BlueText{fourth} step is the dual variable update, which can be done locally on each bus (see Section~\ref{sec:ald_devlop_net_and_dual_variable_update}).
 The fifth step is the stopping criterion.
 Natural stopping criterions include 1) running the ADMM-DOPF algorithm for a fixed number of iterations, 2) running the ADMM-DOPF algorithm till the decrement between the local- and net variables of each bus $k$ ($||\vec{E}_k\vec{v}- \BlueText{\vec{v}_k}||_2$) is below a predefined threshold 3) running the ADMM-DOPF algorithm till the objective value decrement between two successive iterations is below a predefined threshold.
\noindent In the sequel, we discuss in detail the algorithm steps~\eqref{eq:ADMM_private_variable_update_generator_bus}-\eqref{eq:ADMM_dual_variable_update}.

\subsection{The subproblems: Private variable update} \label{sec:ald_devlop_private_variable_update}


\BlueText{
 In this section, problem~\eqref{eq:ADMM_private_variable_update_generator_bus} is considered. 
 Since Problem~\eqref{eq:ADMM_private_variable_update_generator_bus} is NP-hard, only exponentially complex global methods can guarantee its optimality. 
 We capitalize on sequential convex approximations~\cite{Boyd_SCP} to design an algorithm, which is efficient compared to global methods. Similar techniques are used in~\cite{opf_form_IV} for centralized OPF, which we use as basis for designing our subproblem algorithm.

}

We start by noting that constraints~\eqref{eq:local_variables_set_sub}, \eqref{eq:affine_sub}, \eqref{eq:convex_inequality_sub}, and \eqref{eq:convex_nonlinear_inequality_sub} are convex as opposed to constraints \eqref{eq:nonconvex_inequality_sub}, \eqref{eq:nonlinear_equality1_sub}, and \eqref{eq:nonlinear_equality2_sub}, which are clearly nonconvex. The idea is to approximate the nonconvex constraints.

\begin{figure}[t]
        \centering
        \begin{subfigure}[b]{0.3\linewidth}
             \centering
             \begin{tikzpicture}

               \def\ringa{( 0,0) circle (0.4) ( 0,0) circle (1)}
               \begin{scope}[thick,font=\scriptsize]

                 \draw [->] (-1.3,0) -- (1.3,0) node [above ]  {\tiny$(\vec{v}_k\sr{Re})_r$};
                 \draw [->] (0,-1.3) -- (0,1.3) node [ right] {\tiny$(\vec{v}_k\sr{Im})_r$};

                  \draw ( -0.4,3pt) -- ( -0.4,-3pt);
                  \draw ( -1,3pt) -- ( -1,-3pt); 
                \node at (-1, 0.2) { $A$}; 
                \node at (-0.4, 0.2) { $B$}; 

               \end{scope}

               \begin{scope}[even odd rule, fill opacity=0.2]
                  \fill[fill=orange] \ringa;
                   \end{scope}
                   \draw \ringa;

                \end{tikzpicture}
                \caption{\subCaptionFontSize $\mathcal{X}^k_r$}
                \label{fig:subproblems_doughnut_set}
        \end{subfigure}
         ~
        \begin{subfigure}[b]{0.3\linewidth}
             \centering
             \begin{tikzpicture}
                 \begin{scope}[thick,font=\scriptsize]

                 \draw [->] (-1.3,0) -- (1.3,0) node [above ]  {\tiny$(\vec{v}_k\sr{Re})_r$};
                 \draw [->] (0,-1.3) -- (0,1.3) node [ right] {\tiny$(\vec{v}_k\sr{Im})_r$};

                  \draw ( -0.4,3pt) -- ( -0.4,-3pt);
                  \draw ( -1,3pt) -- ( -1,-3pt);
                \node at (-1, 0.2) { $A$}; 
                \node at (-0.4, 0.2) { $B$};

               \end{scope}
               \def\ringa{( 0,0) circle (0.4) ( 0,0) circle (1)}
               \def\lineb{  (-0.55, 1.1157) -- (1.2, -0.6343)}
               \def\linec{ (-0.55, 1.1157) --  (1.2,-0.6343)--(1.6,1.6071)--(-0.55, 1.1157)}
               \def\pointc{  (-1.5, 1.5) }

                   \begin{scope}[even odd rule, fill opacity=0.2]
                      \clip \ringa;
                      \fill[fill=orange] \linec;

                   \end{scope}
                   \draw \ringa;
                    \draw \lineb;
                    \fill ( 0.5 , 0.5) circle [radius=0.05];
                    \node [above] at (0.5,0.45) {\tiny $\check{\vec{v}}_l$};

                    \draw [-] (0,0) -- ( 0.5 , 0.5); 

                    \draw [-] (0.283,0.283) -- ( 0.183 , 0.383); 
                    \draw [-] (0.383,0.383) -- ( 0.283 , 0.483); 
                    \draw [-] ( 0.183 , 0.383)-- ( 0.283 , 0.483); 

                    \fill ( 0.283 , 0.283) circle [radius=0.05];
                    \node [right] at (0.283,0.283) {\tiny $C$};

                \end{tikzpicture}
                \caption{\subCaptionFontSize $\check{\mathcal{X}}^k_r$}
                \label{fig:subproblems_subset_XX}
        \end{subfigure}
        ~
        \begin{subfigure}[b]{0.3\linewidth}
             \centering
             \begin{tikzpicture}
                 \begin{scope}[thick,font=\scriptsize]

                 \draw [->] (-1.3,0) -- (1.3,0) node [above ]  {\tiny$(\vec{v}_k\sr{Re})_r$};
                 \draw [->] (0,-1.3) -- (0,1.3) node [ right] {\tiny${(\vec{v}_k\sr{Im})_r}$};

                  \draw ( -0.4,3pt) -- ( -0.4,-3pt);
                  \draw ( -1,3pt) -- ( -1,-3pt);
                \node at (-1, 0.2) { $A$}; 
                \node at (-0.4, 0.2) { $B$};

               \end{scope}
               \def\ringa{( 0,0) circle (0.4) ( 0,0) circle (1)}
               \def\lineb{  (-0.55, 1.1157) -- (1.2, -0.6343)}
               \def\linec{ (-0.55, 1.1157) --  (1.2, -0.6343)--(1.6,1.6071)--(-0.55, 1.1157)}
               \def\pointc{  (-1.5, 1.5) }
               \def\octagon{  (1,-0.4) -- (1,0.4)--(0.4,1)--(-0.4,1)--(-1,0.4)--(-1,-0.4)-- (-0.4,-1)--(0.4,-1)--(1,-0.4)}

                   \begin{scope}[even odd rule, fill opacity=0.2]
                      \clip \octagon;;
                      \fill[fill=orange] \linec;

                   \end{scope}
                   \draw \ringa;
                    \draw \lineb;
                    \draw \octagon;
                    \fill ( 0.5 , 0.5) circle [radius=0.05];
                    \node [above] at (0.5,0.45) {\tiny $\check{\vec{v}}_l$};

                \end{tikzpicture}
                \caption{\subCaptionFontSize $\check{\mathcal{X}}^k_r$}
                \label{fig:subproblems_subset_XXX}
        \end{subfigure}
        \caption{\captionFontSize The feasible set of $((\vec{v}\sr{Re}_k)_r,(\vec{v}\sr{Im}_k)_r)$, where $A=(\vec{v}_k\sr{max})_r$ and $B=(\vec{v}_k\sr{min})_r$. }
        \label{f}
        \vspace{-5mm}
\end{figure}
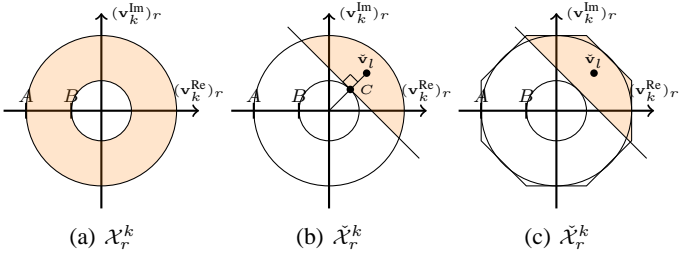

In the case of~\eqref{eq:nonconvex_inequality_sub}, we note that for any $r\in \{ 1,\ldots , |\mathcal{N}_k| \}$, the values of \BlueText{$(\vec{v}_k\sr{Re})_l$} and \BlueText{$(\vec{v}_k\sr{Im})_l$} represent a donut, see Fig.~\ref{fig:subproblems_doughnut_set}. In other words, the 2-dimensional set
 \be
     \mathcal{X}^k_r {=} \big\{  \left( \BlueText{(\vec{v}\sr{Re}_k)_r,(\vec{v}\sr{Im}_k)_r}  \right) \in \R^2   \big|  \BlueText{ (\vec{v}\sr{min}_k)_r^2{\leq} (\vec{v}\sr{Re}_k)_r^2  {+} (\vec{v}\sr{Im}_k)_r^2 {\leq}   (\vec{v}\sr{max}_k)_r^2 }  \big\} \notag
 \ee
 is a donut, which is clearly nonconvex. 

We approximate the nonconvex set $\mathcal{X}^k_r$ by considering a convex subset of $\mathcal{X}^k_r$ instead, which we denote by $\check{\mathcal{X}}^k_r$, see Fig.~\ref{fig:subproblems_subset_XX}.
 To do this, we simply consider the hyperplane tangent to the inner circle of the donut at the point $C$ in Fig.~\ref{fig:subproblems_subset_XX}.
 Specifically, given a point $\BlueText{((\check{\vec{v}}\sr{Re}_k)_r,(\check{\vec{v}}\sr{Im}_k)_r)}\in\mathcal{X}^k_r$, $\check{\mathcal{X}}^k_r$ is the intersection of $\mathcal{X}^k_r$ and the halfspace
 \begin{equation}\label{eq:halfspace}
      \big\{  \left( \BlueText{(\vec{v}\sr{Re}_k)_r,(\vec{v}\sr{Im}_k)_r } \right) {\in} \R^2  \left|  \ a_r  \BlueText{(\vec{v}_k\sr{Re})_r} {+}b_r \BlueText{(\vec{v}_k\sr{Im})_r} {\geq} c_r   \right.   \big\}  ,
 \end{equation}
 where \vspace{-2mm}
 \begin{align}
     a_r  = \text{sign} \BlueText{ (  \check{\vec{v}}\sr{Re}_k)_r }\ \sqrt{   \frac{ \BlueText{ (\vec{v}_k\sr{min})_r^2  }}{1+ ( \BlueText{ (\check{\vec{v}}\sr{Im}_k)_r}/ \BlueText{(\check{\vec{v}}\sr{Re}_k)_r )^2}} }, \notag \\
     b_r  = a_r\left(  \frac{ \BlueText{ (\check{\vec{v}}\sr{Im}_k)_r } }{ \BlueText{ (\check{\vec{v}}\sr{Re}_k)_r  } }  \right) \ ,\hspace{1cm}  
     c_r = \BlueText{ (\vec{v}_k\sr{min})_r^2 }\ , \notag
 \end{align}
 if $\BlueText{(\check{\vec{v}}\sr{Re}_k)_r} \neq 0$ and 
  \be
     a_r = 0, \hspace{1cm} 
     b_r  = \text{sign} \BlueText{ (\check{\vec{v}}\sr{Im}_k)_r }, \hspace{1cm}
     c_r =  \BlueText{ (\vec{v}_k\sr{min})_r }, \notag
 \ee
 if $\BlueText{(\check{\vec{v}}\sr{Re}_k)_r}= 0$.
In the case of nonlinear nonconvex constraints \eqref{eq:nonlinear_equality1_sub} and \eqref{eq:nonlinear_equality2_sub}, we capitalized on the well known Taylor's approximation. Specifically, given a point $\hat{\vec{z}}_k$, we denote by $\hat{\boldsymbol{\lambda}}_k^{\hat{\vec{z}}_k}$ the first order Taylor's approximation of $\boldsymbol{\lambda}_k$ at $\hat{\vec{z}}_k$. Similarly, we denote by $\hat{\boldsymbol{\mu}}_k^{\hat{\vec{z}}_k}$ the first order Taylor's approximation of $\boldsymbol{\mu}_k$ at $\hat{\vec{z}}_k$. The approximation is refined in an iterative manner until a stopping criterion is satisfied.

 It is worth noting that to construct the functions $\hat{\boldsymbol{\mu}}_k^{\hat{\vec{z}}_k}$ and $\hat{\boldsymbol{\lambda}}_k^{\hat{\vec{z}}_k}$ one only needs the values of \BlueText{$\hat{\vec{v}}_k$}, where \BlueText{$\hat{\vec{v}}_k$} is the component of $\hat{\vec{z}}_k$ corresponding to \BlueText{$(\vec{v}_k\sr{Re},\vec{v}_k\sr{Im})$}.
  \footnote{This follows directly from the definition of the first order Taylor approximation and equations~\eqref{eq:dist_formulation_current_as_a_function_of_voltages} and~\eqref{eq:dist_formulation_psysical_5}.}


By using the constraint approximations discussed above, we design a subroutine to perform step~2 of the ADMM-DOPF algorithm. The outline of this successive approximation algorithm is given as follows.

 \noindent\rule{\linewidth}{0.3mm}
\\
\emph{Algorithm 2}: \ \textsc{\small{Subroutine for step~2 of the ADMM-DOPF}}
\begin{enumerate}
   \item Initialize: Given $\vec{v}$ and $\vec{y}_k$ from ADMM-DOPF $n$th iteration. Set $(\vec{v}\sr{Re},\vec{v}\sr{Im})=\vec{v}$. For all $r\in\{1,\ldots,|\mathcal{N}_k|\}$, set $( \BlueText{(\check{\vec{v}}\sr{Re}_k)_r,(\check{\vec{v}}\sr{Im}_k)_r})=(( \vec{E}_k \vec{v}\sr{Re})_r,(\vec{E}_k \vec{v}\sr{Im})_r)$ and construct~$\check{\mathcal{X}}^k_r$. Let $m=1$ and initialize~$\hat{\vec{z}}_k$.

   \item \label{al:subproblem_algorithm_private_variable_update}
            Solve the approximated subproblem:
                \begin{subequations}  \label{eq:the_approximated_subproblem}
					      \begin{align}
					        & {\text{min}}
					          & &   f_k\sr{G} (p_k\sr{G}) {+} \vec{y}_k\tran (\BlueText{\vec{v}_k}{-}\bec{E}_k\vec{v}) {+}  \frac{\rho}{2}||\BlueText{\vec{v}_k}{-}\bec{E}_k \vec{v}||_2^2 \label{eq:the_approximated_subproblem_obj}\\
					            & \text{s. t.}
          & &  \vec{z}_k=({p}_k\sr{G}, {q}_k\sr{G}, {p}_k, {q}_k, {i}_k\sr{Re}, {i}_k\sr{Im}, \BlueText{\vec{v}_k\sr{Re}}, \BlueText{\vec{v}_k\sr{Im}}, \nonumber\\
          & & &  \hspace{2.6cm} \bec{i}_{k}\sr{Re}, \bec{i}_{k}\sr{Im}, \bec{p}_{k}, \bec{q}_{k}) && \label{eq:local_variables_set_sub_arrox} \\
          & & &  \boldsymbol{\alpha}_k(\vec{z}_k)=\vec{0}    &&    \label{eq:affine_sub_approx}\\
        & & &     \hat{\boldsymbol{\lambda}}_k^{\hat{\vec{z}}_k}(\vec{z}_k)=\vec{0}  && \label{eq:nonlinear_equality1_sub_approx} \\
       & & &\hat{\boldsymbol{\mu}}_k^{\hat{\vec{z}}_k}(\vec{z}_k)=\vec{0} &&  \label{eq:nonlinear_equality2_sub_approx}\\
       & & &    \boldsymbol{\beta}_k(\vec{z}_k)\leq \vec{0}                                      &&  \label{eq:convex_inequality_sub_approx} \\
       & & &    \boldsymbol{\gamma}_k(\vec{z}_k)\leq \vec{0}                                      &&  \label{eq:convex_nonlinear_inequality_sub_approx} \\
       & & & ( \BlueText{(\vec{v}\sr{Re}_k)_r,(\vec{v}\sr{Im}_k)_r})\in\check{\mathcal{X}}^k_r,  ~ r{=}1,\ldots,|\mathcal{N}_k|,&  \label{eq:nonconvex_inequality_sub_approx} 
					      \end{align}
					\end{subequations}
	    where the variables are ${p}_k\sr{G}$, ${q}_k\sr{G}$, ${p}_k$, ${q}_k$, ${i}_k\sr{Re}$, ${i}_k\sr{Im}$, \BlueText{$\vec{v}_k\sr{Re}$}, \BlueText{$\vec{v}_k\sr{Im}$}, $\bec{i}_{k}\sr{Re}$, $\bec{i}_{k}\sr{Im}$, $\bec{p}_{k}$, $\bec{q}_{k}$, and $\vec{z}_k$. The solution corresponding to the variable  $\vec{z}_k$ is denoted by $\vec{z}^{(m)}_k$ and all the dual optimal variables are denoted by $\vec{u}^{(m)}_k$.
     \item Stopping criterion:  If stopping criterion is not met, set $\hat{\vec{z}}_k=\vec{z}^{(m)}_k$, $m\asign m+1$ and go to step~\ref{al:subproblem_algorithm_private_variable_update}. Otherwise \textsc{stop} and return $(\vec{z}^{(m)}_k, \vec{u}^{(m)}_k)$.
\end{enumerate}
\vspace{-3mm}
\rule{\linewidth}{0.3mm}

   The initialization in the first step is done by setting $\BlueText{\hat{\vec{v}}_k}=(\BlueText{\check{\vec{v}}\sr{Re}_k,\check{\vec{v}}\sr{Im}_k})$, where \BlueText{$\hat{\vec{v}}_k$} is the component of $\hat{\vec{z}}_k$ corresponding to the variable \BlueText{$(\vec{v}_k\sr{Re},\vec{v}_k\sr{Im})$}.
 The rest of the vector $\hat{\vec{z}}_k$ is then initialized according to equations~\eqref{eq:dist_formulation_current_as_a_function_of_voltages}-\eqref{eq:dist_formulation_power_as_a_function_of_voltage_and_current_flow_line}, which have a unique solution when \BlueText{$\hat{\vec{v}}_k$} is given.
   The second step involves solving a convex optimization problem and the third step is the stopping criterion.
  A natural stopping criterion is to run the algorithm till the decrement between two successive iterations is below a certain predefined threshold,  i.e., \BlueText{$||\vec{z}^{(m+1)}_k-\vec{z}^{\BlueText{(m)}}_k|| < \epsilon$} for a given $\epsilon>0$.
  However, since \BlueText{ $\vec{z}_k$ only depends on $\vec{v}_k$}, the component related to \BlueText{$(\vec{v}_k\sr{Re},\vec{v}_k\sr{Im})$}, we use
 \begin{equation}\label{eq:sub_problem_stopping_criterion_epsilon}
    \BlueText{
    || \vec{v}_k^{(m+1)} -  \vec{v}_k^{(m)}||_2 < \epsilon\sr{sub},
    }
 \end{equation}
 where \BlueText{ $\vec{v}_k^{(m)}$ and $\vec{v}_k^{(m+1)}$ are the components of $\vec{z}_k^{(m)}$ and $\vec{z}^{(m+1)}_k$}, respectively, corresponding to the variable \BlueText{$\vec{v}_k$} and $\epsilon\sr{sub}>0$ is a given threshold.
 Furthermore, we do not need to reach the minimum accuracy in every ADMM iteration, but only as the ADMM method \BlueText{progresses}.
 Therefore,  it might be practical to set an upper bound on the number of iterations, i.e
 \begin{equation} \label{eq:sub_problem_stopping_criterion_max_iter}
              m \geq \texttt{max\_iter},
  \end{equation}
  for some $\texttt{max\_iter} \in \mathbb{N}$.

\subsection{On the use of quadratic programming (QP) solvers} \label{sec:ald_devlop_qp_solver}
      Problem~\eqref{eq:the_approximated_subproblem} can be efficiently solved by using general interior-point algorithms for convex problems. However, even higher efficiencies are achieved if problem~\eqref{eq:the_approximated_subproblem} can be handled by \emph{specific} interior-point algorithms. For example, if the objective function~\eqref{eq:the_approximated_subproblem_obj} is quadratic, sophisticated QP solvers can be easily employed. See the appendix for the detail.   


\subsection{Net variables and dual variable updates} \label{sec:ald_devlop_net_and_dual_variable_update}

Note that the net variable $\vec{v}^{(n+1)}$ is the unique solution of the unconstrained convex quadratic optimization problem~\eqref{eq:ADMM_net_variable_update}, and is given by
 \begin{align}
     \vec{v}^{(n+1)}
     =&\left( \ssum{k\in \mathcal{N}} \bec{E}_k\tran \bec{E}_k \right)^{-1} {\ssum{k\in \mathcal{N}}} \bec{E}_k \tran \left( \BlueText{ \vec{v}_k^{(n+1)} }{+} \frac{1}{\rho} \vec{y}_k^{(n)} \right), \label{eq:netvariable_1}\\
     =&\left( \ssum{k\in \mathcal{N}} \bec{E}_k\tran \bec{E}_k \right)^{-1} \ssum{k\in \mathcal{N}} \bec{E}_k \tran \BlueText{\vec{v}_k^{(n+1)} }, \label{eq:netvariable_2}\\
     =& \mbox{diag}\left( \big(\ssum{k\in \mathcal{N}} \vec{E}_k\tran \vec{E}_k\big)^{-1},\big(\ssum{k\in \mathcal{N}} \vec{E}_k\tran \vec{E}_k\big)^{-1}\right) \times \notag \\
     & \hspace{1cm} \bigg(\ssum{k\in \mathcal{N}} \vec{E}_k \tran \BlueText{\vec{v}_k^{\textrm{\tiny{Re}}(n+1)}}, \ssum{k\in \mathcal{N}}\vec{E}_k \tran \BlueText{\vec{v}_k^{\textrm{\tiny{Im}}(n+1)} } \bigg), \label{eq:netvariable_3} \displaybreak[1]\\
     =& \bigg(\underbrace{\mbox{diag}\bigg(\frac{1}{|\mathcal{N}_1|},\ldots,\frac{1}{|\mathcal{N}_N|}\bigg)\ssum{k\in \mathcal{N}} \vec{E}_k \tran  \BlueText{\vec{v}_k^{\textrm{\tiny{Re}}(n+1)}}}_{\vec{v}^{\textrm{\tiny{Re}}(n+1)}}, \notag \\
   & \hspace{0.2cm} \underbrace{\mbox{diag}\bigg(\frac{1}{|\mathcal{N}_1|},\ldots,\frac{1}{|\mathcal{N}_N|}\bigg)\ssum{k\in \mathcal{N}} \vec{E}_k \tran  \BlueText{\vec{v}_k^{\textrm{\tiny{Im}}(n+1)}}}_{\vec{v}^{\textrm{\tiny{Im}}(n+1)}}\bigg), \label{eq:netvariable_4}
 \end{align}
where \eqref{eq:netvariable_1} follows trivially from the differentiation of the objective function of problem~\eqref{eq:ADMM_net_variable_update}, \eqref{eq:netvariable_2} follows by invoking the optimality conditions for problem~\eqref{eq:ADMM_net_variable_update}, i.e., $\sum_{k\in \mathcal{N}}   \bec{E}_k\tran   \vec{y}_k=0$, \eqref{eq:netvariable_3} follows from that $\bec{E}_k=\diag(\vec{E}_k, \vec{E}_k)$, and \eqref{eq:netvariable_4} follows from that $\ssum{k\in \mathcal{N}} \vec{E}_k\tran \vec{E}_k=\diag({\BlueText{|\mathcal{N}_1|},\ldots,|\mathcal{N}_N|})$. From~\eqref{eq:netvariable_4}, it is not difficult to see that any net variable component update is equivalently obtained by averaging its copies maintained among the neighbor nodes. Such an averaging can be accomplished by using fully distributed algorithms such as gossiping~\cite{Boyd-Ran-Gos-Alg-2006}. Therefore, step~3 of ADMM-DOPF algorithm can be  carried out in a fully distributed manner.

The dual variable update~\eqref{eq:ADMM_dual_variable_update} can be carried out in a fully distributed manner, where every bus increment the current dual variables by a (scaled) discrepancy between current net variables and its own copies of those net variables.

 \section{Properties of the distributed solution method} \label{sec:algorithm_properties}

\BlueText{
Recall that the original problem~\eqref{eq:main_formulation} or equivalently problem~\eqref{eq:dist_formulation} is nonconvex and NP-hard. Therefore, ADMM based approaches are not guaranteed to converge~\cite[\S~9]{Boyd:2011:DOS:2185815.2185816}, though general convergence results are available for the \emph{convex} case~\cite[\S~3.2]{Boyd:2011:DOS:2185815.2185816}. }
Nevertheless, in the sequel, we highlight some of the convergence properties of our proposed ADMM-DOPF algorithm. In particular, we first illustrate, by using an example, the possible scenarios that can be encountered by \emph{Algorithm~2}, i.e., step~2 of the ADMM-DOPF algorithm. Then we capitalized on one of the scenario, which is empirically observed to be the most dominant, in order to characterize the solutions of the ADMM-DOPF algorithm.

\subsection{Graphical illustration of \emph{Algorithm~2}}  \label{sec:alg_proper_graph_illus}
We start by focussing on the step~2, the main ingredient of ADMM-DOPF algorithm. To get insights into the subroutine (i.e., \emph{Algorithm 2}) performed at step~2, we first rely on a simple graphical interpretation. Here instead of problem~\eqref{eq:ADMM_private_variable_update_generator_bus}, we consider a small dimensional problem to built the essential ingredient of the analysis. In particular, we consider the convex objective function $f(p,x)$ in the place of \eqref{eq:ADMM_step2_obj}. Moreover, instead of the nonconvex constraints~\eqref{eq:nonlinear_equality1_sub} and \eqref{eq:nonlinear_equality2_sub} [\emph{cf} \eqref{eq:dist_formulation_power_as_a_function_of_voltage_and_current},\eqref{eq:dist_formulation_power_as_a_function_of_voltage_and_current_flow_line}], we consider the constraint
\be
p=g(x) \ ,
\ee
where $g$ is a nonconvex function, which resembles righthand side of \eqref{eq:dist_formulation_power_as_a_function_of_voltage_and_current} and \eqref{eq:dist_formulation_power_as_a_function_of_voltage_and_current_flow_line}. Finally, instead of the remaining constraints \eqref{eq:local_variables_set_sub}, \eqref{eq:affine_sub}, \eqref{eq:convex_inequality_sub}, \eqref{eq:convex_nonlinear_inequality_sub}, and~\eqref{eq:nonconvex_inequality_sub} of problem~\eqref{eq:ADMM_private_variable_update_generator_bus}, we consider the constraint
\be \label{eq:non-convex-set-Z}
(p,x)\in\mathcal{Z} \ ,
\ee
where $\mathcal{Z}$ is not a convex set~[\emph{cf}~\eqref{eq:nonconvex_inequality_sub}]. Thus the smaller dimensional problem, which resembles subproblem~\eqref{eq:ADMM_private_variable_update_generator_bus} is given by
\begin{equation} \label{eq:small_dim_problem}
\begin{array}{ll}
\mbox{minimize} & f(p,x) \\
\mbox{subject to} & p=g(x) \\
&  (p,x)\in \mathcal{Z} \ ,
\end{array}
\end{equation}
where the variables are $p\in \R$ and $x\in \R$. Recall that \emph{Algorithm~2} approximates nonconvex functions in constraints~\eqref{eq:nonlinear_equality1_sub} and \eqref{eq:nonlinear_equality2_sub} of problem~\eqref{eq:ADMM_private_variable_update_generator_bus} by using their first order Taylor's approximations [see \eqref{eq:nonlinear_equality1_sub_approx},\eqref{eq:nonlinear_equality2_sub_approx}] and the nonconvex constraint~\eqref{eq:nonconvex_inequality_sub} by using a convex constraint [see \eqref{eq:nonconvex_inequality_sub_approx}]. Particularized to the smaller dimensional problem~\eqref{eq:small_dim_problem}, the approximations pointed above equivalent to replace $g$ by its first order Taylor's approximation $\hat g$ and to approximate $\mathcal{Z}$ by some convex set $\check{\mathcal{Z}}$, where $\check{\mathcal{Z}}\subseteq {\mathcal{Z}}$. The result is the approximated subproblem given by
\begin{equation} \label{eq:small_dim_problem_approx}
\begin{array}{ll}
\mbox{minimize} & f(p,x) \\
\mbox{subject to} & p=\hat{g}(x)= g(\hat x)+ g'(\hat x)(x-\hat x)\\
&  (p,x)\in \check{\mathcal{Z}} \ ,
\end{array}
\end{equation}
where the variables are $p\in \R$ and $x\in \R$ and $\hat x$ represents the point at which the first order Taylor's approximation is made.

Let us next examine the behavior of \emph{Algorithm~2} by considering, instead of problem~\eqref{eq:the_approximated_subproblem}, the representative smaller dimensional problem~\eqref{eq:small_dim_problem_approx}. Recall that the key idea of \emph{Algorithm~2} is to iteratively refine the first order Taylor's approximations $\hat{\boldsymbol{\lambda}}_k^{\hat{\vec{z}}_k}(\vec{z}_k)$ and $\hat{\boldsymbol{\mu}}_k^{\hat{\vec{z}}_k}(\vec{z}_k)$ [see step~2,3 of \emph{Algorithm~2}], until a stopping criterion is satisfied. This behavior is analogously understood from problem~\eqref{eq:small_dim_problem_approx}, by iteratively refining the first order Taylor's approximation $\hat g$ of $g$.

\begin{figure}[t]
        \centering
        \begin{subfigure}[t]{0.24\textwidth}
                \includegraphics[width=\textwidth]{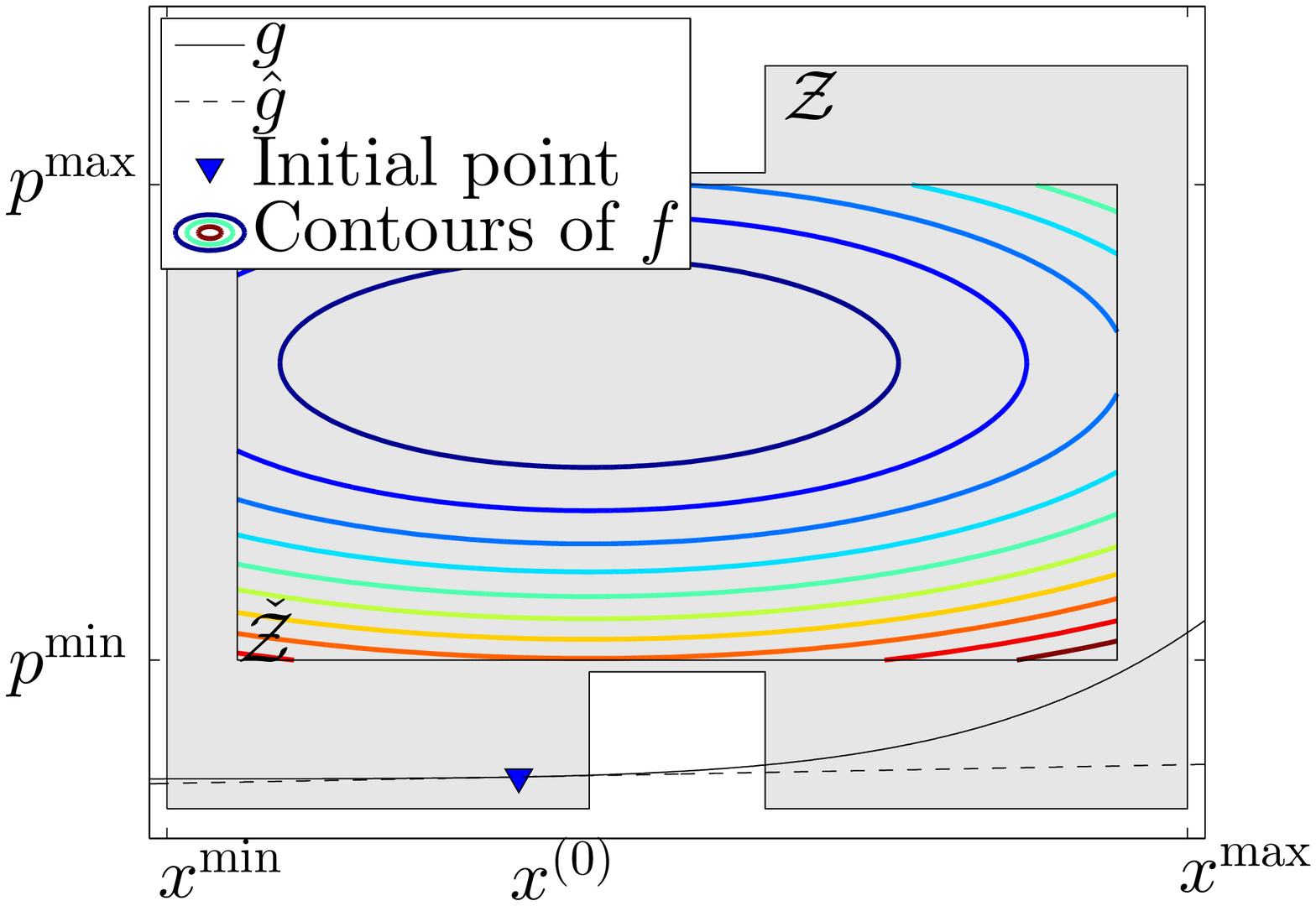}
               \caption{\subCaptionFontSize Scenario 1,  improper approximation of set~$\mathcal{Z}$, see \eqref{eq:non-convex-set-Z}, makes the approximated problem~\eqref{eq:small_dim_problem_approx} infeasible. }
                \label{fig:sub_infeasible_case}
        \end{subfigure}%
        ~
        \begin{subfigure}[t]{0.24\textwidth}
                \includegraphics[width=\textwidth]{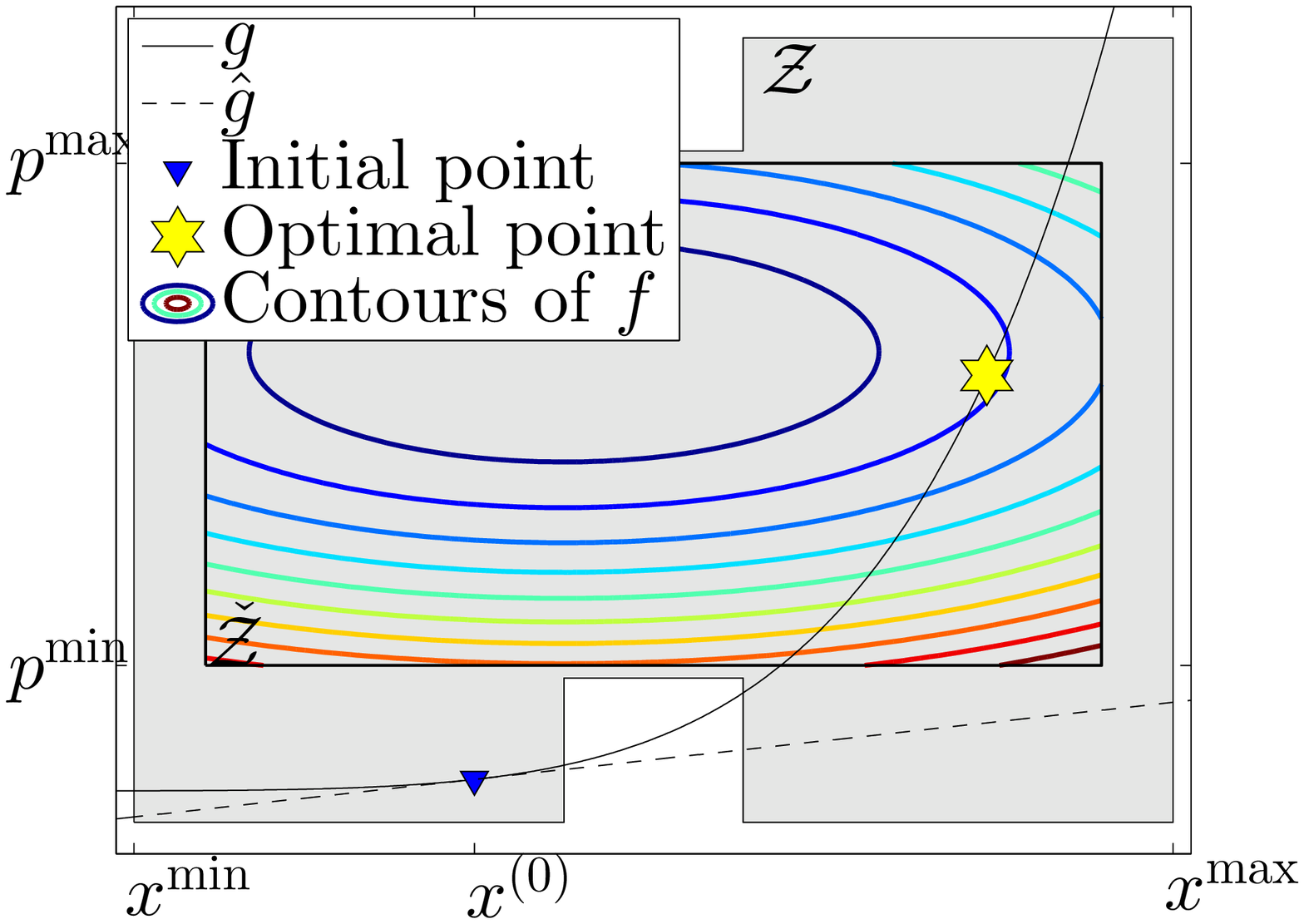}
                 \caption{\subCaptionFontSize Scenario 2, improper choice of the approximation point $\hat x=x^{(0)}$ makes the approximated problem~\eqref{eq:small_dim_problem_approx} infeasible.}
                \label{fig:subproblem_algorithm_infeasible_case}
        \end{subfigure}
        \begin{subfigure}[H]{0.24\textwidth}
                \includegraphics[width=\textwidth]{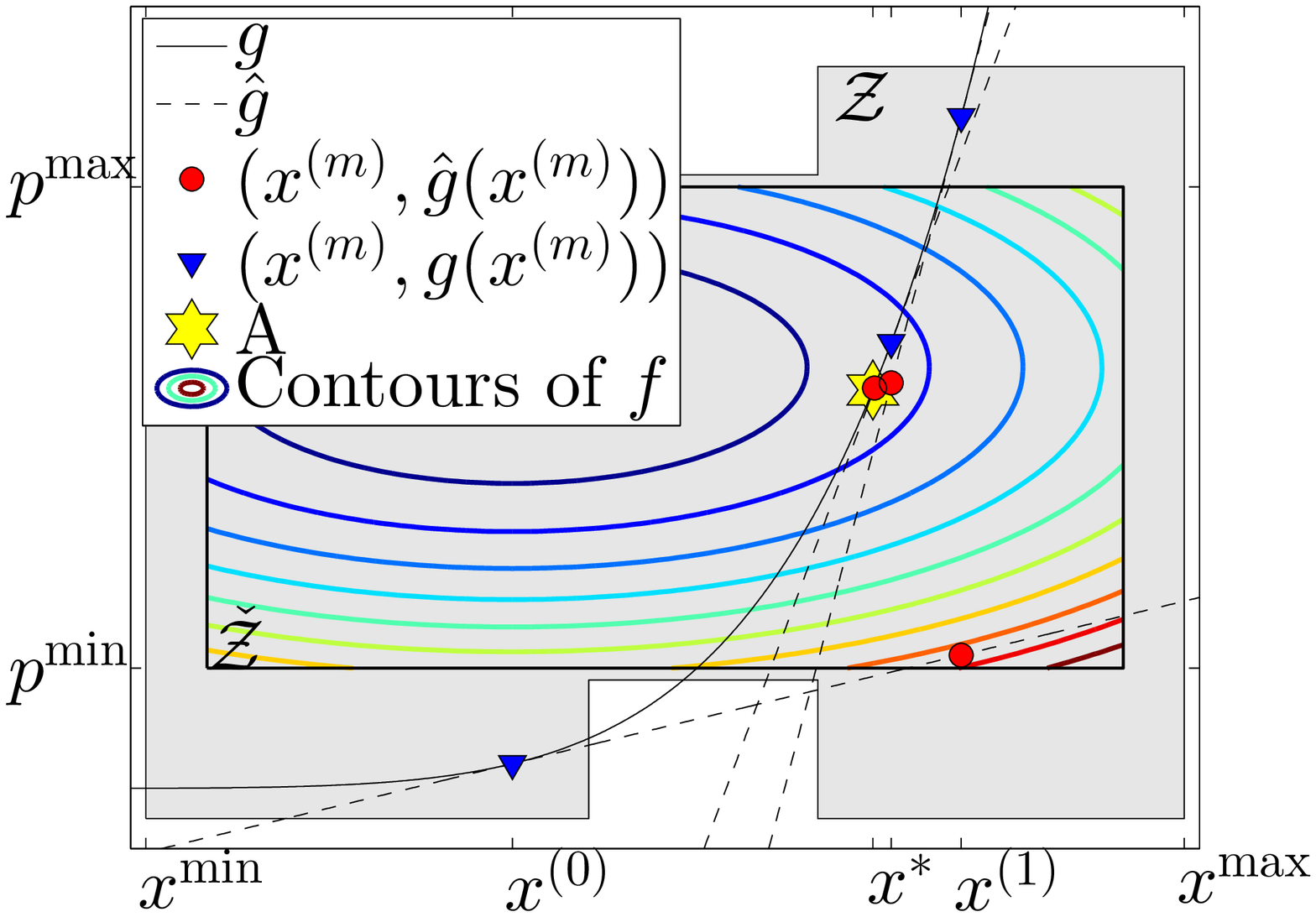}
               \caption{\subCaptionFontSize Scenario 3,  the sequence of approximations eventually converges to a desired point A. }
                \label{fig:subproblem_algorithm_converge_case}
        \end{subfigure}
        \begin{subfigure}[H]{0.24\textwidth}
                \includegraphics[width=\textwidth]{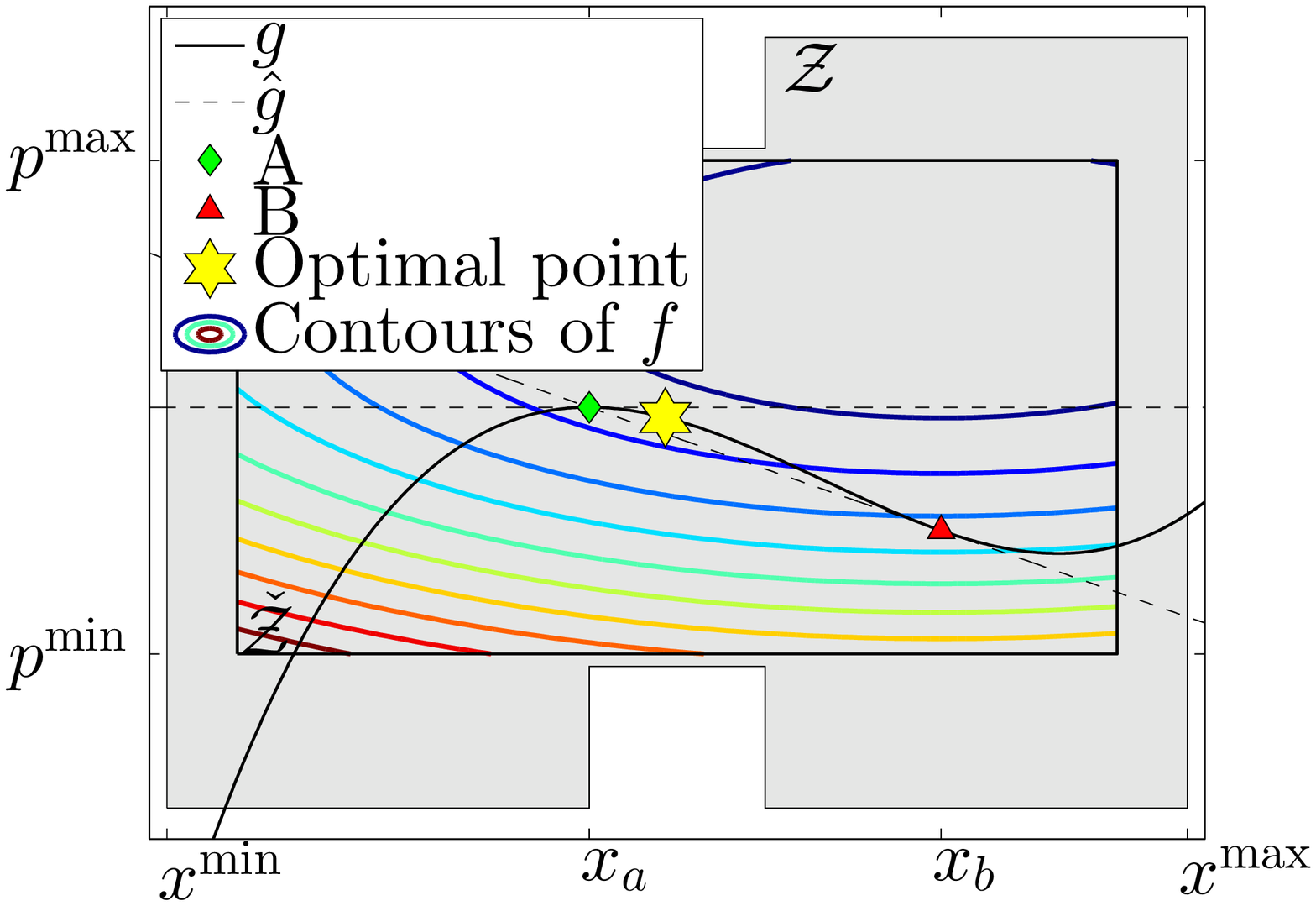}
                \caption{\subCaptionFontSize Scenario 4, the algorithm jumps between points $A=(x_a,g(x_a))$ and $B=(x_b,g(x_b))$.}
                \label{fig:subproblem_algorithm_diverges_case}
        \end{subfigure}
               \caption{ \captionFontSize Graphical illustration of \emph{Algorithm~2}.}
                \label{fig:subproblem_algorithm}
         \vspace{-5mm}
\end{figure}

Fig.~\ref{fig:subproblem_algorithm} illustrate sequential refinement of $\hat g$, where the shaded area represents the set~$\mathcal{Z}$, rectangular box represents the convex set~$\check{\mathcal{Z}}$, the solid curve represent function $g$, dotted curves represents sequential approximations $\hat g$, and thick solid curves represents the contours of $f$. Note that there are several interesting scenarios, which deserve attention to built intuitively the behavior of \emph{Algorithm~2}, see  Fig.~\ref{fig:sub_infeasible_case}- \ref{fig:subproblem_algorithm_diverges_case}.
Fig.~\ref{fig:sub_infeasible_case} shows the first scenario, where \emph{an improper approximation of set~$\mathcal{Z}$} makes the approximated problem~\eqref{eq:small_dim_problem_approx} infeasible, irrespective of the choice of $\hat x$. In contrast, Fig.~\ref{fig:subproblem_algorithm_infeasible_case} depict a scenario, where \emph{an improper choice of the approximation point $\hat x$} makes the approximated problem~\eqref{eq:small_dim_problem_approx} infeasible. Fig.~\ref{fig:subproblem_algorithm_converge_case} shows a sequence of approximations, which eventually converges to a desired point `A' that would have obtained without having the first order Taylor's approximations on $g$. Finally, Fig.~\ref{fig:subproblem_algorithm_diverges_case} shows a scenario, where a sequence of approximations switch between two points `A' and `B', i.e., there is no convergence. Any other scenario can be constructed by combining cases from Fig.~\ref{fig:sub_infeasible_case}- \ref{fig:subproblem_algorithm_diverges_case}.

Analogously, the discussion above suggests that the approximation points $({\hat{\vec{z}}_k})_{k\in\mathcal{N}}$ [\emph{cf} $\hat x$] used when constructing $\hat{\boldsymbol{\lambda}}_k^{\hat{\vec{z}}_k}(\vec{z}_k)$ and $\hat{\boldsymbol{\mu}}_k^{\hat{\vec{z}}_k}(\vec{z}_k)$ [\emph{cf} $\hat g (x)$] and the approximations used in the set $(\check{\mathcal{X}}^k_r)_{r=1,\ldots,|\mathcal{N}_k|}$ [\emph{cf} $\check{\mathcal{Z}}$], can heavily influence the performance of \emph{Algorithm~2}. Therefore, especially if the scenario~1 and 2 depicted in Fig.~\ref{fig:sub_infeasible_case}  and Fig.~\ref{fig:subproblem_algorithm_infeasible_case} occurs, during the algorithm iterations, they have to be avoided by changing the initializations. However, extensive numerical experiments show that there are specific choices of ${\hat{\vec{z}}_k}$ and $\check{\mathcal{X}}^k_r$ can make \emph{Algorithm~2} often converge to a point as depicted in Fig.~\ref{fig:subproblem_algorithm_converge_case} and barely encounters the scenarios depicted in Fig.~\ref{fig:sub_infeasible_case}, Fig.~\ref{fig:subproblem_algorithm_infeasible_case}, and Fig.~\ref{fig:subproblem_algorithm_diverges_case} see Section~\ref{sec:results_prop_of_alg2} for details.

\subsection{Optimality properties of \emph{Algorithm~2} solution} \label{sec:alg_proper_optim_props_of_algor_2}

Results obtained in this section are based on the empirical observations~(see \S~\ref{sec:results}) that scenario~3 depicted in Fig.~\ref{fig:subproblem_algorithm_converge_case} is more dominant compared to others. In particular, we make the following assumptions.

\begin{assumption} \label{assumption_1}
For any $k\in\mathcal{N}$, there exists $(\vec{z}^\star_k,\vec{u}^\star_k)$, to which \emph{Algorithm~2} can converge. Specifically, there exists $(\vec{z}^\star_k,\vec{u}^\star_k)$, where $\lim_{m\rightarrow\infty} (\vec{z}^{(m)}_k,\vec{u}^{(m)}_k)=(\vec{z}^\star_k,\vec{u}^\star_k)$ for all $k\in\mathcal{N}$. In addition, for all $k\in\mathcal{N}$, the components $(\BlueText{\vec{v}\sr{Re}_k,\vec{v}\sr{Im}_k})$ of $\vec{z}^\star_k$, \emph{strictly} satisfy the constraint~\eqref{eq:nonconvex_inequality_sub_approx}.
\end{assumption}

Under \emph{Assumption~1}, the following assertion can be made:
\begin{prop}\label{lemma:subproblem_optimality}
Suppose \emph{Assumption~1} holds. Then the output $(\vec{z}^\star_k,\vec{u}^\star_k)$ of \emph{Algorithm~2} satisfy Karush-Kuhn-Tucker (KKT) conditions for problem~\eqref{eq:ADMM_private_variable_update_generator_bus}.
\end{prop}
\begin{IEEEproof}
   See Appendix~\ref{app:proof_of_prop_1}.
\end{IEEEproof}

Combined with our empirical observations that \emph{Algorithm~2} almost always converges to a point as depicted in Fig.~\ref{fig:subproblem_algorithm_converge_case} (i.e., \emph{Assumption~1} holds usually), \emph{Proposition~\ref{lemma:subproblem_optimality}} claims that that point is \emph{locally optimal or globally optimal}.

\subsection{Optimality properties of ADMM-DOPF solution} \label{sec:alg_props_properties_of_ADMM-DOPF}
As we already pointed out, there is no guarantee that the eventual output $\big(\vec{z}, \vec{v}, \vec{u}, \vec{y}\big)$ of ADMM-DOPF is optimal, or even feasible to the original problem~\eqref{eq:dist_formulation}, because the problem is NP-hard. However, \emph{Proposition~\ref{lemma:subproblem_optimality}}, asserts that the eventual output $\big(\vec{z}_k, \vec{u}_k\big)$ of \emph{Algorithm~2} is a KKT point for problem~\eqref{eq:ADMM_private_variable_update_generator_bus} solved at step~2 of ADMM-DOPF. One can easily relate this result to characterize properties of $(\vec{z}, \vec{u})$ of ADMM-DOPF output, as we will see later. However, the properties of the remaining output $(\vec{v}, \vec{y})$ is still to be investigated. In this section, combined with the results of \emph{Proposition~\ref{lemma:subproblem_optimality}}, we analyze the optimality properties of ADMM-DOPF output.

To quantify formally the optimality properties of ADMM-DOPF, we rely on the following definition:
\begin{defin}[$(\delta,\epsilon$)-KKT optimality]\label{def:KKT_optimality}
Consider the possibly nonconvex problem of the form
\vspace{-0mm}
\begin{equation} \label{eq:optimization_prob_general}
\begin{array}{ll}
\mbox{minimize} & f_0({\vec x}) \\
\mbox{subject to} & f_i({\vec x})\leq 0, \ i = 1,\ldots,q\\
&   h_i({\vec x}) = 0, \ i = 1,\ldots,p \\
&   r_i({\vec x}) = 0, \ i = 1,\ldots,s \ ,
\end{array}
\end{equation}
where $f_0: \R^n \rightarrow \R$ is the {objective function}, $f_i: \R^n \rightarrow \R,\ i=1,\ldots,q$ are the associated inequality constraint functions, $h_i: \R^n \rightarrow \R,\ i=1,\ldots,p$ and $r_i: \R^n \rightarrow \R,\ i=1,\ldots,s$ are the equality constraint functions, and ${\vec x}\in\R^n$ is the {optimization variable}. Moreover, let $\lambda_i$ denote the dual variable associated with constraint $f_i({\vec x})\leq 0$, $\nu_i$ and $\omega_i$ denote the dual variable associated with constraint $h_i({\vec x}) = 0$ and $r_i({\vec x}) = 0$, respectively. Then an arbitrary point $({\vec x}^\star,\lambda_1^\star,\ldots,\lambda_q^\star,\nu_1^\star,\ldots,\nu_p^\star,\omega_1^\star,\ldots,\omega_p^\star)$ is called $(\delta,\epsilon$)-KKT optimal, if
\begin{align}
      f_i({\vec x}^\star)
     &\leq  0, \ i = 1,\ldots,q \label{eq:KKT1}\\
     \textstyle h_i({\vec x}^\star) &= 0, i = 1,\ldots,q \label{eq:KKT2}\\
  \textstyle (1/s)\sum_{i=1}^{s}||r_i({\vec x}^\star)||^2_2 &= \delta \label{eq:KKT3}\\
  \lambda_i^\star & \geq 0, \ i = 1,\ldots,q \label{eq:KKT4}\\
   \lambda_i^\star f_i({\vec x}^\star) & = 0, \ i = 1,\ldots,q \label{eq:KKT5}
   \end{align}
   \begin{multline}\label{eq:KKT6}
   \textstyle(1/n)||\nabla_{\vec{x}} f_0({\vec x}^\star) + \sum_{i=1}^{q}\lambda_i^\star \nabla_{\vec{x}} f_i({\vec x}^\star)\\
    + \sum_{i=1}^{p}\nu_i^\star \nabla_{\vec{x}} h_i({\vec x}^\star) + \sum_{i=1}^{s}\omega_i^\star \nabla_{\vec{x}} r_i({\vec x}^\star)||_2^2 = \epsilon \ .
 \end{multline}
\end{defin}


Note that \eqref{eq:KKT1}-\eqref{eq:KKT6} are closely related to the well known KKT optimality criterions, see~\cite[\S~5.5.3]{convex_boyd}. It suggest that smaller the $\delta$ and $\epsilon$, better the point $({\vec x}^\star,\lambda_1^\star,\ldots,\lambda_q^\star,\nu_1^\star,\ldots,\nu_p^\star,\omega_1^\star,\ldots,\omega_p^\star)$ to its local optimality. We use \emph{Definition~1} to formally analyze the optimality properties of ADMM-DOPF as discussed in the sequel.

Recall that, we have used $\vec{z}=(\vec{z}_k)_{k\in\mathcal{N}}$ to denote the vector of all the local primal variables in~\eqref{eq:local_vec}, $\vec{v}=(\vec{v}\sr{Re},\vec{v}\sr{Im})$ to denote the vector of all net variables, $\vec{u}=(\vec{u}_k)_{k\in\mathcal{N}}$ to denote the dual variables associated with constraints~\eqref{eq:local_variables_set}-\eqref{eq:nonconvex_inequality}, and finally $\vec{y}$ to denote the dual variables associated with constraint~\eqref{eq:dist_formulation_consistency_constraints}.

 Let us assume that at the termination of ADMM-DOPF, the output corresponding to $\vec{v}$ and $\vec{y}$ is $\vec{v}^\star$ and $\vec{y}^\star$, respectively. The output of ADMM-DOPF corresponding to $\vec{z}$ and $\vec{u}$ are simply the output of \emph{Algorithm~2} given by $\vec{z}^\star=(\vec{z}^\star_k)_{k\in\mathcal{N}}$ and $\vec{u}^\star=(\vec{u}^\star_k)_{k\in\mathcal{N}}$. However, unlike in convex problems, in the case of problem~\eqref{eq:dist_formulation}, one cannot take as granted that the consistency constraint~\eqref{eq:dist_formulation_consistency_constraints} is satisfied (\emph{cf} \cite[\S~3.2.1]{Boyd:2011:DOS:2185815.2185816}). In particular, $||\BlueText{\vec{v}^{(n)}_k}-\bar{\vec{E}}_k\vec{v}^\star||^2_2\rightarrow 0$ does not necessarily hold when $n\rightarrow\infty$, where $k\in\mathcal{N}$ and $n$ is the ADMM-DOPF iteration index. However, appropriate choice of the penalty parameter $\rho$ in the ADMM-DOPF algorithm usually allows finding outputs, where the consistency constraints are almost satisfied with a small error floor, which is negligible in real practical implementations as we will see empirically in \S~\ref{sec:results}. For latter use, let us quantify this error floor from $\boldsymbol{\delta}_k$,~i.e.,
\be\label{eq:discrepancy_k}
\boldsymbol{\delta}_k = \BlueText{\vec{v}^\star_k}-\bar{\vec{E}}_k\vec{v}^\star, \ k\in\mathcal{N} \ .
\ee

Now we can formally establish the optimality properties of ADMM-DOPF as follows:
\begin{prop}\label{prop:KKT_sub_optimality}
Given \emph{Assumption~1} holds, the output $(\vec{z}^\star,\vec{v}^\star, \vec{u}^\star, \vec{y}^\star)$ at the termination of  ADMM-DOPF is $(a^{-1}\bar{\delta},b^{-1}\rho^2\bar{\delta})$-KKT optimal, where $\bar{\delta}=\sum_{k\in\mathcal{N}}||{\boldsymbol \delta}_k||^2_2$, $\rho$ is the penalty parameter used in the ADMM-DOPF iterations, and $a=\mbox{len}((\boldsymbol{\delta}_k)_{k\in\mathcal{N}})$, $b=\mbox{len}(\vec{z}^\star,\vec{v}^\star)$ are normalization factors.
\end{prop}
\begin{IEEEproof}
  See appendix~\ref{app:proof_of_prop_2}.
\end{IEEEproof}

We note that deriving an analytical expression of $\delta$ by using problem~\eqref{eq:dist_formulation} data is very difficult. \BlueText{However, we can numerically compute $\delta$ and $\epsilon$ given in \emph{Proposition~\ref{prop:KKT_sub_optimality}} as
\be \label{eq:numerical-delta-and-eps}
\delta=a^{-1}\bar{\delta}, \quad  \epsilon=b^{-1}\rho^2\bar{\delta} . 
\ee
}
Extensive numerical experiences show that we usually have very small values for $\delta$. For example, for all considered simulations with $\rho=10^6$ [see \S~\ref{sec:results}], we have $\delta$ on the order of $10^{-12}$ (or smaller) and $\epsilon$ on the order of $10^{-1}$ (or smaller) after 5000 ADMM-DOPF iterations.
\addtocounter{equation}{2}

\section{Numerical results } \label{sec:results}

 In this section we present numerical experiments to illustrate the proposed algorithm.
 We compare our algorithm with the branch and bound algorithm \cite{Gopalakrishnan-BandB-OPF-2012}, centralized OPF solver provided by Matpower \cite{Zimmerman-Mat-ste-sta-ope-pla-and-ana-too-for-pow-sys-rea-and-edu-2011}, and the SDP relaxation from \cite{Lavaei-2012-OPF-TPS}.
\BlueText{ In order to study the convergence properties of the algorithm we evaluate it on four examples that have a (non-zero) duality gap, see Table~\ref{table:Introduction_to_test_problems}, rows 1-4. 
 These four examples come from~\cite{Lesieutre-Exa-the-Lim-of-App-of-Sem-Pro-to-Pow-Flo-Pro-2011,Gopalakrishnan-BandB-OPF-2012} and are  obtained by making a small modification, see Table~\ref{table:Introduction_to_test_problems}, column 3.  It is worth noting that the methods based on the SDP relaxation do not apply here due to nonzero duality gap~\cite{Lesieutre-Exa-the-Lim-of-App-of-Sem-Pro-to-Pow-Flo-Pro-2011}.
Moreover, to study the scalability properties of the proposed algorithm, we also evaluate it on two larger examples, see Table~\ref{table:Introduction_to_test_problems}, rows 5-6.      
  The exact specifications of considered examples, such as the objective functions and admittance matrices, are found in Table~\ref{table:Introduction_to_test_problems} and references therein. 
 \begin{table}[h]\scriptsize 
  \centering

  \begin{tabular}{|c|c|c|c|c|c| } 
            \hline
          $N$  &  Original Problem &  Modification  &  \# G    & \# L & FL   \\   \hline
         3   & 3bus~\cite{Lesieutre-Exa-the-Lim-of-App-of-Sem-Pro-to-Pow-Flo-Pro-2011}    &  $s_{23}\sr{max}=s_{32}\sr{max}=50$ MVA    & 3 & 3 &   \eqref{eq:main_formulation_power_limit_flow_line} \\  \hline
         9   & Case9~\cite{Zimmerman-Mat-ste-sta-ope-pla-and-ana-too-for-pow-sys-rea-and-edu-2011} &   $q_k\sr{G,min}=10$ MVAr, $\forall k\in \mathcal{G}$    & 3     & 3 & None    \\   
             &                                                                     &   $p_k\sr{D}=1.1 \bar{p}_k\sr{D}$,  $\forall k\in \mathcal{N}$             & & &     \\   \hline
        14  & IEEE14~\cite{IEEE_test_networks}          &    $q_k\sr{G,min}=0$ MVAr, $\forall k\in \mathcal{G}$      & 5  & 11 &  None   \\   
              &                                                                    &    ${p}\sr{D}_k=0.1 \bar{p}\sr{D}_k,\forall k\in \mathcal{N}$         &   & &    \\   \hline
        30 & IEEE30~\cite{IEEE_test_networks}           &  ${p}\sr{D}_k=0.5 \bar{p}_k\sr{D},_k\forall k\in \mathcal{N}$  &   6     & 21& None      \\ 
             &                                                                    &  ${q}\sr{D}_k=0.1 \bar{q}\sr{D}_k,\forall k\in \mathcal{N}$     &                & &      \\  \hline
         118    &  IEEE118~\cite{IEEE_test_networks}  &  None     &     54      &  99 &  None    \\  \hline
         300   &  IEEE300~\cite{IEEE_test_networks}   &  None   &        69      & 201 &  None    \\  \hline
  \end{tabular}
 \caption{\BlueText{  \captionFontSize Specifications of the test problems. The first column indicates the number of buses. The second column gives the reference to the original problems.  The third column show how we modify the original problem ($ \bar{p}_k\sr{D}$ and  $\bar{q}_k\sr{D}$ indicate  the original problem data associated with the power demands).  The fourth and fifth columns specify the number of generators and loads respectively.  The sixth column specifies the type of flow line limit (FL) used, if any, i.e., which of constraints ~\eqref{eq:main_formulation_current_limit_flow_line}, \eqref{eq:main_formulation_power_limit_flow_line}, and \eqref{eq:main_formulation_real_power_limit_flow_line}  are included. } }
  \label{table:Introduction_to_test_problems}
  \vspace{-4mm}
\end{table}
  
The units of real power, reactive power,  apparent power, voltage magnitude, and the objective function values are MW, MVAr, MVA, p.u. (per unit)~\footnote{ \BlueText {The voltages base is 400 kV.}}, and \$/hour, respectively. 
   In all 6 problems the average power demand of the loads is in the range 10-100 MW and 1-10 MVAr, for the real and reactive powers, respectively. 

}

 The simulations were executed in a sequential computational environment, using matlab version 8.1.0.604 (R2013a) \cite{MATLAB:2013}.
 The convex problem~\eqref{eq:the_approximated_subproblem} is solved with the convex solution method presented in Section \ref{sec:ald_devlop_qp_solver} together with the built in matlab QP solver quadprog.
 As a stopping criterion for \emph{Algorithm~2} we use $\epsilon\sr{sub}=10^{-10}$ and \BlueText{$\texttt{max\_iter}=20$, unless stated otherwise}.
 For the ADMM method we use $\vec{v}\sr{Re}=(1,\cdots,1)$ $\vec{v}\sr{Im}=(0,\cdots,0)$ and  $\vec{y} = (0,\cdots, 0)$ as an initial point.

  \subsection{Properties of \emph{Algorithm~2}} \label{sec:results_prop_of_alg2}

 \BlueText{  In this subsection we investigate the convergence properties of \emph{Algorithm~2}. 
                   In particular, we relate the convergence behavior of \emph{Algorithm~2} to  the analysis in Sections~\ref{sec:alg_proper_graph_illus} and~\ref{sec:alg_proper_optim_props_of_algor_2}, where 4 scenarios, or possible outcomes, of \emph{Algorithm~2} [Fig.~\ref{fig:subproblem_algorithm}] were identified.

  During the numerical evaluations, \emph{Algorithm~2} was executed 11060000 times. 
  Scenarios 1 and 2 [Fig.~\ref{fig:sub_infeasible_case} and~\ref{fig:subproblem_algorithm_infeasible_case}, respectively], where the approximated subproblem is infeasible, occurred only 6 times ($\approx 0.00008\%$). These occurrences  happened at one of the buses in the 300 bus example, during ADMM iterations 137-143. This suggests even if a particular bus fails to converge with \emph{Algorithm~2} in consecutive ADMM iterations, the bus can recover to find a solution in latter ADMM iterations.
  
To numerically study the occurrences of Scenarios 3 and 4, we run \emph{Algorithm~2} with 1) $\texttt{max\_iter}=20$, 2) $\texttt{max\_iter}=1000$, and 3) $\texttt{max\_iter}=25000$. In all considered cases, we use $\epsilon\sr{sub}=10^{-10}$ [compare with \eqref{eq:sub_problem_stopping_criterion_epsilon}]. The following table summarizes how frequently the stopping criteria~\eqref{eq:sub_problem_stopping_criterion_max_iter} of \emph{Algorithm~2} is met. 
\begin{table}[h]\scriptsize 
  \centering 
  \begin{tabular}{| c || c | c | c | c | c | } 
            \hline
             \multirow{2}{*}{$N$ }         &  \multicolumn{4}{ c |}{\texttt{max\_iter}}\\
             & \multicolumn{1}{ c }{$10$}   & \multicolumn{1}{ c }{$20$}  & \multicolumn{1}{ c }{$1000$} & \multicolumn{1}{ c |}{$25000$} \\  \hline
            $3$     & $2\times 10^{-2}\%$ & $5 \times 10^{-3} \%$   & $5 \times 10^{-3} \%$  & $5 \times 10^{-3} \%$     \\     
            $9$     & $0\%$ & $0\%$                      & $0\%$  & $0\%$   \\      
            $14$   & $2\times 10^{-1}\%$        &$0\%$ & $0\%$  & $0\%$     \\      
            $30$   & $2 \times 10^{-5}\%$ &$2 \times 10^{-5}\%$ & $2 \times 10^{-5}\%$   & $2 \times 10^{-5}\%$        \\       
            $118$ & $8 \times 10^{-5}\%$   &$8 \times 10^{-5}\%$ & $8 \times 10^{-5}\%$   & $8 \times 10^{-5}\%$ \\      
            $300$ & $4 \times 10^{-1}\%$ &$3 \times 10^{-2}\%$ & $9 \times 10^{-3}\%$  & $8 \times 10^{-3} \%$ \\      \hline
  \end{tabular}
 \caption{\captionFontSize  Frequency of the  termination of \emph{Algorithm~2}  from the stopping criteria~\eqref{eq:sub_problem_stopping_criterion_max_iter}.}
  \label{table:comparison_frequency}
  \vspace{-4mm}
\end{table}

Note that the entries of Table~\ref{table:comparison_frequency} suggest an upper bounds on the frequencies of Senario~4. Therefore, Table~\ref{table:comparison_frequency} shows that the frequencies of the Scenario~4 is decreasing (or unchanged) as $\texttt{max\_iter}$ is increased. However the effects are marginal as indicated in the table, specially for $\texttt{max\_iter}\geq 20$. On the other hand, recall that $\epsilon\sr{sub}=10^{-10}$, i.e., the decrement of voltages between two successive iterations is below $10^{-10}$ [compare with \eqref{eq:sub_problem_stopping_criterion_epsilon}].  Such an infinitesimal accuracy in the stopping criteria~\eqref{eq:sub_problem_stopping_criterion_epsilon} suggests the algorithm's convergence, see Scenario 3, Fig.~\ref{fig:subproblem_algorithm_converge_case}. For example, consider the case $N=3$ and $\texttt{max\_iter}=20$ in Table~\ref{table:comparison_frequency}. From the results, the frequency of termination of  \emph{Algorithm~2}  from the stopping criteria~\eqref{eq:sub_problem_stopping_criterion_epsilon}  (i.e., Scenario~3) become $99.995\%$.  Thus, from \emph{Proposition~\ref{lemma:subproblem_optimality}}, it follows that when $\texttt{max\_iter}=20$, $99.995\%$ of the cases \emph{Algorithm~2} converges to a point satisfying the KKT conditions. It is worth noting that for all the considered cases the convergence of the algorithm is in the range $99.99\%$-$100\%$. Results also suggest that the convergence properties of \emph{Algorithm~2}  can be improved (see the case $N=300$) or remain intact  (see the cases $N=3$,  $9$, $14$, $30$, and $118$)  at the expense of the increase in $\texttt{max\_iter}$.

 }

   \begin{figure}[t]
      \centering
        \includegraphics[width=0.25\textwidth]{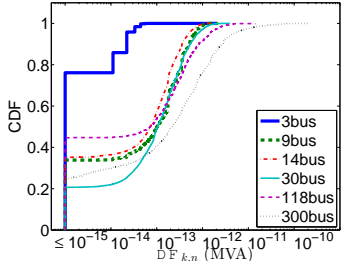}
        \caption{\captionFontSize CDF displaying $\textsc{DF}_{k,n}$ of Eq. \ref{eq:degree_of_feasiblility} for every subproblem $k$ and every ADMM iteration $n$ for each of the four examples.}
        \label{fig:subproblem_CDF}
       \vspace{-5mm}
   \end{figure}

   \begin{figure}[t]
        \centering
        \begin{subfigure}[t]{0.24\textwidth}
                \includegraphics[width=\textwidth]{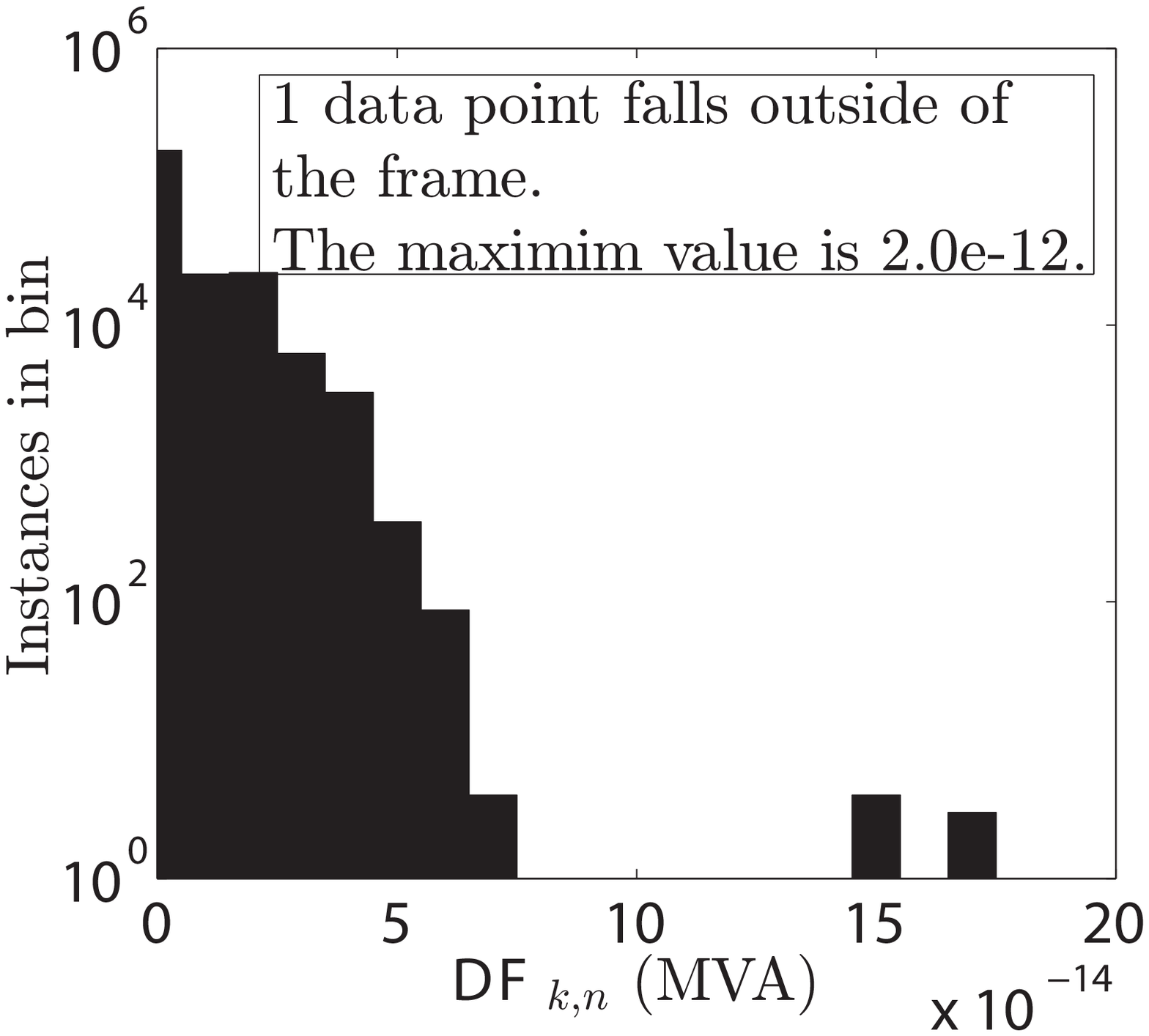}
                \caption{\subCaptionFontSize 3 bus networks }
                \label{fig:subproblem_histogram_3bus}
        \end{subfigure}%
        ~
        \begin{subfigure}[t]{0.24\textwidth}
                \includegraphics[width=\textwidth]{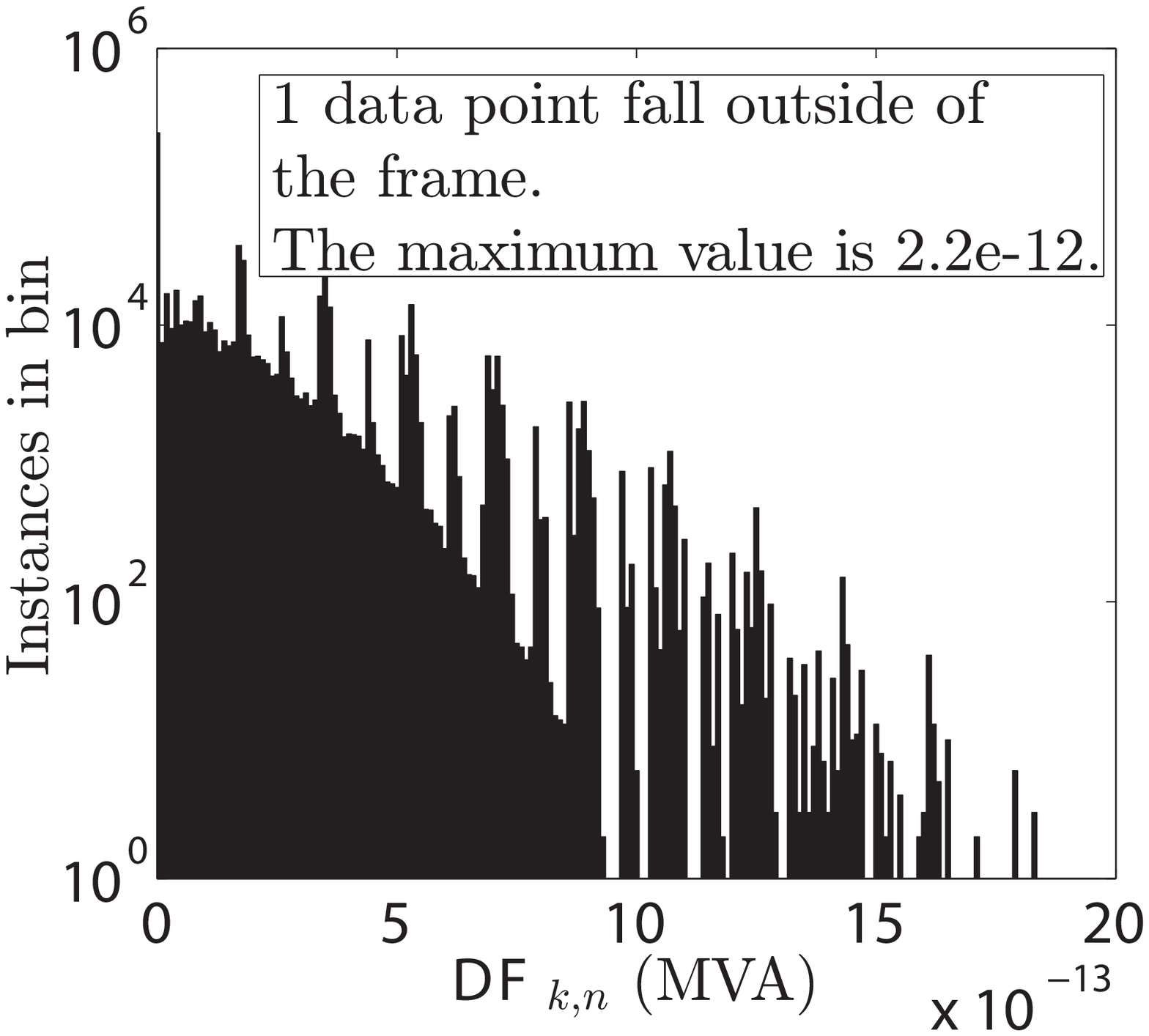}
                \caption{\subCaptionFontSize 9 bus network }
                \label{fig:subproblem_histogram_9bus}
        \end{subfigure}
        \begin{subfigure}[H]{0.24\textwidth}
                \includegraphics[width=\textwidth]{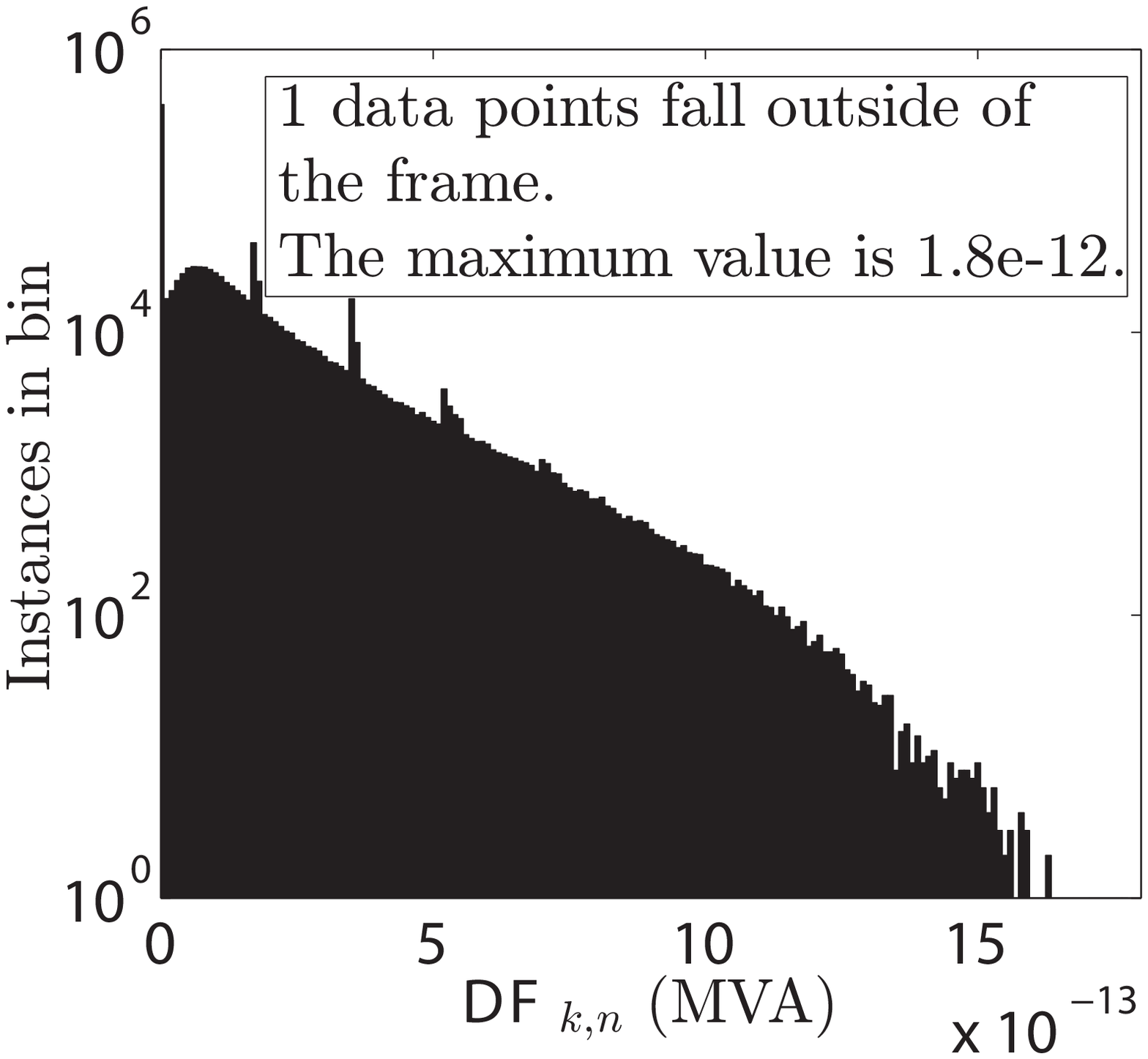}
                \caption{\subCaptionFontSize 14 bus network}
                \label{fig:subproblem_histogram_14bus}
        \end{subfigure}
        \begin{subfigure}[H]{0.24\textwidth}
                \includegraphics[width=\textwidth]{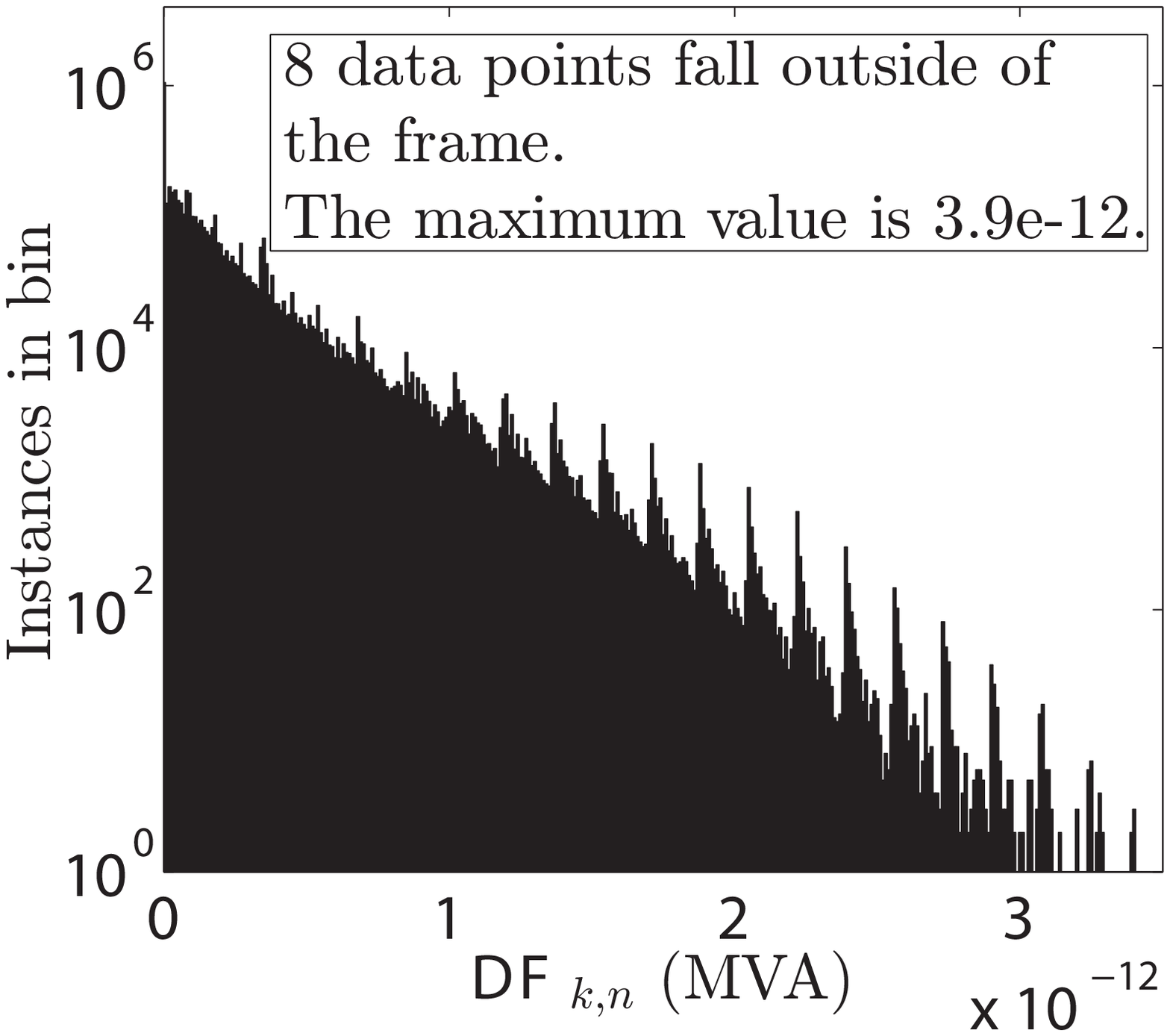}
                \caption{\subCaptionFontSize 30 bus network}
                \label{fig:subproblem_histogram_30bus}
        \end{subfigure}
        \begin{subfigure}[H]{0.24\textwidth}
                \includegraphics[width=\textwidth]{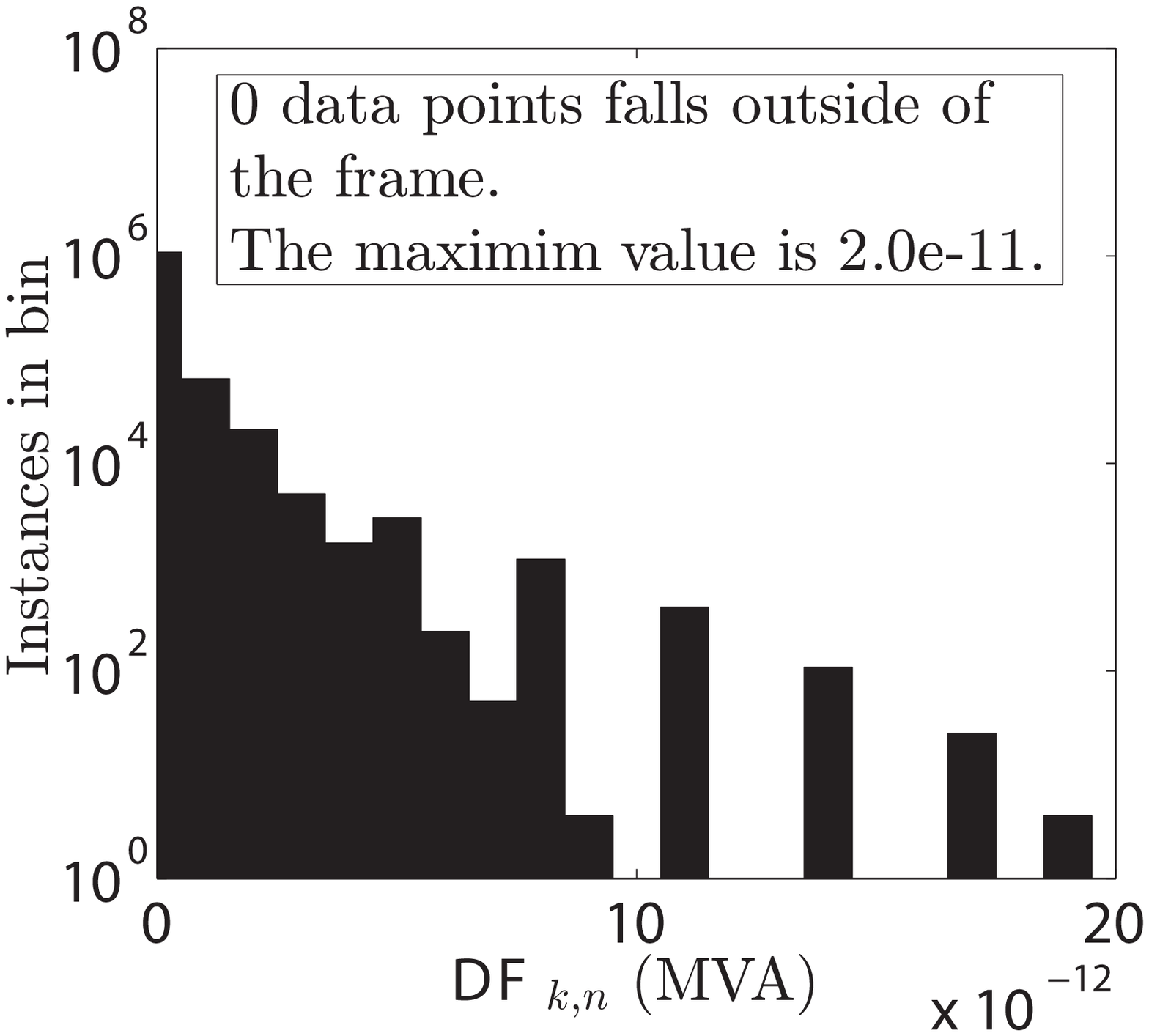}
                \caption{\subCaptionFontSize 118 bus network}
                \label{fig:subproblem_histogram_118bus}
        \end{subfigure}
        \begin{subfigure}[H]{0.24\textwidth}
                \includegraphics[width=\textwidth]{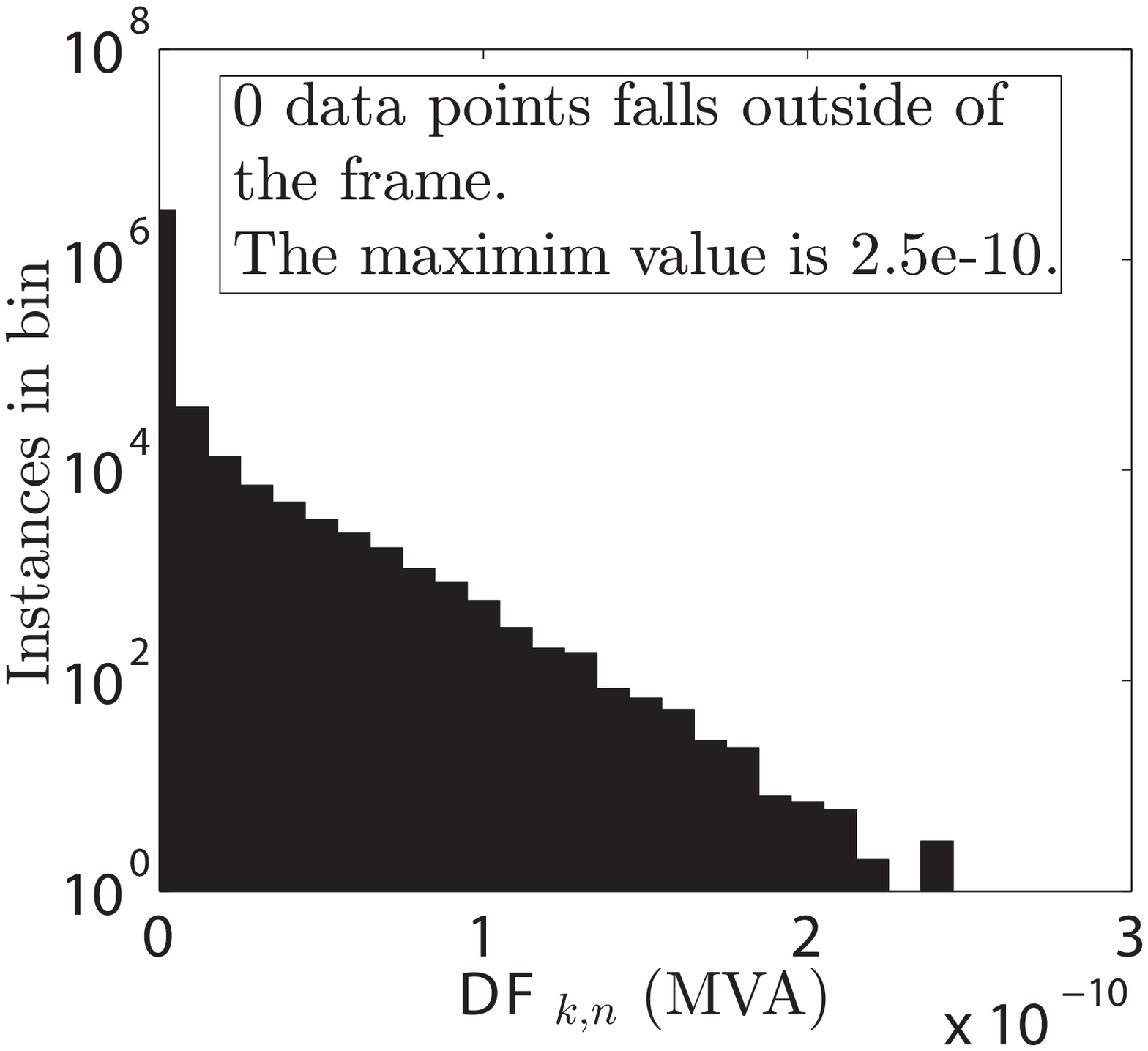}
                \caption{\subCaptionFontSize 300 bus network}
                \label{fig:subproblem_histogram_300bus}
        \end{subfigure}
               \caption{\captionFontSize Histograms displaying $\textsc{DF}_{k,n}$ in Eq. \ref{eq:degree_of_feasiblility} for every subproblem $k$ and every ADMM iteration $n$ for \BlueText{the considered test networks}.}
                \label{fig:subproblem_histogram}
              \vspace{-5mm}
  \end{figure}

\BlueText{
  Note that the voltages returned by \emph{Algorithm~2} are always feasible to problem~\eqref{eq:ADMM_private_variable_update_generator_bus}, i.e., satisfies~\eqref{eq:nonconvex_inequality_sub}. 
 However, the resulting power injections might be infeasible, [compare with~\eqref{eq:affine_sub}-\eqref{eq:convex_nonlinear_inequality_sub}]. 
 To measure the feasibility of the returned power injections of \emph{Algorithm~2}, we define the following metric called the degree of feasibility (DF):
}
  \begin{equation} \label{eq:degree_of_feasiblility}
         \textsc{DF}_{k,n} =  \min_{p + j q\in \mathcal{S}_k} | p_k^{(n)} + j q_k^{(n)}-(p + j q) |  ,
  \end{equation}
  where $k$ and $n$ indicate the bus and ADMM iteration, respectively, $p_k^{(n)} + j q_k^{(n)}$ is the returned power injection, and
  \begin{equation}
       \mathcal{S}_k =  \left\{ z\in \C \left|  \begin{array}{c}
                                                                            p_k\sr{G,min}-p_k\sr{D} \leq \Re(z) \leq  p_k\sr{G,max}-p_k\sr{D} \\
                                                                            q_k\sr{G,min}-q_k\sr{D} \leq \Im(z) \leq  q_k\sr{G,max}-q_k\sr{D}
                                                                      \end{array} \right\} \right..                         
  \end{equation}
 \BlueText{  The unit of measurement for $\textsc{DF}_{k,n}$ is MVA. }
 In order to provide a statistical description of $\textsc{DF}_{k,n}$ for every execution of \emph{Algorithm~2}, we consider an empirical cumulative distribution function (CDF) [Fig.~\ref{fig:subproblem_CDF}] and a histogram [Fig.~\ref{fig:subproblem_histogram}], for each example separately. 
 These results suggest that \emph{Algorithm~2} returns a feasible solution with high accuracy in all cases, where the worst case accuracy is \BlueText{$2.5\times 10^{-10}$}.

 As a consequence of this promising behavior of \emph{Algorithm~2}, we will proceed under \emph{Assumption~\ref{assumption_1}}.



  \subsection{Connection to \emph{Proposition~\ref{prop:KKT_sub_optimality}} } \label{subsection:num_con_to_prop_2}
      In this section we relate the numerical evaluations to \emph{Proposition~\ref{prop:KKT_sub_optimality}}.
\BlueText{      In particular, we inspect the behavior of $\delta$ [compare with~\eqref{eq:KKT3} and~\eqref{eq:numerical-delta-and-eps}],  and $\epsilon$ [compare with~\eqref{eq:KKT6} and~\eqref{eq:numerical-delta-and-eps}] with respect to $\rho$, which are defined in Section~\ref{sec:alg_props_properties_of_ADMM-DOPF}.  }
\BlueText{
     The unit of measurements for $\delta$ is $p.u.^2$ and
     the unit of $\epsilon$ can be interpreted as the square of the decrease/increase in \$/hour with respect to a small perturbation in the variable $\vec{z}=(\vec{z}_k)_{k\in \mathcal{N}}$.
     }



\begin{figure}
        \centering
        \begin{subfigure}[t]{0.24\textwidth}
                \includegraphics[width=\textwidth]{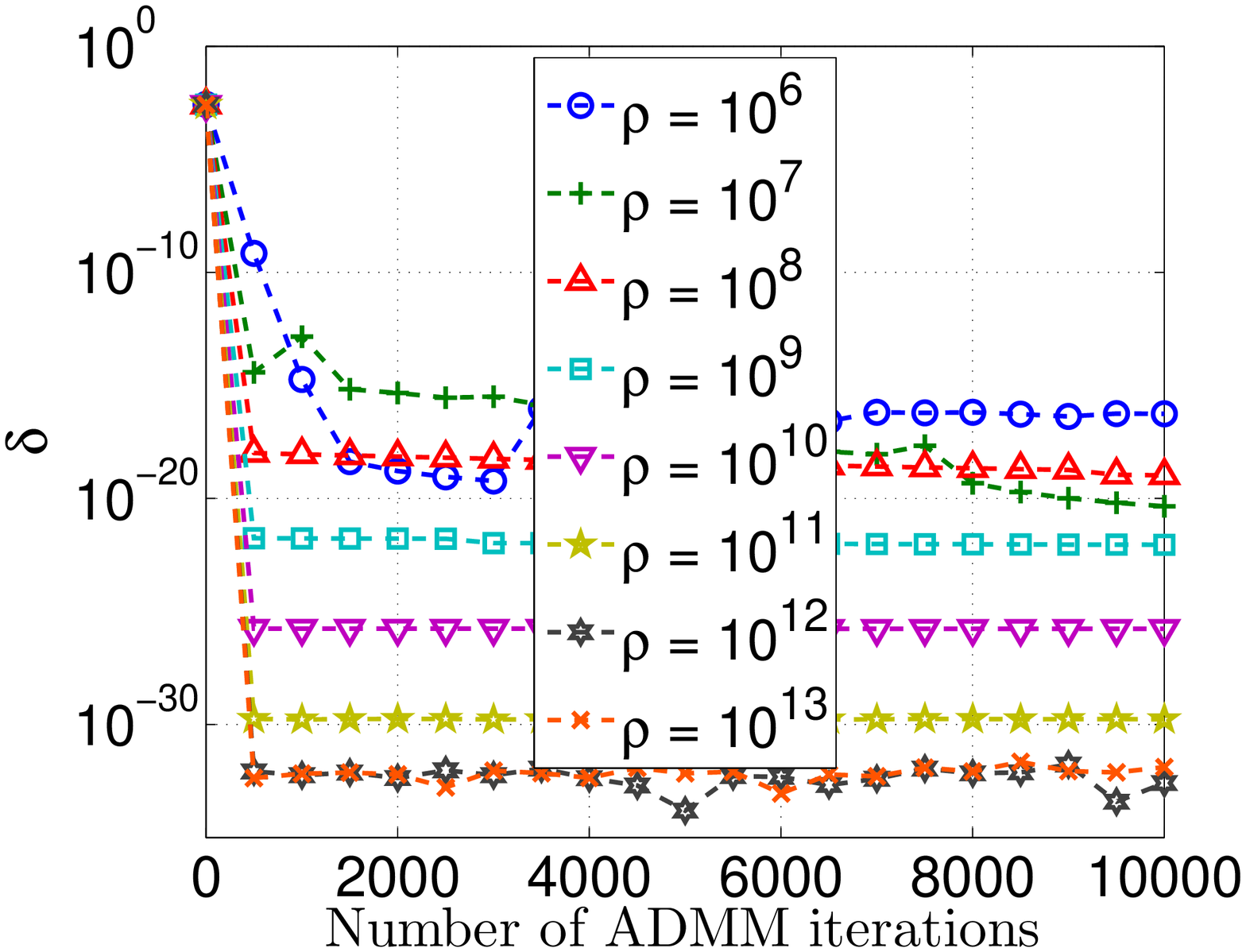}
                \caption{\subCaptionFontSize $3$ bus network }
                \label{fig:ADMM_iter_vs_cons_3bus}
        \end{subfigure}%
        ~
        \begin{subfigure}[t]{0.24\textwidth}
                \includegraphics[width=\textwidth]{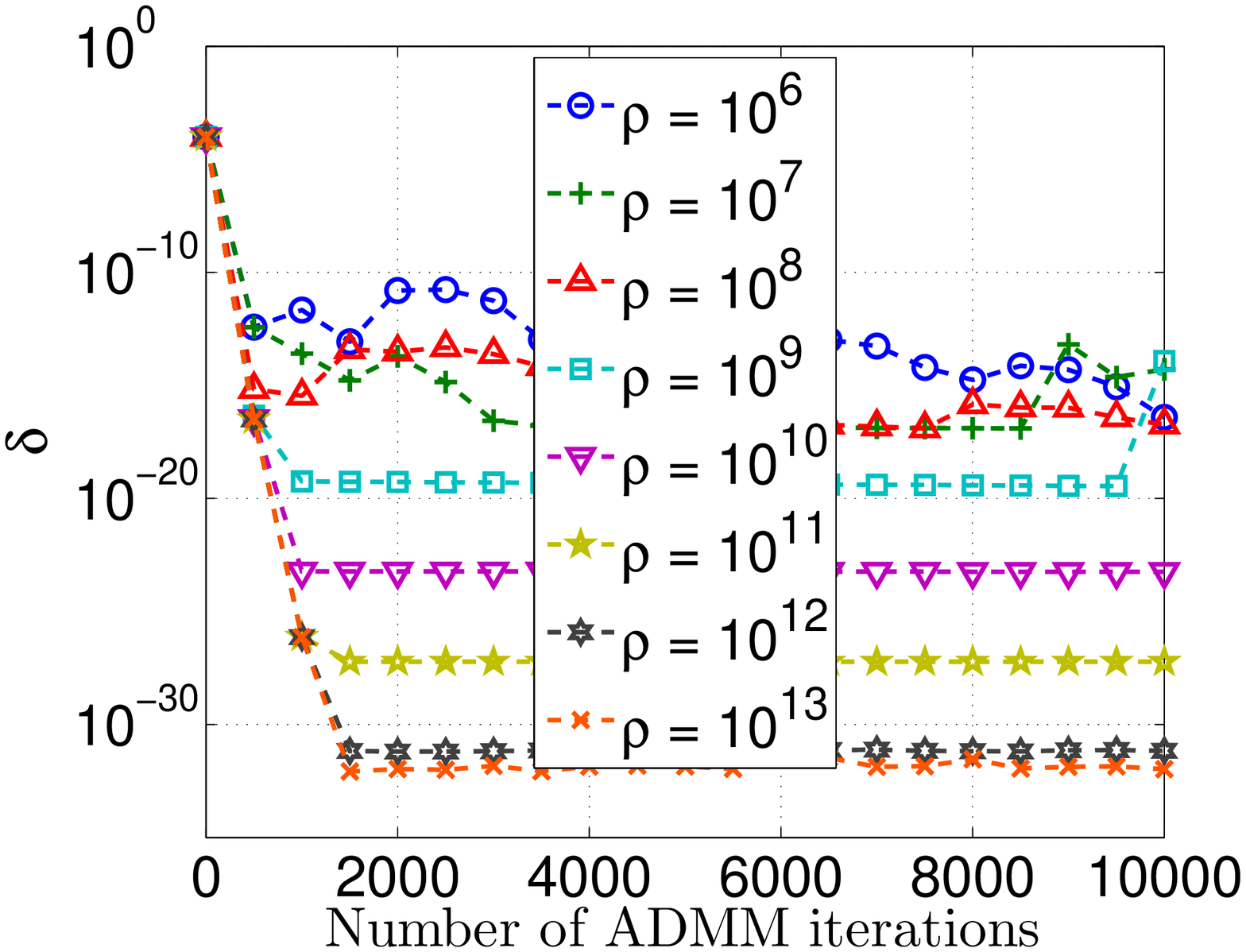}
                \caption{\subCaptionFontSize $9$ bus network}
                \label{fig:ADMM_iter_vs_cons_9bus}
        \end{subfigure}
        \begin{subfigure}[H]{0.24\textwidth}
                \includegraphics[width=\textwidth]{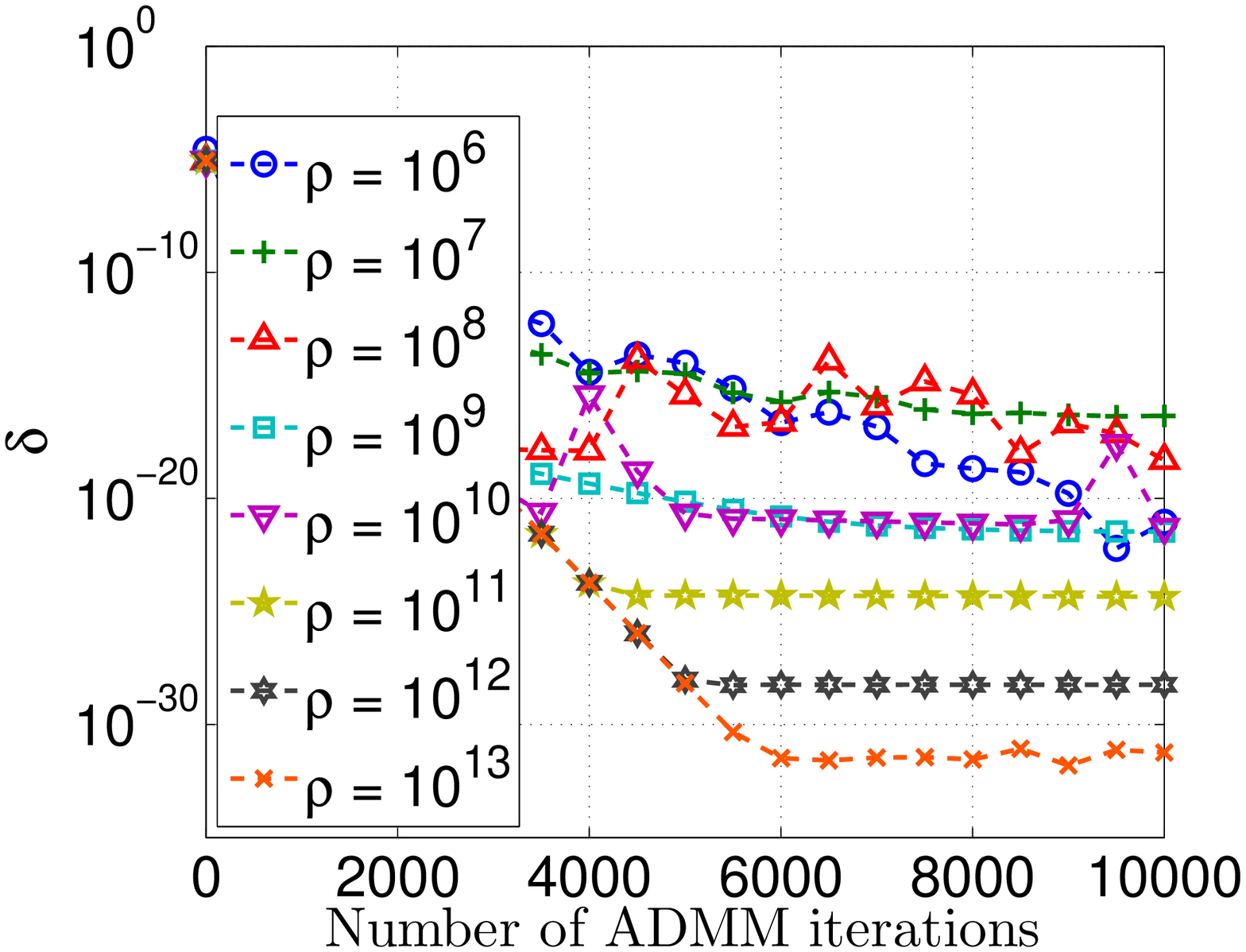}
                \caption{\subCaptionFontSize $14$ bus network}
                \label{fig:ADMM_iter_vs_cons_14bus}
        \end{subfigure}
        \begin{subfigure}[H]{0.24\textwidth}
                \includegraphics[width=\textwidth]{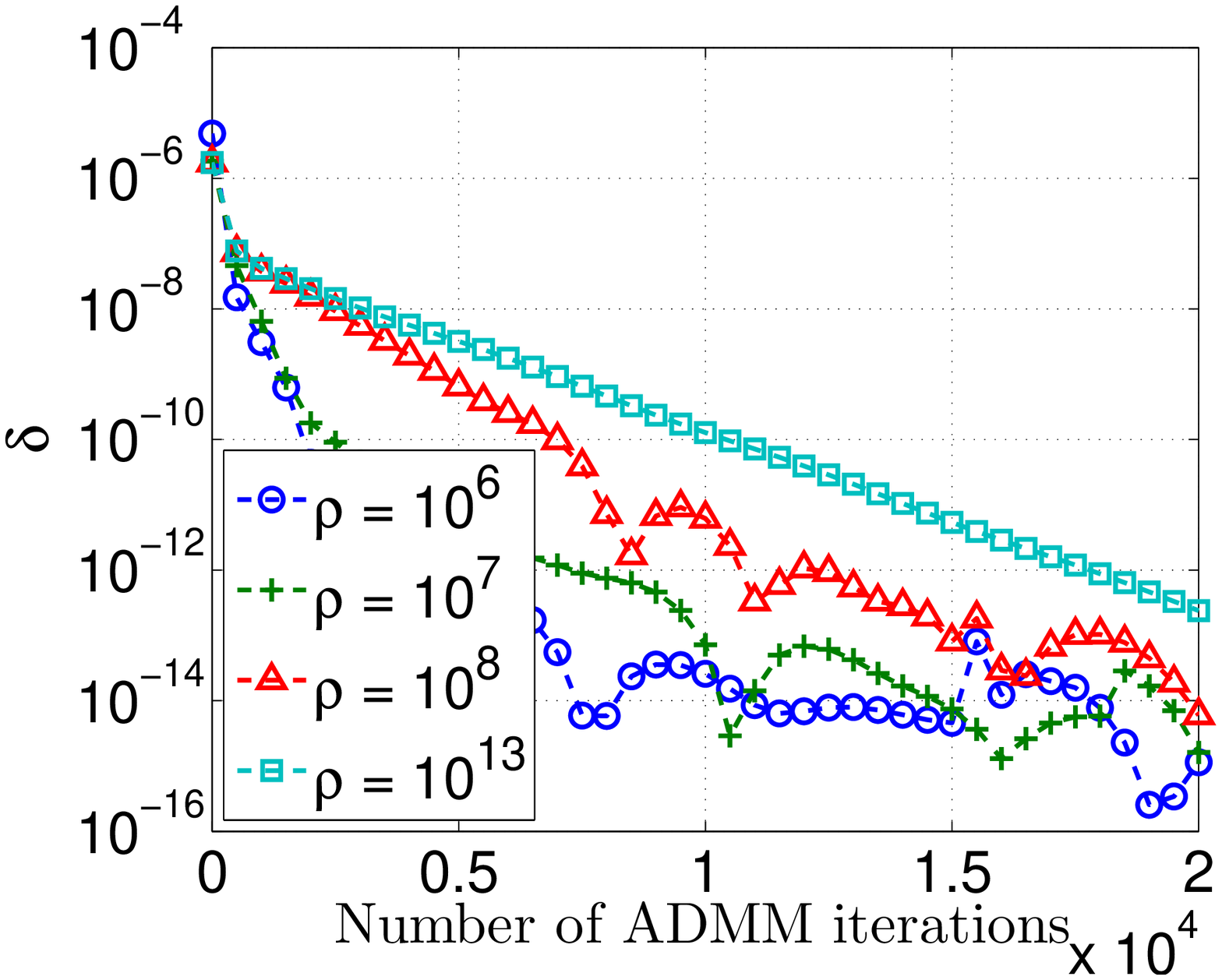}
                \caption{\subCaptionFontSize $30$ bus network}
                \label{fig:ADMM_iter_vs_cons_30bus}
        \end{subfigure}
               \caption{\captionFontSize $\delta$ versus number of ADMM iterations.   }
                \label{fig:ADMM_iter_vs_cons}
            \vspace{-5mm}
  \end{figure}

 Fig.~\ref{fig:ADMM_iter_vs_cons} depicts $\delta$ at every  $500$ ADMM iterations, for $\rho=10^6,\cdots,10^{13}$.
 In the $30$ bus example the results  are almost identical for $\rho=10^9,10^{10},\cdots,10^{13}$ and, accordingly, we only include the results for $\rho=10^6,10^7,10^8,10^{13}$. 
 Since $\delta$ measures the inconsistency between the subproblems, the point returned by ADMM-DOPF can only be considered feasible when $\delta$ has reached acceptable accuracy, i.e., $\delta< \gamma$ for some $\gamma>0$.
 We do not consider any particular threshold $\gamma$, since we are only interested in observing the convergence behavior.
 In this aspect, the result show a promising behavior, as $\delta$ has a decreasing trend in all cases. 
 Furthermore, for the $3$, $9$, and $14$ bus examples, $\delta$ converges to a fixed error floor for the larger values of the penalty parameter $\rho$.
 In particular, as $\rho$ increases, $\delta$ converges to a point closer to zero, which suggests a negative relationship between $\delta$ and $\rho$.
 Thereby, indicating that increasing the penalty parameter enforces higher accuracy of consistency among the subproblems.
 On the contrary to the $3$, $9$, and $14$ bus examples, $\delta$ decreases slower when the penalty parameter increases in the case of the $30$ bus example.
 However, in the case of the $30$ bus example, $\delta$ is still decreasing after the last iteration considered when $\rho=10^9,\cdots,10^{13}$.

\begin{figure}
        \centering
        \begin{subfigure}[t]{0.24\textwidth}
                \includegraphics[width=\textwidth]{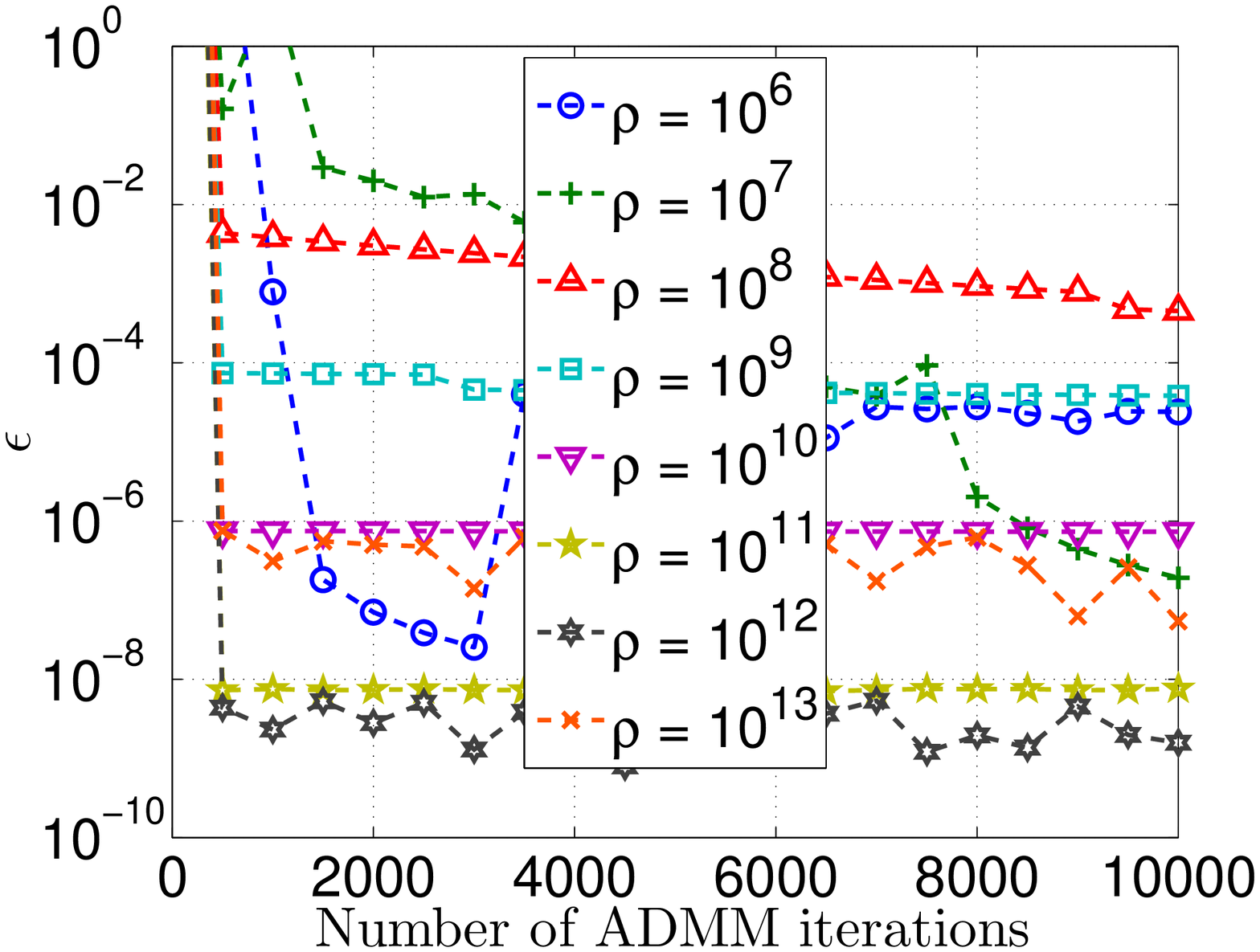}
                \caption{\subCaptionFontSize $3$ bus network }
                \label{fig:ADMM_iter_vs_cons_3bus}
        \end{subfigure}%
        ~
        \begin{subfigure}[t]{0.24\textwidth}
                \includegraphics[width=\textwidth]{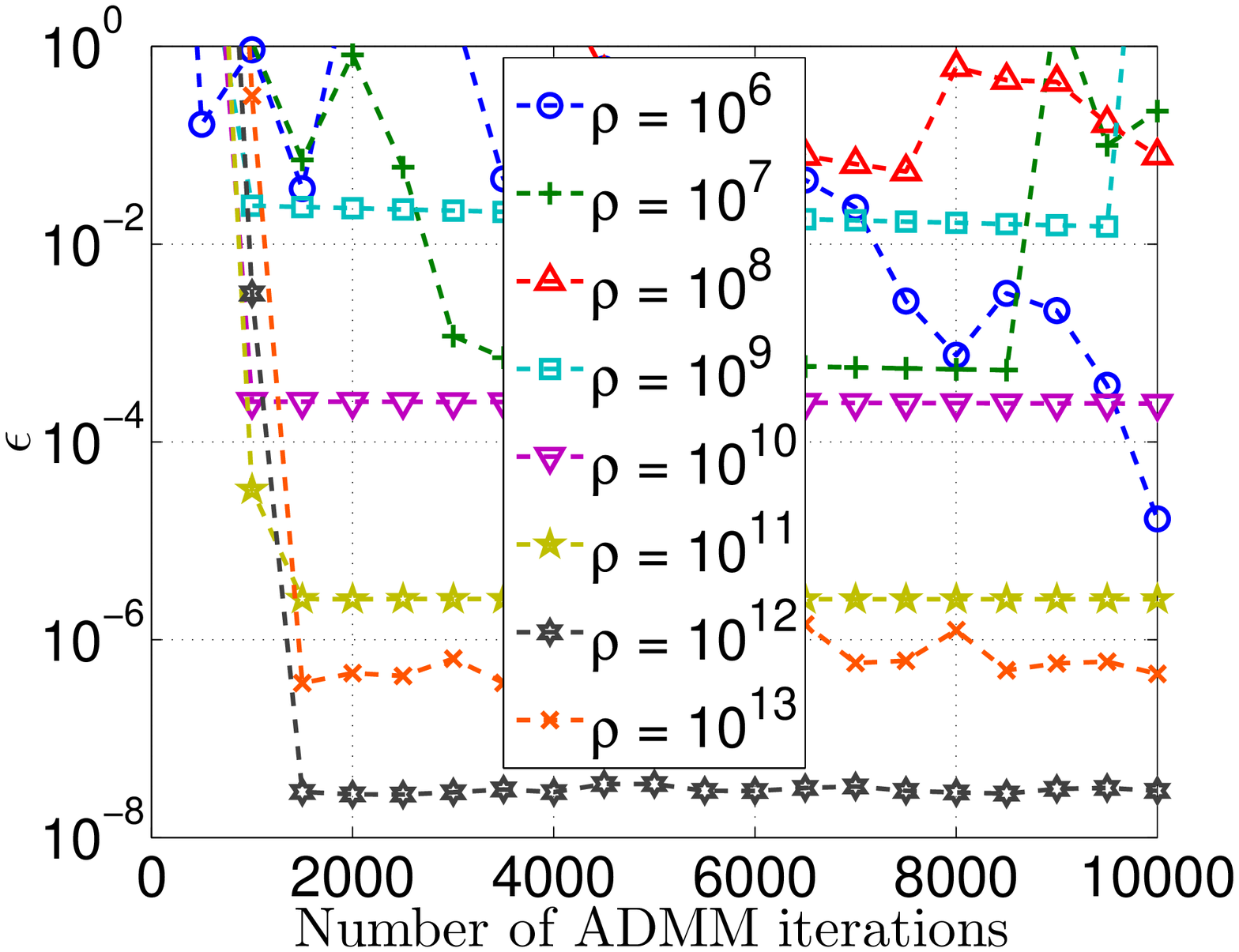}
                \caption{\subCaptionFontSize $9$ bus network}
                \label{fig:ADMM_iter_vs_cons_9bus}
        \end{subfigure}
        \begin{subfigure}[H]{0.24\textwidth}
                \includegraphics[width=\textwidth]{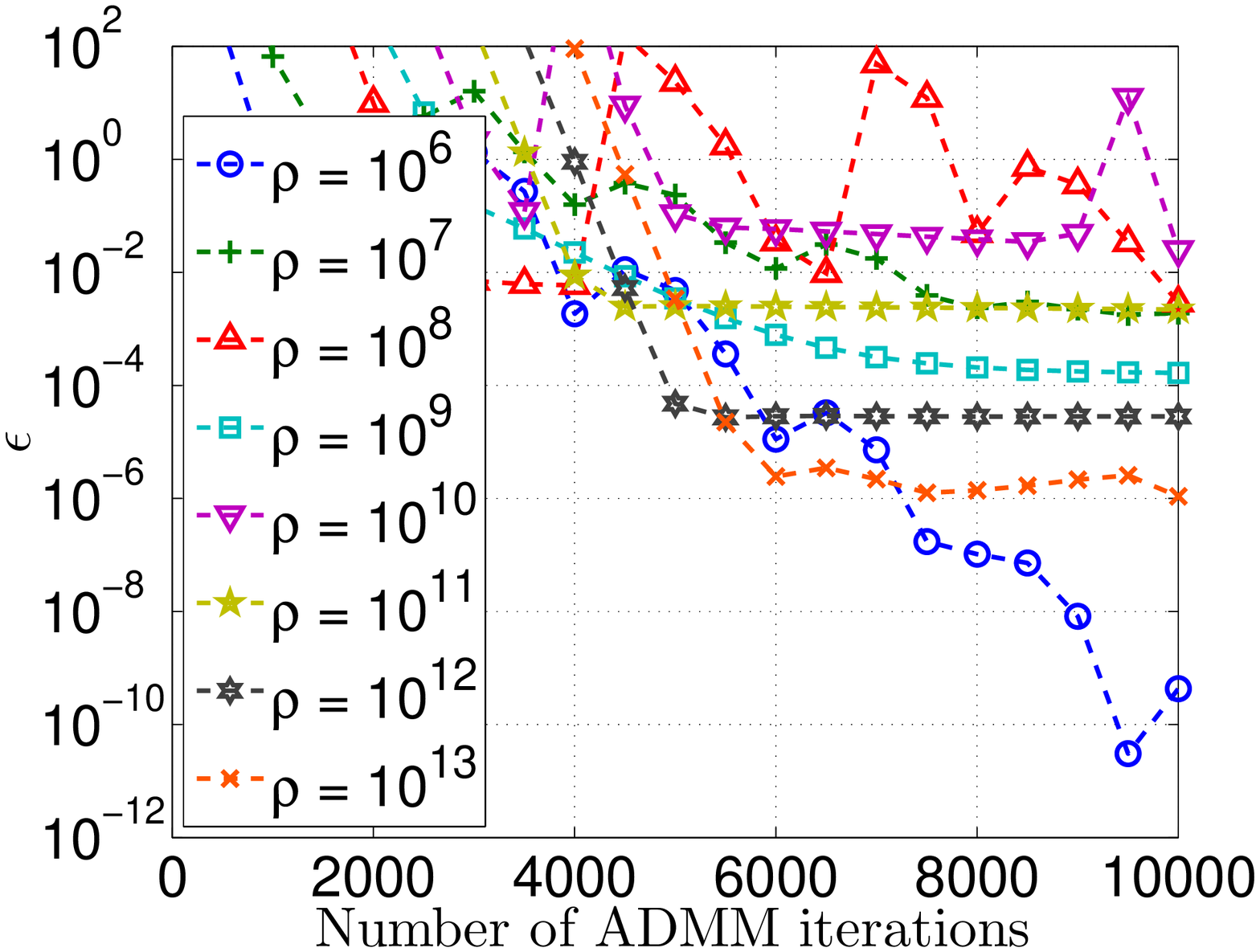}
                \caption{\subCaptionFontSize $14$ bus network}
                \label{fig:ADMM_iter_vs_cons_14bus}
        \end{subfigure}
        \begin{subfigure}[H]{0.24\textwidth}
                \includegraphics[width=\textwidth]{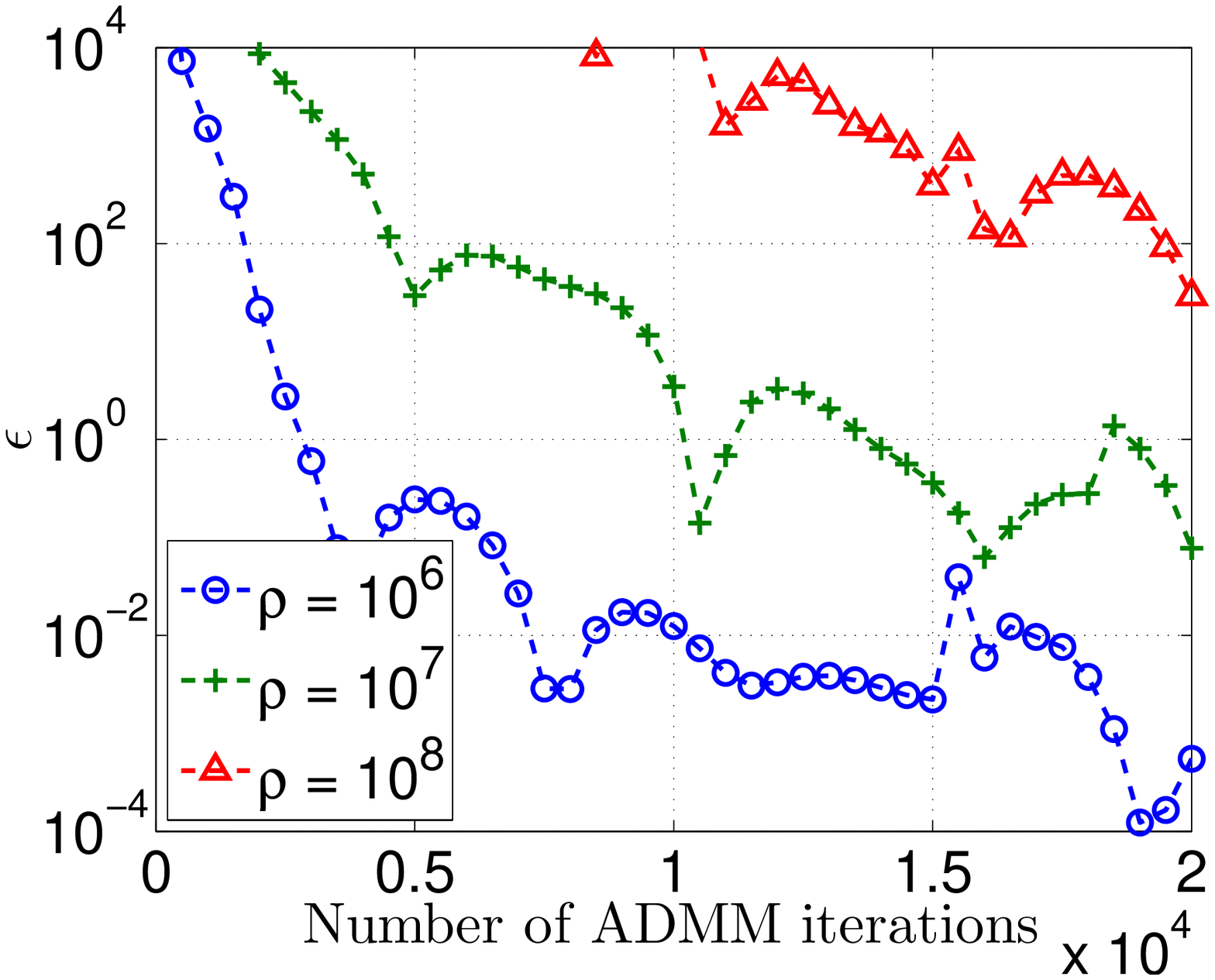}
                \caption{\subCaptionFontSize $30$ bus network}
                \label{fig:ADMM_iter_vs_cons_30bus}
        \end{subfigure}
               \caption{\captionFontSize $\epsilon$ versus number of ADMM iterations.   }
                \label{fig:ADMM_iter_vs_epsilon}
            \vspace{-5mm}
  \end{figure}

 Fig.~\ref{fig:ADMM_iter_vs_epsilon} depicts $\epsilon$ at every $500$'th ADMM iterations, for different $\rho$'s.
   In contrast to $\delta$, decreasing trend in $\epsilon$ is not necessary to obtain a feasible solution to the problem.
   However, under \emph{Assumption~\ref{assumption_1}}, as $\delta$ and $\epsilon$ go to zero, the algorithm converges to KKT optimal point. 
   Therefore, the decreasing trend in $\epsilon$, which is observed from the results,  is desired.
   In the case of the $3$, $9$ and $14$ bus $\epsilon$ reaches values between $10^{-2}$ and $10^{-11}$ in almost every case.
   However, in the $30$ bus example, only when $\rho=10^6$ does epsilon reach below $10^{-2}$.

\BlueText{

 \subsection{Convergence and scalability properties}


By convention, the  objective value of problem~\eqref{eq:dist_formulation} is $\infty$ if the problem is infeasible and is given by~\eqref{eq:dist_formulation_obj} if it is feasible. Therefore, when computing the objective function of the problem, one has to verify whether the constraints~\eqref{eq:local_variables_set}-\eqref{eq:dist_formulation_consistency_constraints}  are feasible or not. Based on Fig.~\ref{fig:subproblem_histogram}, the feasibility of the subproblem variables including $p_k\sr{G}$  is on the order of $10^{-10}$ in the worst case for every ADMM iteration [compare with~\eqref{eq:degree_of_feasiblility}]. 
 In other words, $(p_k\sr{G})_ {k\in \mathcal{G}}$ returned by \emph{Algorithm~2} in every ADMM iteration is feasible (with very high precision). Therefore, our proposed \emph{Algorithm~1}, which includes \emph{Algorithm~2} as a subroutine, ensures the feasibility of constraints~\eqref{eq:local_variables_set} -\eqref{eq:nonconvex_inequality} (with very high precision).
However, the feasibility of the remaining constraint~\eqref{eq:dist_formulation_consistency_constraints} has to be verified, in order to compute a sensible operating point.  In the sequel, we numerically analyse the feasibility of the constraint~\eqref{eq:dist_formulation_consistency_constraints} together with the objective value computed by using~\eqref{eq:dist_formulation_obj}.

  Fig.~\ref{fig:ADMM_obj_vs_cons} shows the objective value versus $\delta$,  at every $100$ or $200$' ADMM iterations.
   \begin{figure}
        \centering
        \begin{subfigure}[t]{0.24\textwidth}
                \includegraphics[width=\textwidth]{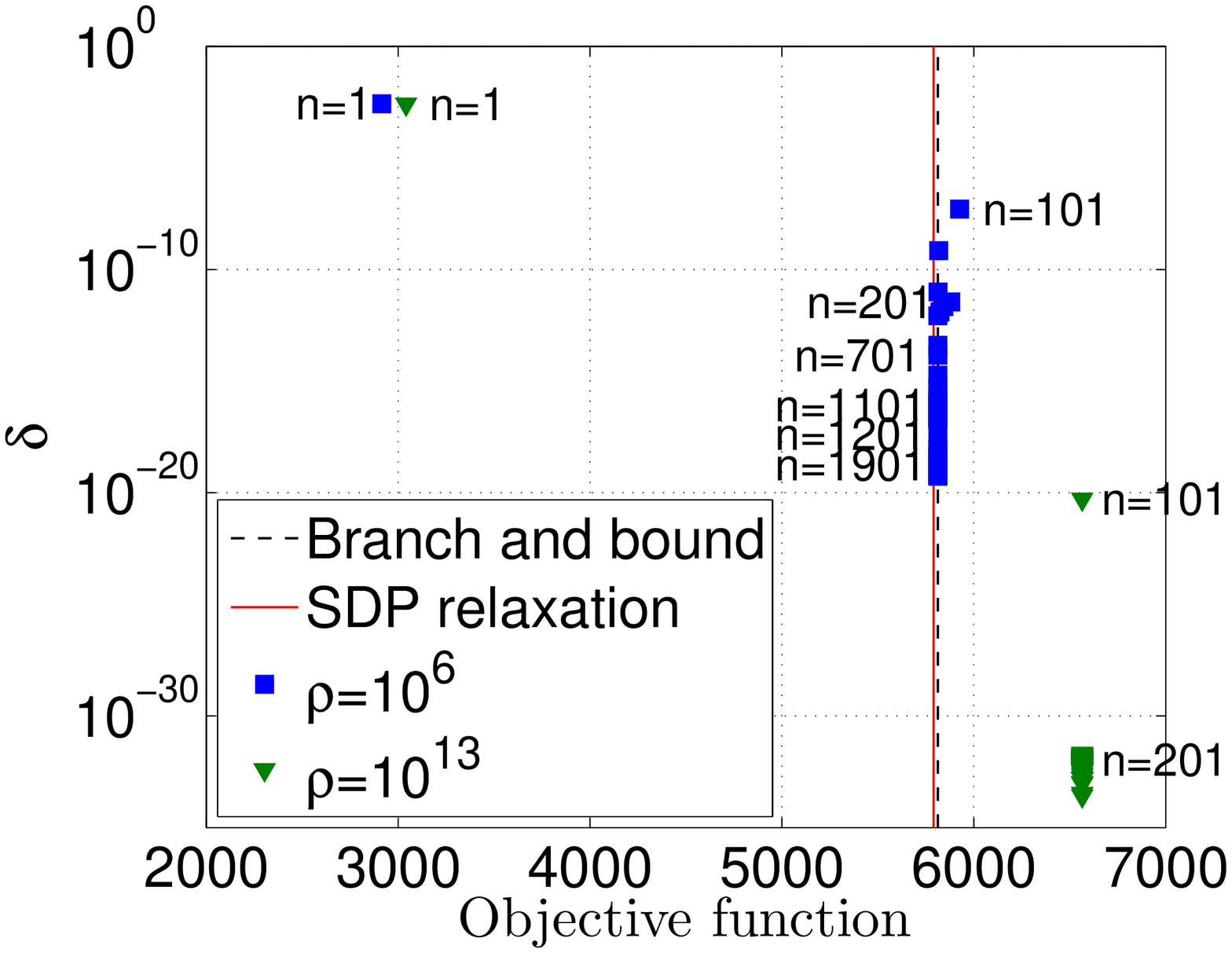}
                \caption{\subCaptionFontSize $3$ bus network }
                \label{fig:ADMM_obj_vs_cons_3bus}
        \end{subfigure}%
        ~
        \begin{subfigure}[t]{0.24\textwidth}
                \includegraphics[width=\textwidth]{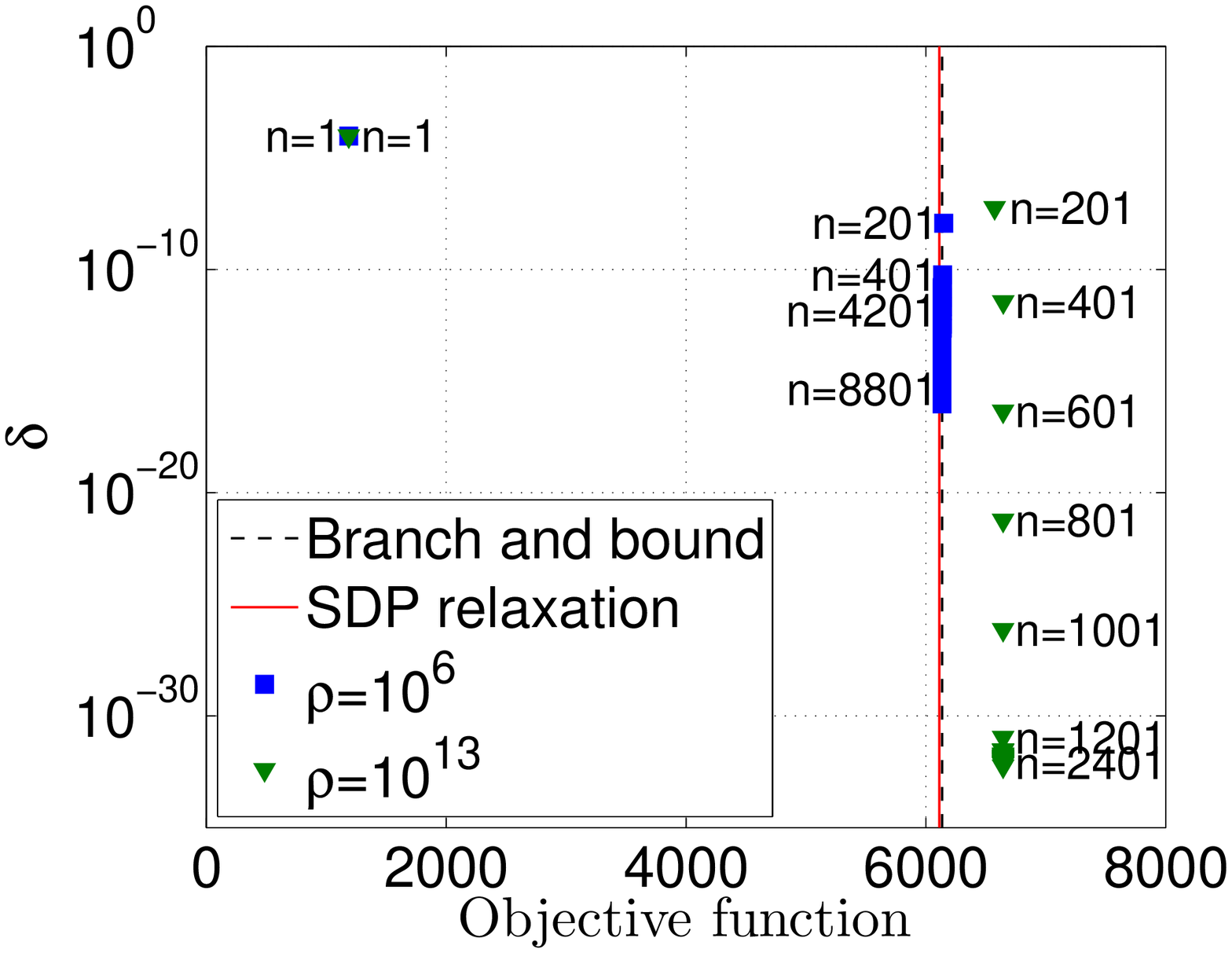}
                \caption{\subCaptionFontSize $9$ bus network}
                \label{fig:ADMM_obj_vs_cons_9bus}
        \end{subfigure}
        \begin{subfigure}[H]{0.24\textwidth}
                \includegraphics[width=\textwidth]{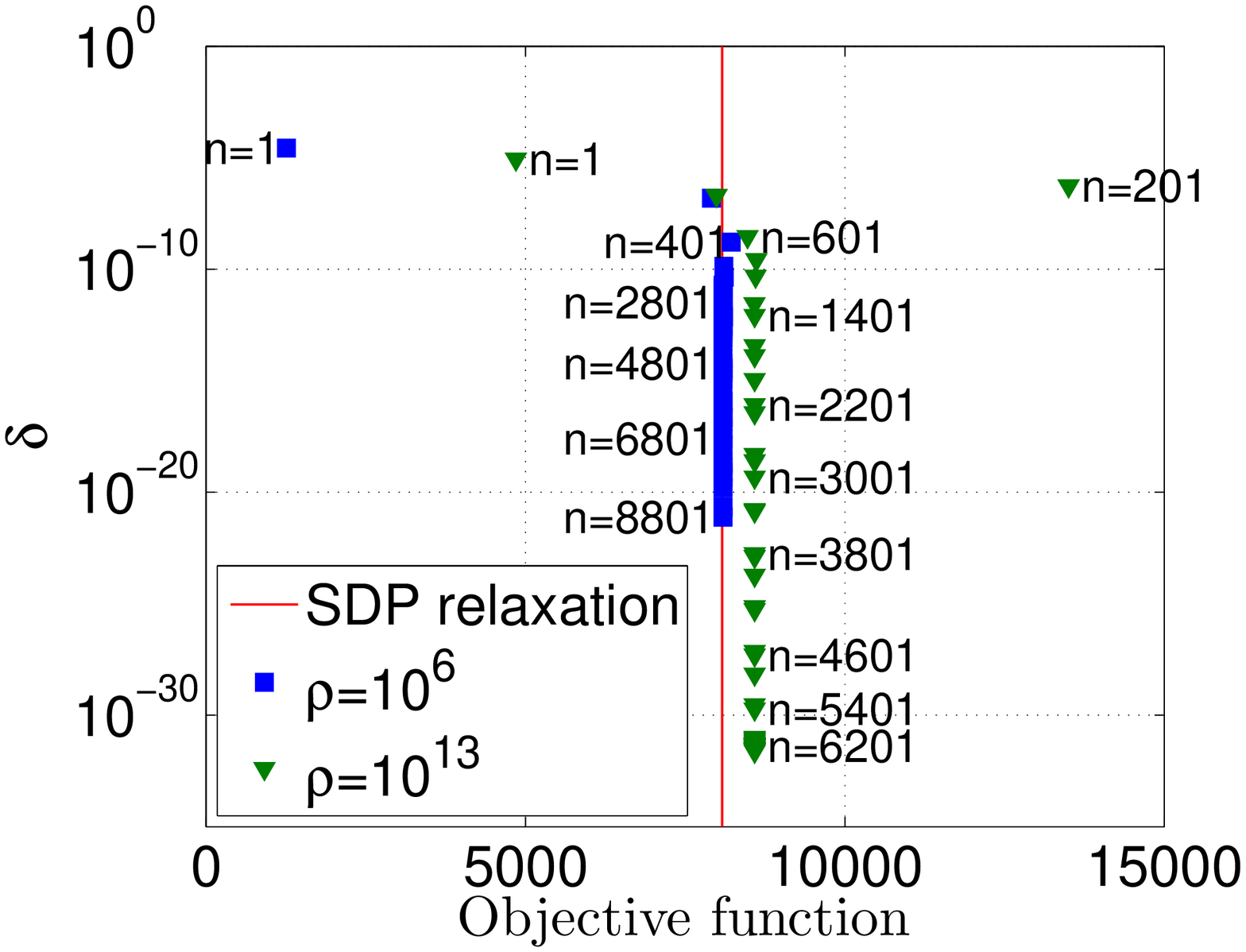}
                \caption{\subCaptionFontSize $14$ bus network}
                \label{fig:ADMM_obj_vs_cons_14bus}
        \end{subfigure}
        \begin{subfigure}[H]{0.24\textwidth}
                \includegraphics[width=\textwidth]{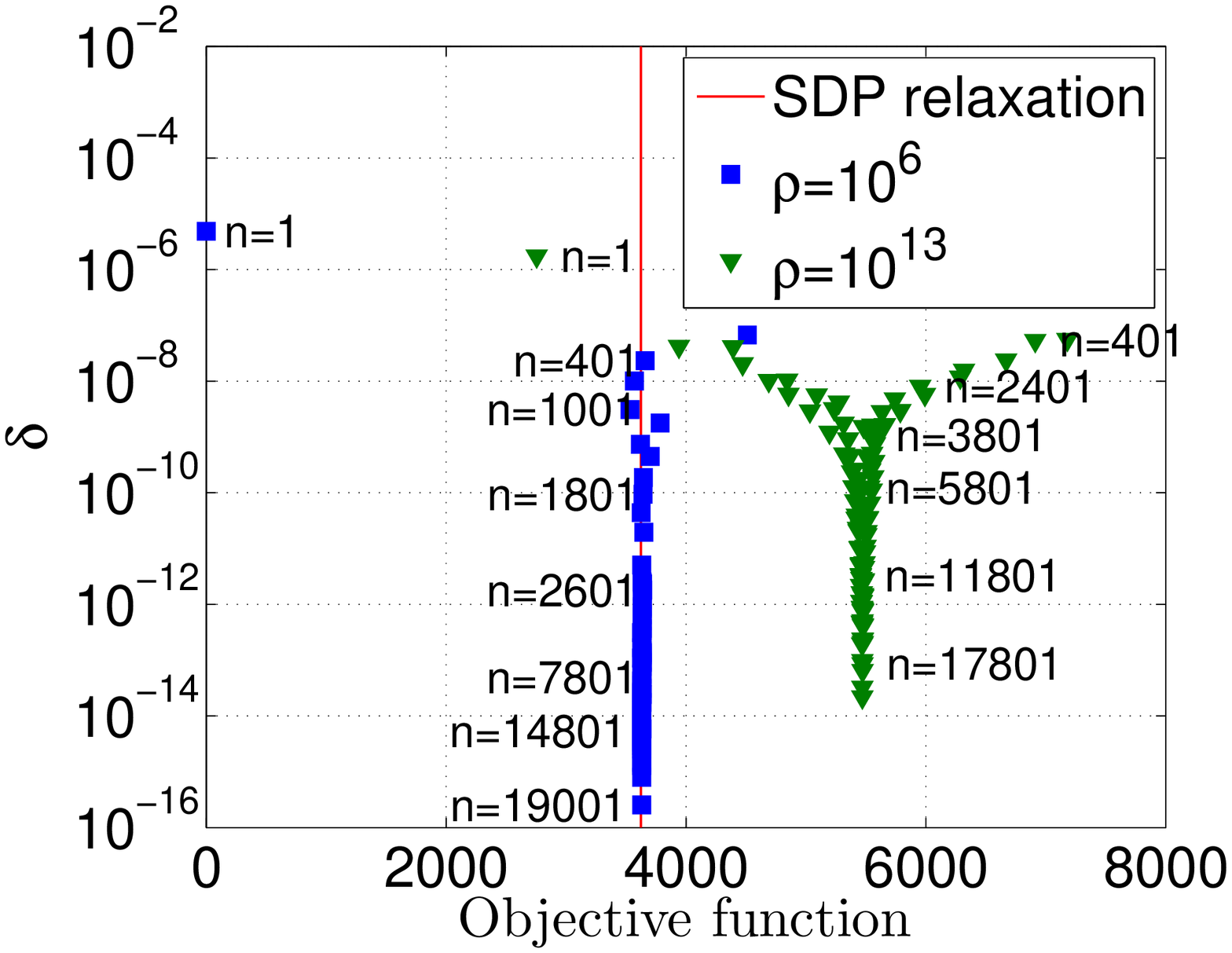}
                \caption{\subCaptionFontSize $30$ bus network}
                \label{fig:ADMM_obj_vs_cons_30bus}
        \end{subfigure}
               \caption{\captionFontSize $\delta$ versus the objective function value.  }
                \label{fig:ADMM_obj_vs_cons}
              \vspace{-5mm}
  \end{figure}
 In the case of the $3$ and $9$ bus examples we compare the objective value with the branch and bound algorithm from \cite{Gopalakrishnan-BandB-OPF-2012} where the relative tolerance, the difference between the best upper and lower bounds relative to the best upper bound, is $0.001$.
 The upper bound is obtained from Matpower and the lower bound is obtained by using the matlab toolbox YALMIP \cite{YALMIP} and the solver SEDUMI \cite{S98guide} to solve the dual of the SDP relaxation.
 In the case of the $14$ and $30$ bus examples, the branch and bound algorithm failed due to memory errors.
 In all cases we compare our results with the SDP relaxation from \cite{Lavaei-2012-OPF-TPS}. 
 The results show that the algorithm converges to some objective value in relatively few iterations, which can even be optimal with an appropriate choice of $\rho$. For example, for the considered cases $\rho=10^6$ yields almost optimal objective values.
 Moreover, as desired the consistency metric $\delta$ is driven towards zero as the number of ADMM iterations increases.

 Fig.~\ref{fig:rel_obj} depicts the relative objective function $(|f-f^{\star}|/f^{\star})$ for the 3 and 9 bus examples.  The results are consistent with thous presented in Fig.~\ref{fig:ADMM_obj_vs_cons}. For example, in the case of $\rho=10^6$ the relative objective function value is on the order of $10^{-6}$. Results suggest that a proper choice of $\rho$ is beneficial to achieve a good network operating point.
   \begin{figure}
        \begin{subfigure}[t]{0.23\textwidth}
                \includegraphics[width=\textwidth]{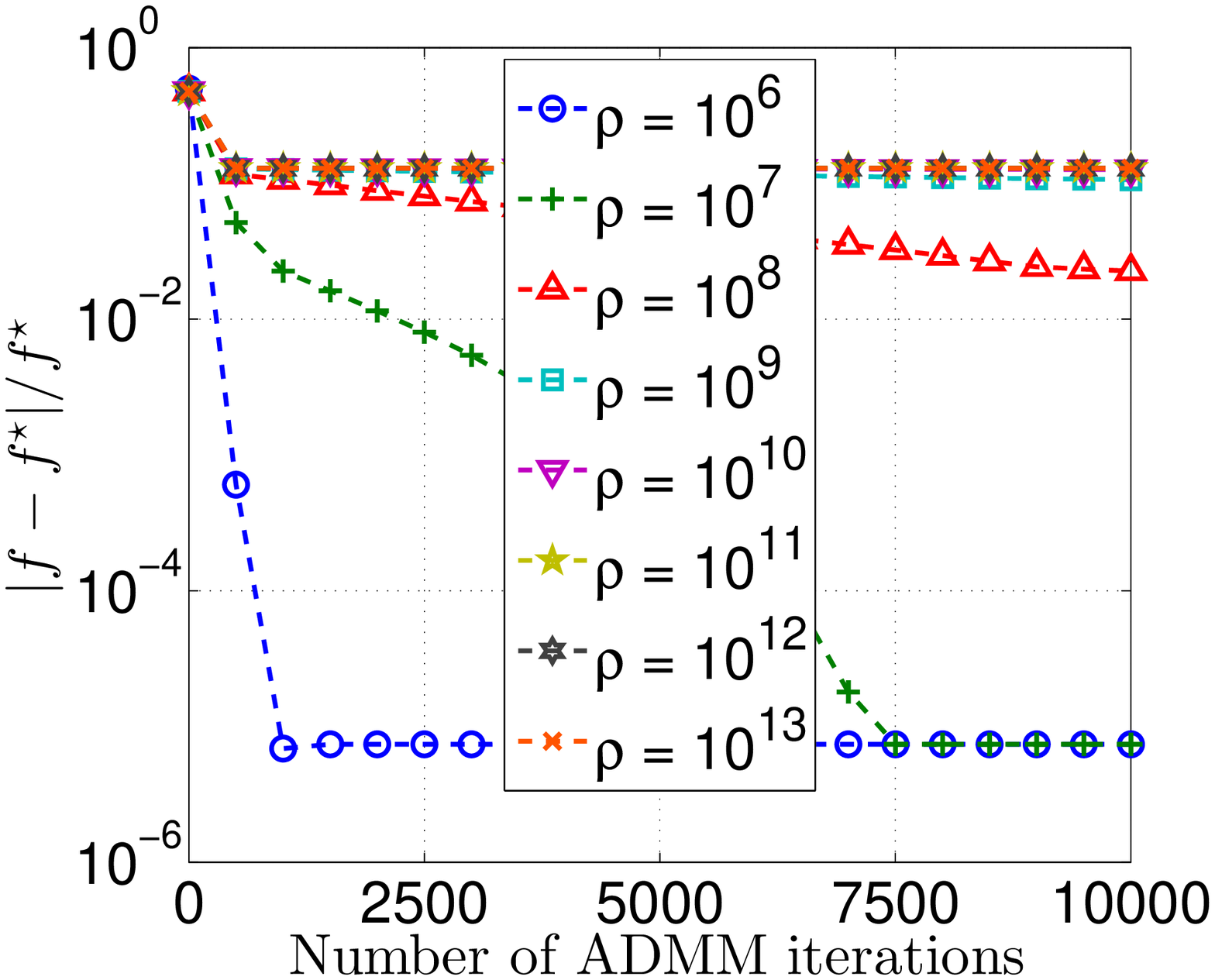}
                \caption{3 bus example}
                \label{fig:r_rel_3_buses}
        \end{subfigure}
        ~
        \begin{subfigure}[t]{0.23\textwidth}
                \includegraphics[width=\textwidth]{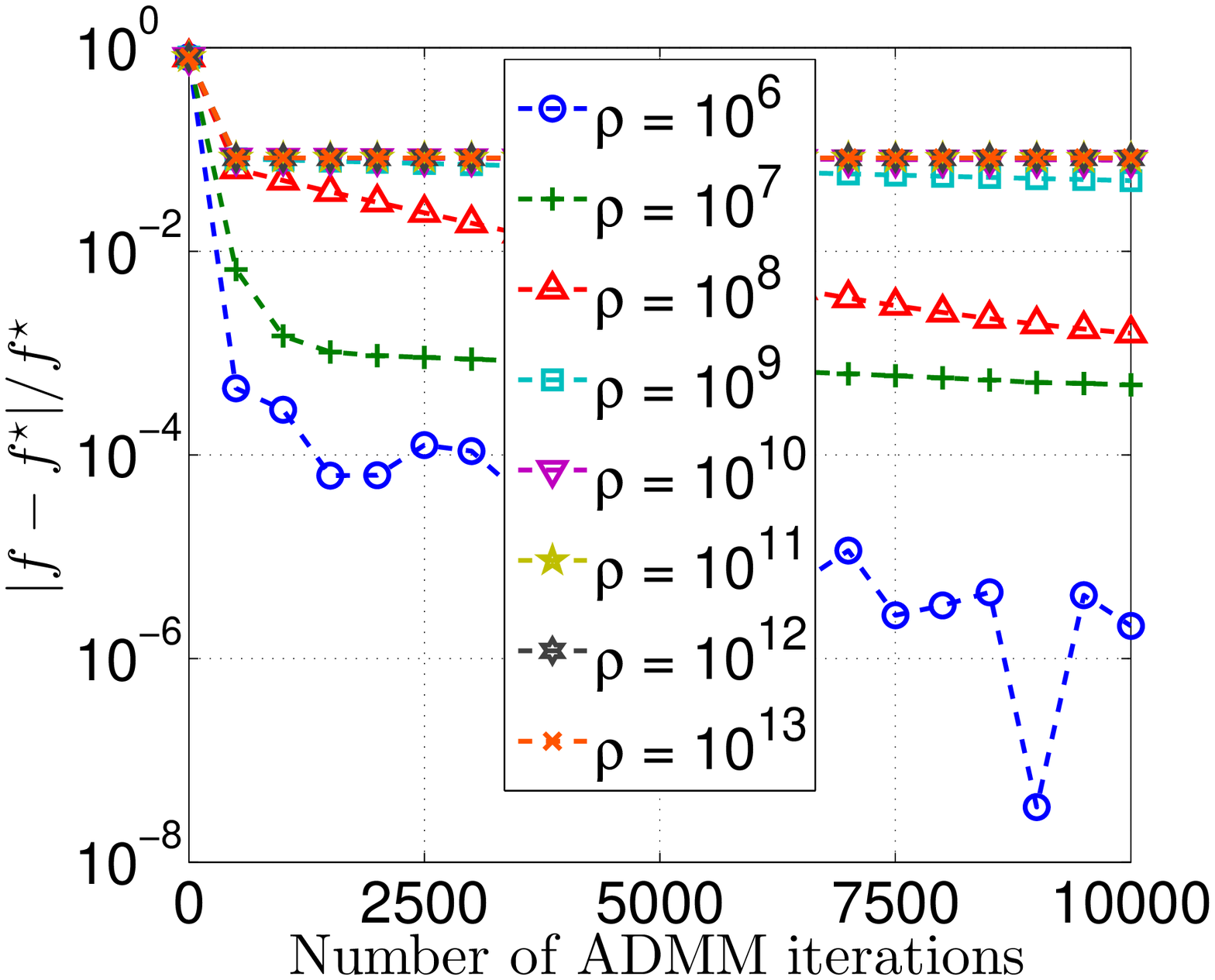}
                \caption{9 bus example }
                \label{fig:r_rel_3_buses}
        \end{subfigure}
        \caption{\captionFontSize  Relative objective function.}
        \label{fig:rel_obj}
  \end{figure}


  To study the scalability properties of the proposed algorithm, we compute the CPU time, relative objective value, $\delta$, and $\epsilon$ of the algorithm for all the considered examples. In the case of $3$, $9$, $14$, and $30$ bus examples, we choose $\rho=10^6$ and in the case of $118$ and $300$ bus examples we chose $\rho=10^7$.

 Fig.~\ref{fig:time_comparison_admm_times} shows the parallel running times $T_p$ versus ADMM iterations. In particular, we define $T_p=T_s/|\mathcal{N}|$,  where $T_s$ is the sequential CPU time. 
The behavior of the plots in the case of $9$, $14$, $30$, $118$, and $300$ bus examples are very similar. In other words, the parallel running time $T_p$ is independent of the number of busses, indicating promising scalability of the proposed algorithm.  Note that the $3$ bus example has to handle more variables per subproblem, compared with the other examples, see column 6 of Table~\ref{table:Introduction_to_test_problems}. This is clearly reflected in the plot of $3$ bus example, as an increase of the associated parallel running time $T_p$.

 Fig.~\ref{fig:time_comparison_relative_objective_function} depicts the relative objective function, $|f-f^{\star}|/f^{\star}$.
  In the case of $3$, $9$, $118$, and $300$ bus examples the global optimum $f^{\star}$ is found by branch and bound algorithm. However, in the case of  $14$ and $30$ bus examples,  branch and bound algorithm failed, and therefore the best known objective value found by Matpower was considered as $f^{\star}$.~\footnote{ Since the $118$ and $300$ bus examples have \emph{zero duality gap}, the branch and bound algorithm worked efficiently. However, in the case of  $14$ and $30$ bus examples, where there is \emph{nonzero duality gap}, the branch and bound algorithm failed.} Results show that for large and relatively large test examples (e.g., $N=30$, $118$, $300$) the relative objective function value is on the order of $10^{-3}$ and is not affected by the network size. A similar independence of the performance is observed for very small test examples as well (e.g., $N=3$, $9$, $14$) with relative objective function values on the order of $10^{-6}$. The reduction of the relative objective function values of smaller networks compared with larger network examples are intuitively expected due to substantial size differences of those networks.


 Fig.~\ref{fig:time_comparison_consistency_constraints} and Fig.~\ref{fig:time_comparison_eps} depict $\epsilon$ and $\delta$ as a function of ADMM iterations. Results show that, irrespective of the number of buses, the metrics  $\delta$ and $\epsilon$ are decreasing as desired. Results further suggest that those values are driven towards small values as ADMM iterations increase.

   \begin{figure}
        \begin{subfigure}[t]{0.23\textwidth}
                \includegraphics[width=\textwidth]{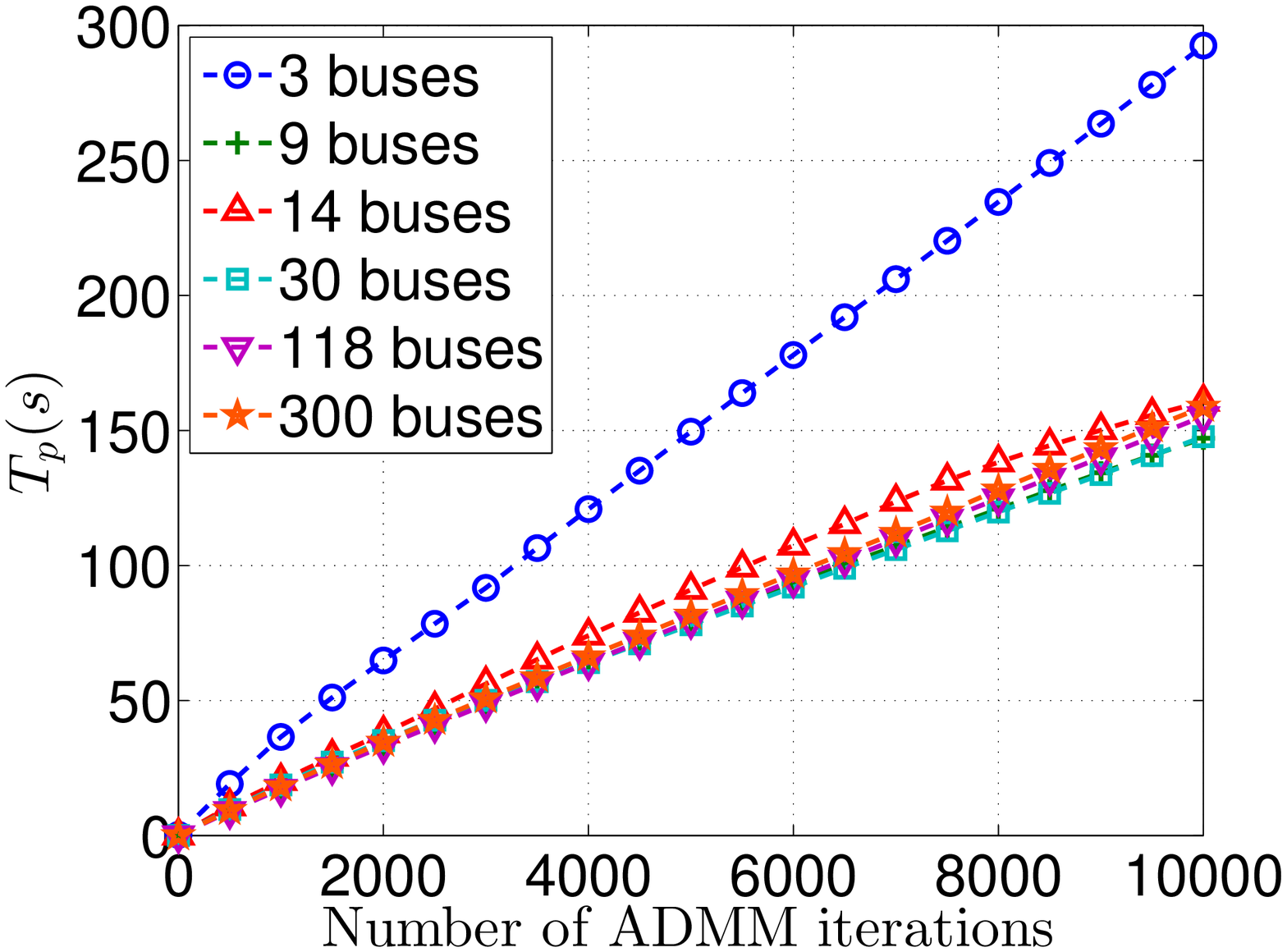}
                \caption{\subCaptionFontSize $T_p$ }
                \label{fig:time_comparison_admm_times} 
        \end{subfigure}
        ~
        \begin{subfigure}[t]{0.23\textwidth}
                \includegraphics[width=\textwidth]{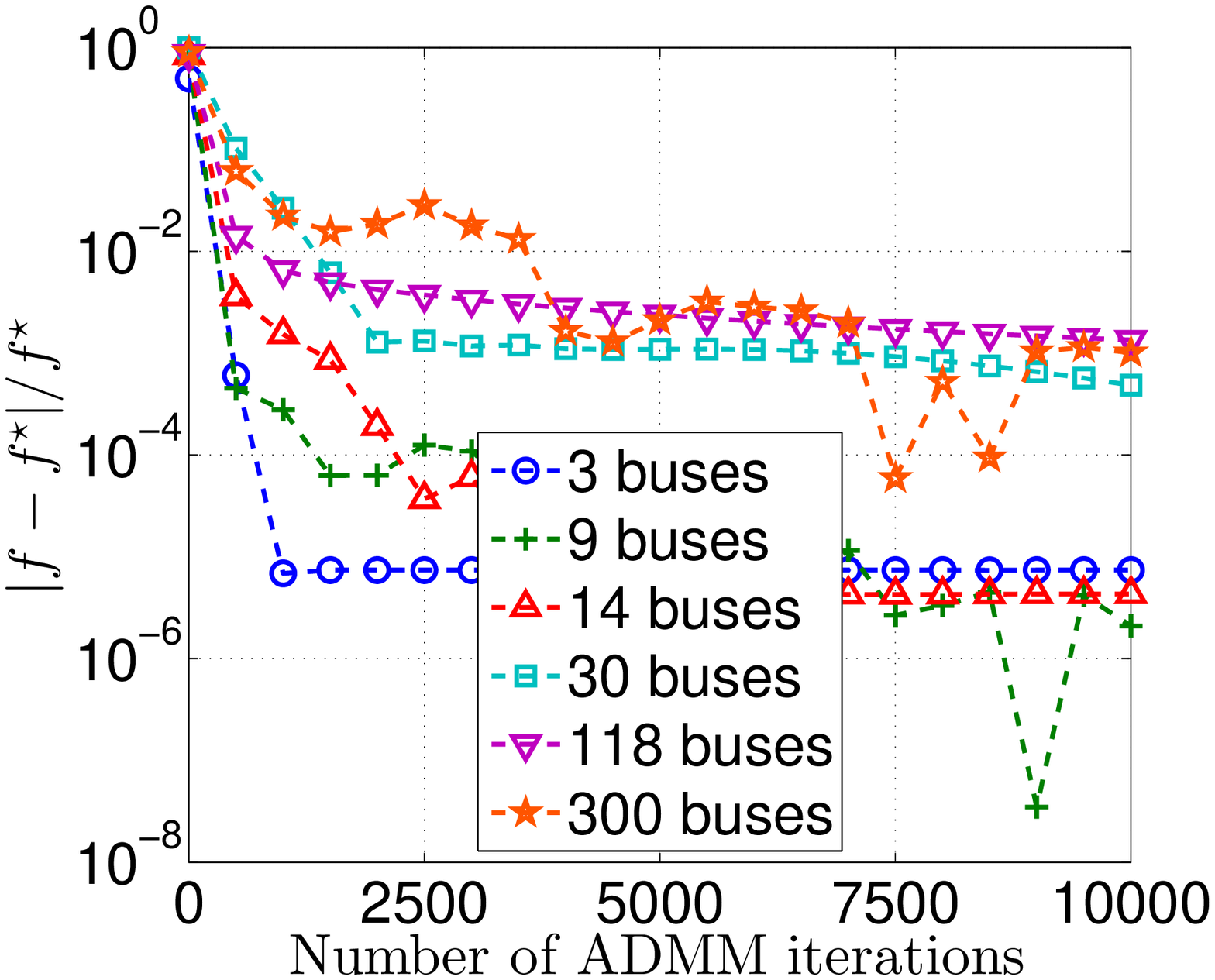}
                \caption{\subCaptionFontSize $|f-f^{\star}|/ f^{\star}$   }
                \label{fig:time_comparison_relative_objective_function}
        \end{subfigure}\\ \centering
        \begin{subfigure}[t]{0.23\textwidth}
                \includegraphics[width=\textwidth]{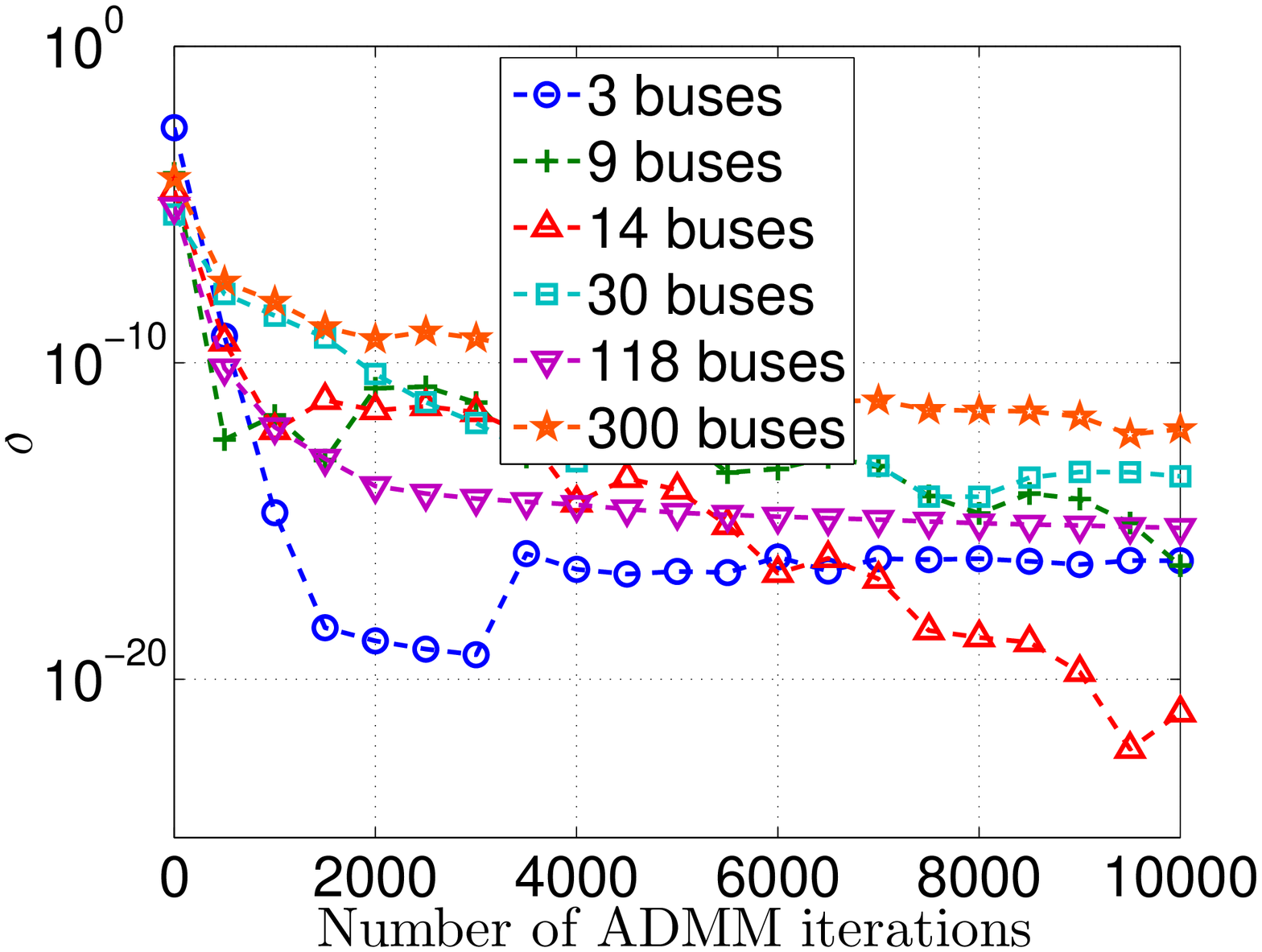}
                \caption{\subCaptionFontSize $\delta$  }
                \label{fig:time_comparison_consistency_constraints}
        \end{subfigure}
        \begin{subfigure}[t]{0.23\textwidth}
                \includegraphics[width=\textwidth]{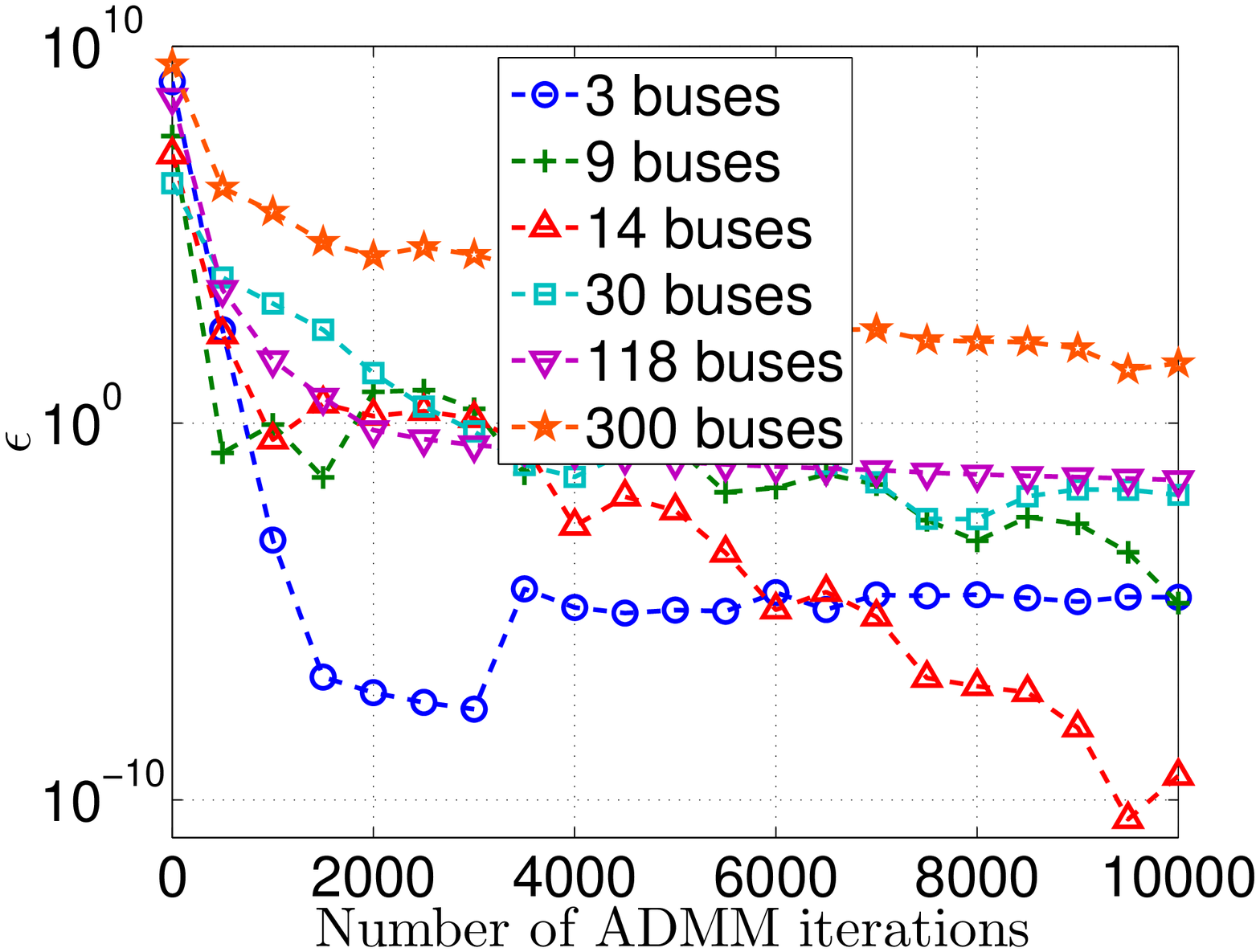}
                \caption{\subCaptionFontSize $\epsilon$  }
                \label{fig:time_comparison_eps}
        \end{subfigure}
        \caption{\captionFontSize $T_p$, $|f-f^{\star}|/ f^{\star}$, $\delta$, and $\epsilon$ as a function of ADMM iterations. }
        \label{fig:time_comparison}
         \vspace{-5mm}
  \end{figure}


\begin{table*}[h]\scriptsize 
  \centering 
 \caption{\captionFontSize  Comparison of output of ADMM-DOPF when at ADMM iteration $n=3000$ and $n=10000$ with the SDP relaxation~\cite{Lavaei-2012-OPF-TPS}, branch and bound~\cite{Gopalakrishnan-BandB-OPF-2012}, and Matpower~\cite{Zimmerman-Mat-ste-sta-ope-pla-and-ana-too-for-pow-sys-rea-and-edu-2011}. $T_s$, $T_p$, obj, and \# iter indicate the sequential and parallel running times in seconds, the objective value in \$/hours, and the number of iterations, respectively.   }
  \label{table:comparison_benchmarks_vs_proposed_algorithm}
  \begin{tabular}{| l || l | l || l | l || l | l || l | l || l | l |} 
            \hline
           \multirow{2}{*}{Example}        &   \multicolumn{2}{ c ||}{SDP}   & \multicolumn{2}{ c ||}{ Branch \& Bound}& \multicolumn{2}{ c ||}{Matpower} & \multicolumn{2}{ c ||}{ADMM-DOPF, n=3000} & \multicolumn{2}{ c |}{ADMM-DOPF, n=10000}  \\   
                                &      \multicolumn{1}{ c }{obj}          & $T_s$   &    \multicolumn{1}{ c }{obj} & $T_s$                             &    \multicolumn{1}{ c }{obj}          & $T_s$  & \multicolumn{1}{ c }{obj} & \multicolumn{1}{ c ||}{$T_p$} &               \multicolumn{1}{ c }{obj} & \multicolumn{1}{ c |}{$T_p$}   \\   \hline
            $3$ buses    &   $5789.9$  &  $5.09$    & $5812.6$ &  8                            &  $5812.6$  & 0.08   &   $5812.6$ & 91.73 & 5812.6   & 292.54  \\  
            $9$ buses    &   $6113.2$  &  $0.48$    & $6135.2$ &  $1.6\times 10^5$ &  $6135.2$  &  0.08  & $6135.9$   &  49.19   & 6135.2 & 147.16    \\   
            $14$ buses  &   $8079.6$  &  $0.44$    & $-$           &       $-$                     &  $8092.4$ & 0.1    & $8092.9$    & 56.12 & 8092.4  & 160.71   \\   
            $30$ buses  &   $3624.0$  &  $0.85$    & $-$           &       $-$                    &  $3630.7$  & 0.12  & $3634.9 $    & 51.53 & 3632.5  & 152.95 \\    
            $118$ buses  &  $129654.1$  &  $10.30$  & $129660.6$  &   $1.8\times 10^1$  &  $129660.6$   & 0.2    & $130094.3$    &  48.62 &  129835.2 & 155.29  \\    
            $300$ buses  &  $719711.7$ & $154.35$ &  $719 725.1$  &    $ 4.6\times 10^3$&  $719 725.1$   & 0.42  & $732629.1$    &  51.19 &  720449.4 & 158.97 \\   \hline 
  \end{tabular}
  \vspace{-4mm}
\end{table*}

 Table \ref{table:comparison_benchmarks_vs_proposed_algorithm} shows the running time and the objective value obtained by different approaches. As benchmarks, we consider the centralized algorithms, SDP relaxation~\cite{Lavaei-2012-OPF-TPS}, branch and bound~\cite{Gopalakrishnan-BandB-OPF-2012}, and Matpower~\cite{Zimmerman-Mat-ste-sta-ope-pla-and-ana-too-for-pow-sys-rea-and-edu-2011}. Table~\ref{table:comparison_benchmarks_vs_proposed_algorithm}  shows that our proposed method yields \emph{network operating points}, which are almost optimal, where the discrepancy with respect to the optimal is on the order of $0.1\%$ (respectively, $1\%$) or less with $10000$ (respectively, $3000$) ADMM iterations. Note that the running time of ADMM-DOPF is insensitive to the network size, see Fig.~\ref{fig:time_comparison_admm_times} for more details. However, even in small networks (e.g., case with $N=14$, $N=30$), the running time of  branch and bound algorithm can explode. This is expected because the worst case complexity of the branch and bound algorithm grows exponentially with the problem size. Results further suggest that the running time of the centralized algorithm SDP relaxation increases as the network size grows, unlike the proposed  ADMM-DOPF. Note that the running time of ADMM-DOPF is large compared to the centralized Matpower. However, those values can further be reduced if ADMM-DOPF is deployed in a parallel computation environment, where every subproblem is handled at a dedicated set of resources, including processors, memory, among others.  In addition to the centralized benchmarks, we also consider the decentralized one proposed in \cite{Sun_2013}, which employs ADMM for general nonconvex OPF. However, the results of \cite{Sun_2013} are not documented in Table~\ref{table:comparison_benchmarks_vs_proposed_algorithm}, because for all the considered examples, the algorithm therein did not converge. This agrees with the numerical results of \cite{Sun_2013}, where the authors mentioned that the convergence of their algorithm is more sensitive to the initial point in the case of mesh networks \cite[p.~5]{Sun_2013}, so does ADMM-DOPF. We note that their method converges if it is initialized close to the optimal solution. However, in practice, such an initialization point in unknown, thus limiting dramatically the applicability of the method in \cite{Sun_2013}.
%

Finally, from all our numerical experiments discussed above, we note that the power losses in the flow lines are typically on the order of $4.4\%$ (or less) of the total power flow in the line. Because the losses are not negligible, approximations such as the linearization of power flow equations can be less applicable to compute better network operating points.
 }


\section{conclusions} \label{sec:conclusions}

We proposed a distributed algorithm for the optimal power flow (OPF), by decomposing the OPF problem among the buses that compose the electrical network. A light communication protocol among neighboring buses is needed during the algorithm, resulting in high scalability properties. The subproblems related to each bus capitalize on sequential convex approximations to gracefully manipulate the nonconvexity of the problem. We showed the convergence of subproblem solutions to a local optimum, under mild conditions. Furthermore, by using the local optimality results associated with the subproblems, we quantified the optimality of the overall algorithm. We evaluated the proposed algorithm on a number of test examples to demonstrate its convergence properties and to compare it with the global optimal method. In all considered cases the proposed algorithm achieved close to optimal solutions. Moreover, the proposed algorithm showed appealing scalability properties when tested on larger examples.

  
\appendices

 \section{On the use of quadratic programming QP solvers}
    Note that, not all the constraints of problem~\eqref{eq:the_approximated_subproblem} are affine (or linear). In particular, constraints~\eqref{eq:convex_nonlinear_inequality_sub_approx} and \eqref{eq:nonconvex_inequality_sub_approx} are not affine. Therefore, QP solvers are not directly applied to solve the problem. However, if constraints~\eqref{eq:convex_nonlinear_inequality_sub_approx} and \eqref{eq:nonconvex_inequality_sub_approx} are approximated by using affine constraints, then QP is readily applied to the modified problem.

     Let us start by considering the feasible regions defined by~\eqref{eq:nonconvex_inequality_sub_approx}, which accounts for $\check{\mathcal{X}}^k_r$, $r\in\{1,\ldots,|\mathcal{N}_k|\}$, see Fig.~\ref{fig:subproblems_subset_XX}. Next, we approximate the nonlinear boundary of $\check{\mathcal{X}}^k_r$ by affine functions as depicted in Fig.~\ref{fig:subproblems_subset_XXX}. We denote by $\check{\mathcal{Y}}^k_r$ the approximated polyhedral set. We can apply similar ideas to approximate the feasible regions specified by~\eqref{eq:convex_nonlinear_inequality_sub_approx} [\emph{cf} \eqref{eq:dist_formulation_current_circle}-\eqref{eq:dist_formulation_power_circle}], where we use $\check{\gamma}_k$ to denote the resulting affine function. Finally, the idea is to find the desired optimal solution of problem~\eqref{eq:the_approximated_subproblem} by constructing a series of sets of the form $\check{\mathcal{Y}}^k_r$ and affine functions of the form  $\check{\gamma}_k$  that approximate the feasible set specified by~\eqref{eq:convex_nonlinear_inequality_sub_approx} and \eqref{eq:nonconvex_inequality_sub_approx} in an \emph{increasing precision}. The QP based algorithm to solve problem~\eqref{eq:the_approximated_subproblem} can be summarized as follows.
 \noindent\rule{\linewidth}{0.3mm}
\\
\emph{Algorithm 3}: \ \textsc{\small{QP to solve Problem~\eqref{eq:the_approximated_subproblem}}}
\begin{enumerate}
   \item Initialize: Given the initial approximated set $\check{\mathcal{Y}}^k_r$ and affine function $\check{\gamma}_k$. Let $\bar m=1$.
   \item \label{al:subproblem_algorithm_private_variable_update_}
            Solve the QP:
                \begin{subequations}  \label{eq:the_approximated_subproblem_QP}
					      \begin{align}
					        & {\text{min}}
					          & &   f_k\sr{G} (p_k\sr{G}) {+}\vec{y}_k\tran(\vec{v}_k-\bec{E}_k\vec{v})  {+} \frac{\rho}{2}||\vec{v}_k{-}\bec{E}_k \vec{v}||_2^2 \label{eq:the_approximated_subproblem_obj_QP}\\
					            & \text{s. t.}
          & &  \vec{z}_k=({p}_k\sr{G}, {q}_k\sr{G}, {p}_k, {q}_k, {i}_k\sr{Re}, {i}_k\sr{Im}, \vec{v}_k\sr{Re}, \vec{v}_k\sr{Im},  \\
          & & &          \hspace{2.6cm} \bec{i}_{k}\sr{Re}, \bec{i}_{k}\sr{Im}, \bec{p}_{k}, \bec{q}_{k}) && \label{eq:local_variables_set_sub_arrox_QP} \\
          & & &  (\boldsymbol{\alpha}_k(\vec{z}_k), \hat{\boldsymbol{\lambda}}_k^{\hat{\vec{z}}_k}(\vec{z}_k),\hat{\boldsymbol{\mu}}_k^{\hat{\vec{z}}_k}(\vec{z}_k) )=\vec{0} \\
       & & &   ( \boldsymbol{\beta}_k(\vec{z}_k) , \check{\boldsymbol{\gamma}}_k(\vec{z}_k))\leq \vec{0}                                    \\
       & & & ((\vec{v}\sr{Re}_k)_r,(\vec{v}\sr{Im}_k)_r)\in\check{\mathcal{Y}}^k_r,  ~ r=1,\ldots,|\mathcal{N}_k|,
					      \end{align}
					\end{subequations}
	    where the variables are ${p}_k\sr{G}$, ${q}_k\sr{G}$, ${p}_k$, ${q}_k$, ${i}_k\sr{Re}$, ${i}_k\sr{Im}$, $\vec{v}_k\sr{Re}$, $\vec{v}_k\sr{Im}$, $\bec{i}_{k}\sr{Re}$, $\bec{i}_{k}\sr{Im}$, $\bec{p}_{k}$, $\bec{q}_{k}$, and $\vec{z}_k$. The solution corresponding to the variable $\vec{z}_k$, $((\vec{v}\sr{Re}_k)_r,(\vec{v}\sr{Im}_k)_r)$ are denoted by $\vec{z}^{(\bar m)}_k$, $\vec{v}_r^{(\bar m)}$, respectively and and all the dual optimal variables are denoted by $\vec{u}^{(\bar m)}_k$.
     \item Stopping criterion:  If $\boldsymbol{\gamma}_k\big(\vec{z}^{(\bar m)}_k\big)\leq 0$ and $\vec{v}_r^{(\bar m)}\in\check{\mathcal{X}}^k_r$ for all $r\in\{1,\ldots,|\mathcal{N}_k|\}$, \textsc{stop} and return $(\vec{z}^{(\bar m)}_k,\vec{u}^{(\bar m)}_k)$. Otherwise, increase the precession of set $\check{\mathcal{Y}}^k_r$ and function $\check{\gamma}_k$ by adding a hyper plane and an affine function, respectively, set $\bar m\asign \bar m +1$ and go to step~2.
\end{enumerate}
\vspace{-3mm}
\rule{\linewidth}{0.3mm}

\noindent The set $\check{\mathcal{Y}}^k_r$ is initialized in the first step by approximating the exterior boundary of the donut $\mathcal{X}^k_r$~[Fig.~\ref{fig:subproblems_doughnut_set}] by an equilateral octagon as shown in Fig.~\ref{fig:subproblems_subset_XXX}, and $\check{\gamma}_k$  is initialized correspondingly.
 The second step simply involves solving a QP programming problem.
 The algorithm terminates in the third step if $\boldsymbol{\gamma}_k(\vec{z}^{(\bar m)})\leq 0$ and $(\vec{v}_k^{(\bar m)})_r\in\check{\mathcal{X}}^k_r$ for all $r\in\{1,\ldots,|\mathcal{N}_k|\}$.
 However, if $(\vec{v}_k^{(\bar m)})_r \in \check{\mathcal{Y}}^k_r\setminus \check{\mathcal{X}}^k_r$ we increase the precession of $\check{\mathcal{Y}}^k_r$ by adding a hyper plane on the exterior boundary of the donut $\mathcal{X}^k_r$, so that $(\vec{v}_k^{(\bar m)})_r \notin \check{\mathcal{Y}}^k_r$.
 In particular, we set $\check{\mathcal{Y}}^k_r = \check{\mathcal{Y}}^k_r \cap \mathcal{W}$ where $\mathcal{W}$ is the halfspace
 \begin{multline}
    \mathcal{W} {=}  \big\{( (\vec{v}\sr{Re}_k)_r,(\vec{v}\sr{Im}_k)_r  ) {\in} \R^2  \big|  \ \alpha_r (\vec{v}_k\sr{Re})_r {+} \beta_r (\vec{v}_k\sr{Im})_r {\leq}\gamma_r     \big\}  ,
 \end{multline}
 where \vspace{-2mm}
 \begin{align}
     &\alpha_r  {=} \text{sign}(\Re((\vec{v}_k^{(\bar m)})_r) )\ \sqrt{   \frac{(\vec{v}_k\sr{max})_r^2 }{1{+} ( (\Im((\vec{v}_k^{(\bar m)})_r)/ (\Re((\vec{v}_k^{(\bar m)})_r) )^2} } , \notag \\
     &\beta_r   = a_r\left(  \frac{(\Im((\vec{v}_k^{(\bar m)})_r)}{(\Re((\vec{v}_k^{(\bar m)})_r)}  \right),  \hspace{1cm}
     \gamma_r = (\vec{v}_k\sr{max})_r^2 , \notag
 \end{align}
 if $\Re(x_r^\star) \neq 0$ and  
  \be
     \alpha_r  = 0,  \hspace{1cm}
     \beta_r   = \text{sign}(\Re((\vec{v}_k^{(\bar m)})_r)), \hspace{1cm}
     \gamma_r = (\vec{v}_k\sr{max})_r. \notag
 \ee
  $\check{\gamma}_k$ can be treated identically.

 \section{Proofs}
 \subsection{Proof of \emph{Proposition~1} } \label{app:proof_of_prop_1}
\begin{IEEEproof}
Obviously, problem~\eqref{eq:the_approximated_subproblem} is convex and in any iteration~$m$ of \emph{Algorithm~2}, $(\vec{z}^{(m)}_k,\vec{u}^{(m)}_k)$ [so is $(\vec{z}^\star_k,\vec{u}^\star_k)$] are primal and dual optimal, with zero duality gap. Thus, $(\vec{z}^\star_k,\vec{u}^\star_k)$ satisfies KKT conditions for problem~\eqref{eq:the_approximated_subproblem}~\cite[\S~5.5.3]{convex_boyd}. However, in order to show that $(\vec{z}^\star_k,\vec{u}^\star_k)$ satisfies KKT conditions for problem~\eqref{eq:ADMM_private_variable_update_generator_bus}, we need to show 1)  $\vec{z}^\star_k$ is primal feasible, 2) ${\vec{u}}^\star_k$ is dual feasible, 3) $\vec{z}^\star_k$ and ${\vec{u}}^\star_k$ satisfy complementary slackness conditions, and 4) derivative of the Lagrangian of problem~\eqref{eq:ADMM_private_variable_update_generator_bus} vanishes with $\vec{z}^\star_k$ and ${\vec{u}}^\star_k$~\cite[\S~5.5.3]{convex_boyd}.

We start by noting that the original functions definitions $\boldsymbol{\lambda}_k(\vec{z}_k)$ and $\boldsymbol{\mu}_k(\vec{z}_k)$ [see~\eqref{eq:nonlinear_equality1_sub} and  \eqref{eq:nonlinear_equality2_sub}] are characterized by using the basic form
\be\label{eq:block_function}
h(\underbrace{p,x_1,x_2,y_1,y_2}_{\vec{z}})= p-x_1y_1-x_2y_2 \ ,
\ee
where $p\in\R$ represents power, $x_1,x_2\in\R$ represent currents, $y_1,y_2\in\R$ represent voltages, and we have denoted $(p,x_1,x_2,y_1,y_2)$ compactly by $\vec{z}$. Let ${\hat h}^{\hat{\vec{z}}}$ denote the first order Taylor's approximation of $h$ at $\hat{\vec{z}}$. That is, ${\hat h}^{\hat{\vec{z}}}$ characterizes the basic form of the first order Taylor's approximation of function definitions $\hat{\boldsymbol{\lambda}}_k^{\hat{\vec{z}}_k}(\vec{z}_k)$ and $\hat{\boldsymbol{\mu}}_k^{\hat{\vec{z}}_k}(\vec{z}_k)$, see~\eqref{eq:nonlinear_equality1_sub_approx} and \eqref{eq:nonlinear_equality2_sub_approx}.
Therefore, without loss of generality, we make our assertions based on $h$ and ${\hat h}^{\hat{\vec{z}}}$ together with the assumption $\lim_{m\rightarrow\infty} \vec{z}^{(m)}=\vec{z}^\star$, where $\vec{z}^\star$ plays the role of $\vec{z}^\star_k$ and $\vec{z}^{(m)}$ plays the role of $\vec{z}^{(m)}_k$.

Let us next summarize some intermediate results, which will be useful later.
\begin{lemma}\label{prop:block_result}
Given the function $h$ definition of the form~\eqref{eq:block_function}, and $\lim_{m\rightarrow\infty} \vec{z}^{(m)}=\vec{z}^\star$, we have 1) $\lim_{m\rightarrow\infty}{\hat h}^{{\vec{z}}^{(m)}}({\vec{z}}^\star)= h({\vec{z}}^\star)$ and 2) $\lim_{m\rightarrow\infty} \nabla_{\vec{z}}{\hat h}^{{\vec{z}}^{(m)}}({\vec{z}}^\star)=\nabla_{\vec{z}} h({\vec{z}}^\star)$.
\end{lemma}
\begin{IEEEproof}
see Appendix~\ref{app:Proof_block_result}.
\end{IEEEproof}

From \emph{Lemma~\ref{prop:block_result}} above, we conclude that
\be\label{eq:KKT_sub_result}
h({\vec{z}}^\star)={\hat h}^{{\vec{z}}^\star}({\vec{z}}^\star) \qquad \mbox{and} \qquad \nabla_{\vec{z}} h({\vec{z}}^\star)=\nabla_{\vec{z}}{\hat h}^{{\vec{z}}^\star}({\vec{z}}^\star) \ .
\ee
By relating the result~\eqref{eq:KKT_sub_result} to our original problems~\eqref{eq:ADMM_private_variable_update_generator_bus} and~\eqref{eq:the_approximated_subproblem}, we can deduce that
\be\label{eq:objective_match}
{\boldsymbol{\lambda}}_k(\vec{z}_k^\star)=\hat{\boldsymbol{\lambda}}_k^{\vec{z}_k^\star}(\vec{z}_k^\star), \qquad {\boldsymbol{\mu}}_k(\vec{z}_k^\star)=\hat{\boldsymbol{\mu}}_k^{\vec{z}_k^\star}(\vec{z}_k^\star) \ , \ee
and
\be\label{eq:objective_derivative_match}
\hspace{-1mm}\bar{\nabla}_{\vec{z}_k}{\boldsymbol{\lambda}}_k(\vec{z}_k^\star){=}\bar{\nabla}_{\vec{z}_k}\hat{\boldsymbol{\lambda}}_k^{\vec{z}_k^\star}(\vec{z}_k^\star), \hspace{2mm} \bar{\nabla}_{\vec{z}_k}{\boldsymbol{\mu}}_k(\vec{z}_k^\star){=}{\bar{\nabla}}_{\vec{z}_k}\hat{\boldsymbol{\mu}}_k^{\vec{z}_k^\star}(\vec{z}_k^\star) \ ,
\ee
where $\bar{\nabla}$ is used to represent component-wise differentiation of associated functions.

Now we can easily conclude that $\vec{z}^\star_k$ is primal feasible for problem~\eqref{eq:ADMM_private_variable_update_generator_bus}. This follows from~\eqref{eq:objective_match}, the fact that constraints~\eqref{eq:local_variables_set_sub}, \eqref{eq:affine_sub}, \eqref{eq:convex_inequality_sub}, and \eqref{eq:convex_nonlinear_inequality_sub} are identical to~\eqref{eq:local_variables_set_sub_arrox}, \eqref{eq:affine_sub_approx}, \eqref{eq:convex_inequality_sub_approx}, and \eqref{eq:convex_nonlinear_inequality_sub_approx}, respectively, and that $\check{\mathcal{X}}^k_r\subseteq\mathcal{X}^k_r$.

Dual feasibility of $\vec{u}^\star_k$ associated with constraints~\eqref{eq:convex_inequality_sub_approx} and \eqref{eq:convex_nonlinear_inequality_sub_approx} affirms the dual feasibility of $\vec{u}^\star_k$ associated with identical constraints~\eqref{eq:convex_inequality_sub} and \eqref{eq:convex_nonlinear_inequality_sub}. In the case of constraint~\eqref{eq:nonconvex_inequality_sub_approx}, recall from~\eqref{eq:halfspace} that $\check{\mathcal{X}}^k_r$ is characterized by $\left( (\vec{v}\sr{Re}_k)_r,(\vec{v}\sr{Im}_k)_r  \right) \in \R^2$ such that
\be
c_r \leq a_r (\vec{v}_k\sr{Re})_r + b_r (\vec{v}_k\sr{Im})_r \hspace{2mm}\mbox{and} \hspace{2mm} (\vec{v}\sr{Re}_k)_r^2 {+} (\vec{v}\sr{Im}_k)_r^2 {\leq}   (\vec{v}\sr{max}_k)_r^2 \ .
\ee
Thus, dual feasibility of $\vec{u}^\star_k$ components associated with first (respectively second) constraint above ensures the dual feasibility of same $\vec{u}^\star_k$ components associated with $(\vec{v}\sr{min}_k)_r^2\leq  (\vec{v}\sr{Re}_k)_r^2 + (\vec{v}\sr{Im}_k)_r^2 $ (respectively $ (\vec{v}\sr{Re}_k)_r^2 + (\vec{v}\sr{Im}_k)_r^2 \leq   (\vec{v}\sr{max}_k)_r^2$) of~\eqref{eq:nonconvex_inequality_sub}. Thus, we conclude $\vec{u}^\star_k$ is dual feasible for problem~\eqref{eq:ADMM_private_variable_update_generator_bus}.

From~\eqref{eq:objective_match}, the fact that constraints~\eqref{eq:local_variables_set_sub}, \eqref{eq:affine_sub}, \eqref{eq:convex_inequality_sub}, and \eqref{eq:convex_nonlinear_inequality_sub} are identical to~\eqref{eq:local_variables_set_sub_arrox}, \eqref{eq:affine_sub_approx}, \eqref{eq:convex_inequality_sub_approx}, and \eqref{eq:convex_nonlinear_inequality_sub_approx}, respectively, and that the components $(\vec{v}\sr{Re}_k,\vec{v}\sr{Im}_k)$ of $\vec{z}^\star_k$, \emph{strictly} satisfy the constraint~\eqref{eq:nonconvex_inequality_sub_approx} [see \emph{Assumption~1}], it follows that $\vec{z}^\star_k$ and $\vec{u}^\star_k$ satisfy complementary slackness conditions for problem~~\eqref{eq:ADMM_private_variable_update_generator_bus}. In addition, \emph{Assumption~1} together with complementary slackness condition ensure that the components of $\vec{u}^\star_k$ associated with constraints~\eqref{eq:nonconvex_inequality_sub_approx} are \emph{identically~zero}.

Finally, recall that $(\vec{z}^\star_k,\vec{u}^\star_k)$ are optimal primal and dual variables for problem~\eqref{eq:the_approximated_subproblem}. Therefore, the derivative of the Lagrangian $\hat{L}_k(\vec{z}_k,\vec{u}_k)$ associated with problem~\eqref{eq:the_approximated_subproblem} vanishes at $(\vec{z}^\star_k,\vec{u}^\star_k)$,~i.e., $\nabla_{\vec{z}_k} \hat{L}_k(\vec{z}^\star_k,\vec{u}^\star_k) = 0$. 
This result combined with~\eqref{eq:objective_derivative_match}, the fact that constraints~\eqref{eq:local_variables_set_sub}, \eqref{eq:affine_sub}, \eqref{eq:convex_inequality_sub}, and \eqref{eq:convex_nonlinear_inequality_sub} are identical to~\eqref{eq:local_variables_set_sub_arrox}, \eqref{eq:affine_sub_approx}, \eqref{eq:convex_inequality_sub_approx}, and \eqref{eq:convex_nonlinear_inequality_sub_approx}, respectively, and the fact that the components of $\vec{u}^\star_k$ associated with constraints~\eqref{eq:nonconvex_inequality_sub_approx} are identically zero, affirms that derivative of the Lagrangian ${L}_k(\vec{z}_k,\vec{u}_k)$ associated with problem~\eqref{eq:ADMM_private_variable_update_generator_bus} vanishes at $\vec{z}^\star_k$ and ${\vec{u}}^\star_k$,~i.e.,
\be\label{eq:Lagrangian_sub}
\nabla_{\vec{z}_k} {L}_k(\vec{z}^\star_k,\vec{u}^\star_k) = 0 \ ,
\ee
which concludes the proof.
\end{IEEEproof}

\begin{table*}[t!]
\normalsize
\vspace{2mm}
\begin{align}\label{eq:derivative_Largrangian_whole}
       \nabla_{\vec{z},\vec{v}}  L(\vec{z}^\star,\vec{v}^\star, \vec{u}^\star, \vec{y}^\star)
                             \hspace{-1mm}=\hspace{-1mm}\left[\begin{array}{c}
                                    \hspace{-1mm}\nabla_{\vec{z}_1} L_1( \vec{z}^\star_1,\vec{u}^\star_1|\vec{y}_k^\star)-\rho \bar{\vec{z}}_1\\
                                    \vdots \\
                                    \nabla_{\vec{z}_N} L_N( \vec{z}^\star_N,\vec{u}^\star_N|\vec{y}_k^\star) - \rho \bar{\vec{z}}_N\\
                                    \nabla_{\vec{v}}\left(\displaystyle \ssum{k\in \mathcal{N}}   -\vec{y}_k\tran  \bec{E}_k\vec{v} \right)
                               \end{array}\right]
                               \hspace{-1mm}{=}\hspace{-1mm}\left[\begin{array}{c}
                                    -\rho\bar{\vec{z}}_1 \\
                                    \vdots \\
                                    -\rho\bar{\vec{z}}_N \\
                                    \vec{0}
                               \end{array}\right],
\end{align}
\be \label{eq:local_vecfd}
\bar{\vec{z}}_k=\big(0, 0, 0, 0, 0, 0,\underbrace{(\vec{x}^{\mbox{\tiny{Re}}\star}_k- \vec{E}_k\vec{v}^{\mbox{\tiny{Re}}\star}) , (\vec{x}^{\mbox{\tiny{Im}}\star}_k- \vec{E}_k\vec{v}^{\mbox{\tiny{Im}}\star}_k)}_{ \vec{x}^\star_k-\bar{\vec{E}}_k\vec{v}^\star=\boldsymbol{\delta}_k}, \vec{0}, \vec{0}, \vec{0}, \vec{0}\big) \ .
\ee


\begin{tabular}{p{17.7cm}} \hline  \\ \end{tabular}
\vspace{-8mm}
\end{table*}
\subsection{Proof of \emph{Proposition~2}} \label{app:proof_of_prop_2}
\begin{IEEEproof}
Given \emph{Assumption~1} holds, \emph{Proposition~\ref{lemma:subproblem_optimality}}, asserts that all constraints, but~\eqref{eq:dist_formulation_consistency_constraints} of problem~\eqref{eq:dist_formulation} are primal feasible. Combined this with~\eqref{eq:discrepancy_k}, it trivially follows that $\delta=a^{-1}\bar{\delta}$, where $a=\mbox{len}((\boldsymbol{\delta}_k)_{k\in\mathcal{N}})$ and $\bar{\delta}=\sum_{k\in\mathcal{N}}||{\boldsymbol \delta}_k||^2_2$~[\emph{cf}~\eqref{eq:KKT3}]. To show that $\epsilon=b^{-1}\rho^2\bar{\delta}$ [\emph{cf}~\eqref{eq:KKT6}], let us consider the Lagrangian ${L}(\vec{z},\vec{v}, \vec{u}, \vec{y})$ associated with problem~\eqref{eq:ADMM_private_variable_update_generator_bus}. Note that ${L}(\vec{z},\vec{v}, \vec{u}, \vec{y})$ is related to the Lagrangian ${L}_k(\vec{z}_k,\vec{u}_k|\vec{y}_k)$ of problem~\eqref{eq:ADMM_private_variable_update_generator_bus} as
 \begin{equation}
       \textstyle{L}(\vec{z},\vec{v}, \vec{u}, \vec{y}) =  \ssum{k\in \mathcal{N}} \Big( {L}_k(\vec{z}_k,\vec{u}_k|\vec{y}_k) - (\rho/2) || \vec{x}_k-\bec{E}_k\vec{v}||_2^2 \Big) \ . \nonumber
\end{equation}
Note the notation used when passing the parameters to ${L}_k$, where we have highlighted the dependence of ${L}_k$ on $\vec{y}_k$ [\emph{cf}~\eqref{eq:ADMM_step2_obj}].
Let us now inspect the derivative of the Lagrangian ${L}(\vec{z},\vec{v}, \vec{u}, \vec{y})$, evaluated at $(\vec{z}^\star,\vec{v}^\star, \vec{u}^\star, \vec{y}^\star)$. In particular, we have~\eqref{eq:derivative_Largrangian_whole} [see top of this page],
where $\bar{\vec{z}}_k$ is given in~\eqref{eq:local_vecfd}~[compare with~\eqref{eq:discrepancy_k}]. Here the first equality follows form standard derivation combined with \emph{Proposition~\ref{lemma:subproblem_optimality}}, the second equality follows from~\eqref{eq:Lagrangian_sub} and by invoking the optimality conditions for problem~\eqref{eq:ADMM_net_variable_update}, i.e., $\sum_{k\in \mathcal{N}}  \bec{E}_k\tran   \vec{y}_k^\star=0$. From~\eqref{eq:derivative_Largrangian_whole}-\eqref{eq:local_vecfd}, we conclude that $\epsilon=b^{-1}\rho^2\bar{\delta}$ [\emph{cf}~\eqref{eq:KKT6}], where $b=\mbox{len}(\vec{z}^\star,\vec{v}^\star)$. Finally, conditions~\eqref{eq:KKT1}, \eqref{eq:KKT2}, \eqref{eq:KKT4}, and \eqref{eq:KKT5}, associated with problem~\eqref{eq:dist_formulation} follow from straightforward arguments, which concludes the proof.
\end{IEEEproof}



\subsection{Proof of Lemma~\ref{prop:block_result}}\label{app:Proof_block_result}
Let $\vec{H}$ denote the Hessian of function $h$. Note that $\vec{H}$ is a matrix with constant entries and thus does not depend on $\vec{z}$. From the definition of the Taylor series expansion at ${\vec{z}}^{m}$ we have
\be\label{eq:taylor}
h({\vec{z}})-{\hat h}^{{\vec{z}}^{(m)}}({\vec{z}})
     =  (1/2)({\vec{z}}^{(m)}-{\vec{z}})\tran \vec{H}({\vec{z}}^{(m)}-{\vec{z}}) \ .
\ee
Moreover, differentiation of \eqref{eq:taylor} yields
\be\label{eq:nabla_taylor}
\nabla h({\vec{z}})-\nabla{\hat h}^{{\vec{z}}^{(m)}}({\vec{z}})
     =  \vec{H}({\vec{z}}^{(m)}-{\vec{z}}) \ .
\ee
To show the case~1 of the proposition, we consider the following relations:
 \begin{multline}
     h({\vec{z}}^\star)-{\hat h}^{{\vec{z}}^{(m)}}({\vec{z}}^\star) 
     =  (1/2)({\vec{z}}^{(m)}-{\vec{z}}^\star)\tran \vec{H}({\vec{z}}^{(m)}-{\vec{z}}^\star), \label{eq:taylor1}
 \end{multline}
 \begin{multline}
  (1/2)\lambda_{\textrm{min}}(\vec{H})||{\vec{z}}^{(m)}-{\vec{z}}^\star||^2_2 \leq h({\vec{z}}^\star)-{\hat h}^{{\vec{z}}^{(m)}}({\vec{z}}^\star) \\
  \leq (1/2)\lambda_{\textrm{max}}(\vec{H})||{\vec{z}}^{(m)}-{\vec{z}}^\star||^2_2 \label{eq:taylor2}\ ,
 \end{multline}
where \eqref{eq:taylor1} follows from~\eqref{eq:taylor} with $\vec{z}={\vec{z}}^\star$ and \eqref{eq:taylor2} follows from \eqref{eq:taylor1} and basics of linear algebra. By letting $m\rightarrow\infty$ in \eqref{eq:taylor2}, we conclude $\lim_{m\rightarrow\infty} {\hat h}^{{\vec{z}}^{(m)}}({\vec{z}}^\star)=h({\vec{z}}^\star)$, since $\lim_{m\rightarrow\infty} \vec{z}^{(m)}_k=\vec{z}^\star_k$. Similarly, by using \eqref{eq:nabla_taylor} and that $\vec{H}({\vec{z}}^{(m)}-{\vec{z}}^\star)\leq\lambda_{\textrm{max}}(\vec{H}\tran\vec{H})||{\vec{z}}^{(m)}-{\vec{z}}^\star||^2_2$ and $\vec{H}({\vec{z}}^{(m)}-{\vec{z}}^\star)\geq\lambda_{\textrm{min}}(\vec{H}\tran\vec{H})||{\vec{z}}^{(m)}-{\vec{z}}^\star||^2_2$, we conclude $\lim_{m\rightarrow\infty}\nabla{\hat h}^{{\vec{z}}^{(m)}}({\vec{z}}^\star)= \nabla h({\vec{z}}^\star)$.

\bibliographystyle{IEEEbib}
\bibliography{references_Sindri,jour_short,conf_short,references,referencesElisabetta_04_22,references_PPO}


\end{document}